\documentclass[10pt,A4]{amsart}

\usepackage{amsmath}
\usepackage{amsthm}
\usepackage{amssymb}
\usepackage{amsxtra}
\usepackage{amsfonts}
\usepackage{mathtools}
\usepackage{bm}
\usepackage{bbold}
\usepackage{braket}
\usepackage{longtable}
\usepackage{multirow}
\usepackage{graphicx}
\usepackage{stackrel}
\usepackage{braket}
\usepackage{enumerate}
\usepackage{svg}
\usepackage[noadjust]{cite}

\usepackage{xcolor}
\definecolor{webgreen}{rgb}{0,.5,0}
\definecolor{webbrown}{rgb}{.6,0,0}
\definecolor{RoyalBlue}{cmyk}{1, 0.50, 0, 0}
\usepackage[colorlinks=true, breaklinks=true, urlcolor=webbrown, linkcolor=RoyalBlue, citecolor=webgreen,backref=page]{hyperref}

\DeclareMathOperator{\Tr}{Tr}

\makeatletter

\usepackage{cite}
\usepackage{array}
\usepackage{color}
\usepackage{t1enc}
\usepackage{alltt}
\usepackage{epsfig}
\usepackage{mathrsfs}
\usepackage{graphicx,ifthen}
\usepackage[scaled]{helvet}
\usepackage{booktabs}
\usepackage{wasysym}
\usepackage{arydshln}
\usepackage{enumerate}
\usepackage{orcidlink}

\theoremstyle{plain}

\newtheorem{corollary}{Corollary}[section]
\newtheorem{lemma}{Lemma}[section]
\newtheorem{proposition}{Proposition}[section]
\theoremstyle{definition}
\newtheorem{definition}{Definition}[section]

\theoremstyle{remark}
\newtheorem{remark}{Remark}[section]

\makeatletter

\allowdisplaybreaks[1]

\newcommand{\C}{\mathbb C}
\newcommand{\R}{\mathbb R}

\newcommand{\N}{\mathbb N}

\newcommand{\Rmnum}[1]{\expandafter\@slowromancap\romannumeral #1@}

\DeclareMathSymbol{\Iota}{\mathalpha}{operators}{"49}
\DeclareMathAlphabet{\mathpzc}{OT1}{pzc}{m}{it}

\begin{document}
\title{Gap probabilities for the Bures-Hall Ensemble and the Cauchy-Laguerre Two-Matrix Model}
\author{N. S. Witte\,\orcidlink{0000-0001-7537-8444}${}^{1,2}$ and L. Wei\orcidlink{0000-0002-8065-0956}${}^{1}$}
\address[1]{Department of Computer Science \\Texas Tech University, Lubbock TX 79409-3104, USA}
\address[2]{School of Mathematics and Statistics, Victoria University of Wellington, PO Box 600 Wellington 6140, New Zealand}
\email{\tt n.s.witte@protonmail.com, luwei@ttu.edu}

\date{\vspace{-2ex}\today}

\begin{abstract}
The Bures metric and the associated Bures-Hall measure is arguably the best choice for studying the spectrum of the quantum mechanical density matrix with no apriori knowledge of the system. 
We investigate the probability of a gap in the spectrum of this model, either at the bottom $ [0,s) $ or at the top $ (s,1] $, 
utilising the connection of this Pfaffian point-process with the allied problem in the determinantal point-process of the two-dimensional Cauchy-Laguerre bi-orthogonal polynomial system,
now deformed with two variables $s,t$.
To this end we develop new general results about Cauchy bi-orthogonal polynomial system for a more general class of weights than the Laguerre densities: 
in particular a new Christoffel-Darboux formula, reproducing kernels and differential equations for the polynomials and their associated functions.
This system is most simply expressed as rank-3 matrix variables and possesses an associated cubic bilinear form.     
Furthermore under specialisation to truncated Laguerre type densities for the weight, of direct relevance to the Cauchy-Laguerre system, 
we construct a closed system of constrained, nonlinear differential equations in two deformation variables $s,t$,
and observe that the recurrence, spectral and deformation derivative structures form a compatible and integrable triplet of Lax equations.
\end{abstract}

\maketitle

\section{Volume Measures of Quantum States}
\setcounter{equation}{0}

We address an outstanding problem of the volume measures for bi-partite, mixed quantum systems with system and environment dimensions $ m,n $, $ n>m $ respectively
and having no additional symmetries.
In 1998 Michael Hall \cite{Hal_1998} posed the question {\it "What statistical ensemble corresponds to minimal prior knowledge about a quantum state?"}.
That a statistical question was framed here is recognition of the fact that the dimensionality of the Hilbert spaces grows exponentially with the number of qubits, or whatever the fundamental units are, 
and consequently the complexity of computing the relative volumes of different class of quantum states. 
Furthermore in \cite{HLW_2006} it is observed {\it "There is no question that random entangled states are far easier to understand than all entangled states."}, 
where they studied mixed-state measures induced by taking the partial trace over a larger system.
From a somewhat different perspective Aubrun posed the related question {\it "Is a random state entangled"} in \cite{Aub_2014}.

In addition the question of a typical state has independent interest for other reasons such as those studies addressing violations of the {B}ell inequalities for random pure states \cite{AZ_2015},
and how often is a random quantum state {$k$}-entangled? \cite{SWZ_2011}. 
And the answer to this question has many significance for applications such as the construction of quantum circuits \cite{NDGD_2006}.
Bounds on the entanglement within a qudit subsystem of a larger pure state and with averages formed from sampling taken over the larger one were studied in \cite{KZM_2002}.

Our object of interest is the quantum mechanical density matrix $ \rho $ with the properties: 
(i) a complex Hermitian $m \times m$ matrix $ \rho  = \rho^{\dagger} $, 
(ii) of unit trace $ {\rm Tr}\rho = 1 $
and (iii) being positive definite $ \rho > 0 $, i.e. $ \forall \left| \psi \rangle \right.\in \C^{m} $, $ \left.\langle \psi \right| \rho \left| \psi \rangle > 0 \right. $.
In particular we will focus on the spectrum of the density matrix $\{\rho_{j}\}^{m}_{j=1}$, $ 0\leq \rho_{j} \leq 1 $, $ \sum_{j=1}^{m} \rho_{j} = 1 $
and our statistical question will be framed with respect to these eigenvalues.
For any proper density matrix of dimension $m$ of a bipartite quantum system, its spectrum is supported in the probability simplex
$ \mathcal{D}=\left\{0\leq\rho_{m}<\ldots<\rho_{1}\leq 1,~~\sum_{i=1}^{m}\rho_{i}=1\right\} $,
where the special cases $ \rho_{1}=1,~~\rho_{2}=0,\dots,\rho_{m}=0 $ and $ \rho_{1}=\rho_{2}=\dots\rho_{m}=1/m $
are referred to as the separable state and the maximally-mixed state, respectively. 

One of the candidates in answer to Hall's question arises from the Hilbert-Schmidt (H-S) distance 
\begin{equation*}
	d^2_{\rm HS}(\rho_A,\rho_B) = {\rm Tr}((\rho_A-\rho_B)^2) ,
\end{equation*} 
and which has the joint probability density function for the eigenvalues
\begin{equation*}
	\mathcal{P}(\rho_{1}, \ldots, \rho_{m}) = \frac{1}{C_{m}} \mathbb{1}_{\sum_{j=1}^{m}\rho_{j}=1}
	\prod_{1\leq j<k \leq m}\left( \rho_{k}-\rho_{j} \right)^2 .
\end{equation*}
The Hilbert-Schmidt metric leads directly to the fixed trace $\beta=2$ Laguerre Ensemble \cite{LZ_2010}.

It has become clear, following the work of Fubini \cite{Fub_1904}, Study \cite{Stu_1905}, Bures \cite{Bur_1969}, 
Uhlmann \cite{Uhl_1976},\cite{Uhl_1986},\cite{Uhl_1992} and collaborators \cite{Hueb_1993},\cite{Hueb_1993a},\cite{Dit_1993} and Hall \cite{Hal_1998}
that the answer to Hall's question is the Bures-Hall measure.
One can find syntheses of these developments from a random matrix perspective in \cite{ZS_2001},\cite{SZ_2003},\cite{ZS_2003},\cite{SZ_2004}.
The Fubini-Study metric on projective Hilbert space $ {\rm C}{\mathbb P}^{m} $ starts with homogenous co-ordinates for $\C^{m+1}$, 
$ {\bf Z} = [Z_{0},Z_{1},\dots,Z_{m}] \in \C^{m+1}\backslash\{0\} $,
and the standard Hermitian metric on $\C^{m+1}$, 
$ ds^2 = dZ_{0} \otimes d\bar{Z}_{0} + \ldots + dZ_{m} \otimes d\bar{Z}_{m} $,
and transforming to affine co-ordinates for $ {\rm C}{\mathbb P}^{m} $, $ z_{j} = \frac{Z_{j}}{Z_{0}}, j = 1, \ldots, m $ in co-ordinate patch $ U_{0}=\{ Z_{0}\neq 0 \} $,
one deduces the squared distance
$ ds^2 = \sum_{j,k=1}^{m} g_{j,k} dz_{j} \otimes d\bar{z}_{k} $
with the metric tensor
\begin{equation*}
	g_{j,k} = \frac{1}{1+\left| {\bf z} \right|^2}\delta_{j,k} - \frac{1}{\left( 1+\left| {\bf z} \right|^2 \right)^2}\bar{z}_{j}z_{k} .
\end{equation*}
The Bures-distance \cite{Bur_1969} is given by
\begin{equation*}
	d^2_{\rm B}(\rho_A,\rho_B) = 2-2{\rm Tr}\sqrt{(\sqrt{\rho_A}\,\rho_B\sqrt{\rho_A})}
\end{equation*}
whose joint probability density function for the eigenvalues was found by H{\"u}bner \cite{Hueb_1993}, Hall \cite{Hal_1998}
\begin{equation}
	\mathcal{P}(\rho_{1}, \ldots, \rho_{m}) = \frac{1}{C_{m}} \mathbb{1}_{\sum_{j=1}^{m}\rho_{j}=1}\prod^{m}_{j=1}\rho_{j}^{-1/2} \times 
	\prod_{1\leq j<k \leq m}\frac{\left( \rho_{k}-\rho_{j} \right)^2}{\left(\rho_{k}+\rho_{j} \right)}, \quad \rho_{j} \in [0,1] .
\label{BuresHallJPDF}
\end{equation}
Here the normalisation is
\begin{equation*}
	 C_{m} = \frac{2^{-m(m-1)}\pi^{m/2}}{\Gamma(m^2/2)}\prod_{i=1}^{m}\Gamma(i+1)
\end{equation*}
The Bures-Hall (B-H) measure distinguishes itself from others in satisfying all four requirements:
being Riemannian \cite{Dit_1993}; 
monotone \cite{PS_1996}; 
Fubini-Study adjusted \cite{Uhl_1976} and 
Fisher adjusted \cite{PS_1996}.
A general background to the geometry of quantum states and the metrics relevant here can be found in \cite{Bengtsson+Zyczkowski_2017},
while a review of random matrix techniques in quantum information theory is given in \cite{CN_2016}.
The Bures-Hall measure is a structureless measure in that it only depends on the total dimension of the Hilbert space and not on the tensor product structure.
Some physically motivated ensembles of structured random states have been studied in \cite{ZPNC_2011} -
see also \cite{OSZ_2010} for a matrix-model realisation of the Bures-Hall measure.
	
The application of ideas from convex geometry, 
where the relative volumes of separable or entangled states are often expressed as bounds, 
can be found in the works \cite{Sza_2005},\cite{Ye_2009},\cite{Ye_2010},\cite{ASY_2014} and this approach features prominently in the recent monograph \cite{Aubrun+Szarek_2017}.
The quantum systems are Hilbert spaces of tensor products of $N $ qubit or qudit states, i.e. $ \C^{d_1} \otimes \cdots \otimes \C^{d_N} $ 
and the bounds are given in terms of the relevant dimensions with constant factors independent of these variables.

A number of statistics have been calculated for the Bures-Hall ensemble which are global in nature and only involve the one or two-point correlations such as
fidelity \cite{ZS_2005b},
quantum purity $ \mathbb{E}[{\rm Tr}(\rho^2)] $ \cite{ZS_2005b},\cite{OSZ_2010},\cite{BN_2011},\cite{LAK_2021} and
averages of the von Neumann entropy $ \mathbb{E}[-{\rm Tr}(\rho\log\rho)] $ \cite{Pag_1993} and its higher moments \cite{SK_2019},\cite{Wei_2020},\cite{Wei_2020a},\cite{LW_2021a}.	
In \cite{CLZ_2010} the smallest eigenvalue distribution of the fixed-trace Laguerre beta-ensemble, i.e. the Hilbert-Schmidt metric, was investigated.
The distribution of bipartite entanglement for random pure states was calculated by Giraud in \cite{Gir_2007} using the H-S metric also.
The large deviations of the von Neumann and Renyi entropies for a random bipartite state, i.e. the H-S metric,
was explored in \cite{NMV_2011} in the large dimension regime using Coulomb fluid methods. 

In contrast we study a generalised gap-probability here, which is a refined local statistic of the density matrix spectrum, in particular of the extreme eigenvalues. 
One specialisation of this generalised gap probability (see \eqref{BHFT_gapGF})\footnote{This generalises \eqref{BuresHallJPDF} from $a=-1/2$ to $a=n-m-1/2$.}
is the multivariate integral over the Bures-Hall fixed trace ensemble (B-HFT)
\begin{equation}
	Z^\text{B-HFT}(t;m,a) := \frac{1}{C^\text{B-HFT}(m,a)}
	\frac{1}{m!}\int_{0}^{t}d\rho_{1} \ldots \int_{0}^{t}d\rho_{m}\; \delta\Big( \sum_{j=1}^{m}\rho_{j}-1 \Big)
	\prod_{1\leq j<k \leq m}\frac{(\rho_{k}-\rho_{j})^{2}}{(\rho_{k}+\rho_{j})} \prod_{j=1}^{m} \rho_{j}^{a} ,
\label{FTBH_gap}
\end{equation}
with a cut-off value $0<t$.
To see how our statistic relates to the problem formulated above consider the following: 
this statistic is formed from the integral over that part of the simplex in the positive orthant intersecting with the hypercube $ [0,t]^m $.
Initially $Z^\text{B-HFT}$ remains at a plateau of zero as $t$ increases from zero until the first threshold is reached at $t=\frac{1}{m}$ 
(the body-diagonal of the cube passes through the simplex), and the first piece-wise contribution is made to the integral rising from zero, 
followed by a cascade of thresholds when lower-dimensional faces pass through the simplex and each makes additional piece-wise contributions at $ t=\frac{1}{m-1}, \frac{1}{m-2}, \ldots $.
Finally at $ t=1 $ the simplex is wholly contained within the cube and $ Z^\text{B-HFT}=1 $.
Thus the asymptotic behaviour of $ Z^\text{B-HFT}(t) $ for $ m\to \infty $ and $ t\to 0^{+} $ is recording the contributions made from maximally mixed states.    

The primary object studied with our methods will be the Cauchy-Laguerre two matrix model (C-L2M) and its associated bi-orthogonal polynomial system 
because of a key identity which connects this to the unit-trace Bures-Hall ensemble via an unconstrained trace version of the latter.
The Cauchy-Laguerre two matrix model has a joint eigenvalue PDF in two sets (or species) of eigenvalues $ x_{j}\in \R_{+}, j=1,\ldots m $ and $y_{j}\in \R_{+}, j=1,\ldots m $
\begin{multline}
	\mathcal{P}^\text{C-L2M}(x_{1},\ldots,x_{m};y_{1},\ldots,y_{m})
		= \frac{1}{C^\text{C-L2M}(m,a,b)(m!)^2}
		  \prod_{j=1}^{m} x_{j}^{a}e^{-x_{j}} y_{j}^{b}e^{-y_{j}}
\\	  \times
		  \frac{\prod_{1\leq j<k \leq m}(x_{k}-x_{j})(y_{k}-y_{j})}{\prod_{1\leq j,k \leq m}(x_{j}+y_{k})}
		  \prod_{1\leq j<k \leq m}(x_{k}-x_{j})(y_{k}-y_{j}) .
\label{CL2M_PDF}
\end{multline}
and the gap probability corresponding to \eqref{FTBH_gap} is the $ t=s $, $ b=a+1$ specialisation (see \eqref{CL2M_gapGF} for the full generating function) of
\begin{equation}
	Z^\text{C-L2M}(s,t;m,a,b):= 
	\int_{0}^{s} dx_{1} \ldots \int_{0}^{s} dx_{m}\; \int_{0}^{t} dy_{1} \ldots \int_{0}^{t} dy_{m}\; \mathcal{P}^\text{C-L2M}(x_{1},\ldots,x_{m};y_{1},\ldots,y_{m})  ,
\label{CL2M_gap}
\end{equation}
where $ s,t\in \R_{+} $ are the gap cutoffs. When $ s=t=\infty $ the statistic reduces to the undeformed case.

The undeformed Cauchy-Laguerre bi-orthogonal system, i.e. when the support is $\R^2_{+}$, was initially formulated as a matrix model in the early literature \cite{BGS_2009}, 
and subsequently the bi-orthogonal polynomials and the four kernels were shown to have explicit evaluations as Meijer-G functions \cite{BGS_2014}.
Intermediate between these two works the Cauchy bi-orthogonal system in a general setting was investigated \cite{BGS_2010} and a number of foundational results were established.
This latter work will be our starting point. 
Our goal will be to derive a closed system of partial differential equations with respect to $s,t$ which uniquely characterise this statistic. 

Previous work on a specialisation of the above gap probability, that of $[\xi^m,\psi^m]Z^\text{C-L2M}$, 
only reported numerical evaluations \cite{BGS_2014} of $Z^\text{C-L2M}$ versus $s,t$,
where this was recast as a Fredholm determinant defined by integral operators with kernels constructed from Meijer-G functions,
and computed using the algorithm of Bornemann \cite{Bor_2010b}.
In \cite{LiLi_2019} the authors treated a different deformation of the Cauchy two matrix model where the univariate weights $ w_{1}(x), w_{2}(y) $ 
(both equal, as they only considered the symmetric case)
acquired an extra factor $ e^{xt} $ with a deformation variable $t$ with the support remaining as $\R^2_{+}$.
They identified the C-Toda lattice and CKP hierarchy as the integrable system characterising the dynamics but made no connection with isomonodromic deformations.
	
The authors of \cite{BGS_2013} applied a nonlinear steepest descent method to a class of 
$3\times 3$ Riemann-Hilbert problems introduced in connection with the Cauchy two-matrix random model. 
However this work is not directly relevant in the deformed case as the support for the bi-orthogonal system is $\R^2_{+}$, 
even though the two equilibrium measures were supported on an arbitrary number of intervals.
They solved the Riemann-Hilbert problem for the outer parametrix in terms of sections of a spinorial line bundle on a three-sheeted Riemann surface of arbitrary genus and established strong asymptotic results for the Cauchy biorthogonal polynomials. 

A further generalisation in another direction to the one we consider here, but parallel to that of the Muttalib-Borodin ensembles, was pursued in \cite{FL_2019},
where the second Vandermonde factor in \eqref{CL2M_PDF} is replaced by a $\theta$ extension
$ \prod_{1\leq j<k \leq m}(x_{k}^{\theta}-x_{j}^{\theta})(y_{k}^{\theta}-y_{j}^{\theta}) $
for $ {\rm Re}(\theta)>0 $. Again the support is $\R^2_{+}$, and now the evaluations of the bi-orthogonal polynomials and kernels generalise to Fox H-functions.

Our key results are expressed in the form of a closed set of fundamental coupled, first-order non-linear differential equations, 
given in Prop.\ref{PQ_partial-derivatives}, \ref{PQ_total-derivatives}, \ref{PiEta_derivatives} and \ref{XY_derivatives},
and subject to constraints given in Prop.\ref{constraints}.
This set is fundamental in the sense that any other derivative can be computed from members of this set.
All of \S\ref{Constrained_Dynamics} can be considered as our summary.
We consider our results as a preliminary step towards understanding the distribution of the spectral edges of the density matrix and the dynamical system characterising this,
and have deferred a number of logically compelling tasks as these are non-trivial exercises in themselves, and beyond the scope of the current work.
Since the results of our work may not be adequately summarised in cursory way we present an overall plan of the logic behind our central theme in the following algorithmic form: 
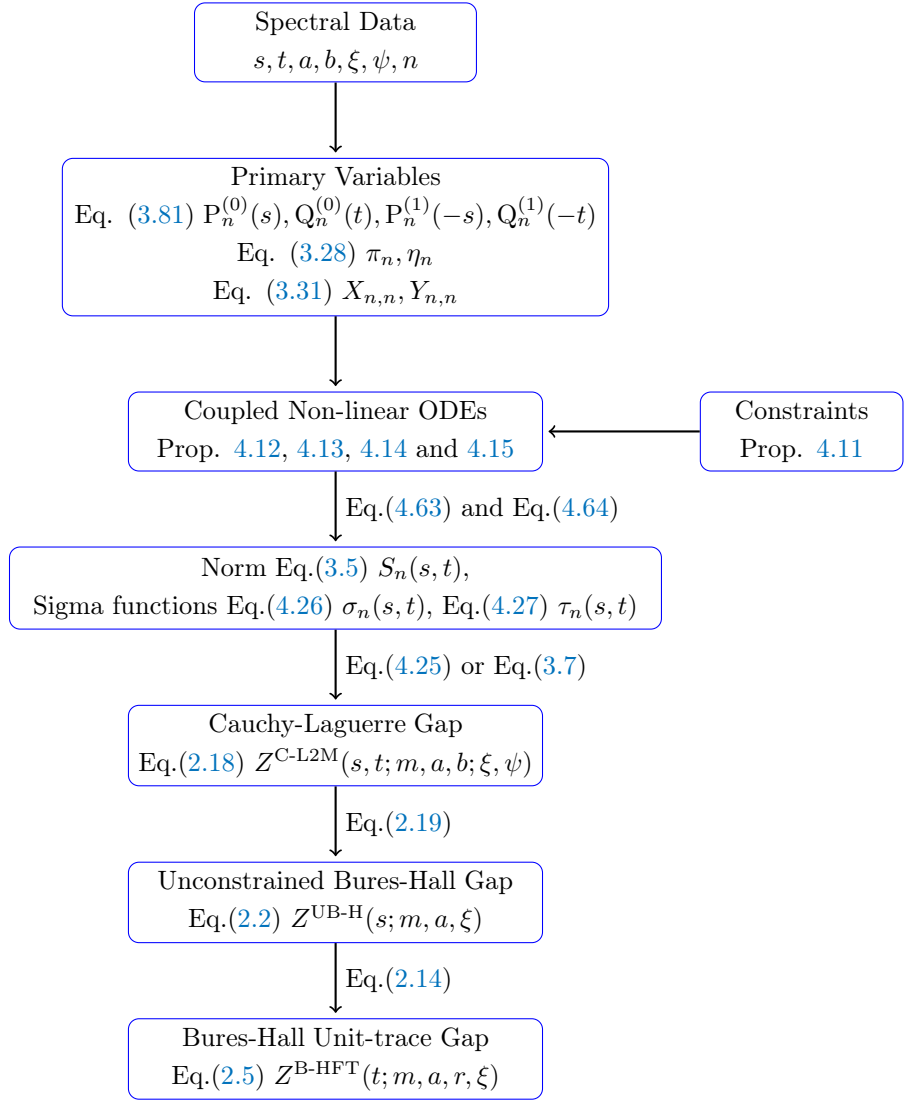
\begin{center}
{
\begin{tikzpicture} [
    auto,
    decision/.style = { diamond, draw=blue, thin, fill=white!20,
                        text width=10em, text badly centered,
                        inner sep=1pt, rounded corners },
    narrowblock/.style    = { rectangle, draw=blue, thin, 
                        fill=white!20, text width=7em, text centered,
                        rounded corners, minimum height=2em },
    block/.style    = { rectangle, draw=blue, thin, 
                        fill=white!20, text width=10em, text centered,
                        rounded corners, minimum height=2em },
    medblock/.style    = { rectangle, draw=blue, thin, 
                        fill=white!20, text width=15em, text centered,
                        rounded corners, minimum height=2em },
    wideblock/.style    = { rectangle, draw=blue, thin, 
                        fill=white!20, text width=20em, text centered,
                        rounded corners, minimum height=2em },
    extrawideblock/.style    = { rectangle, draw=blue, thin, 
                        fill=white!20, text width=24em, text centered,
                        rounded corners, minimum height=2em },
    line/.style     = { draw, thick, ->, shorten >=2pt },
  ]
  \matrix [column sep=5mm, row sep=10mm] {
                    & \node [block] (Data) {Spectral Data \\ $s,t,a,b,\xi,\psi,n$};            & \\
                    & \node [wideblock] (PQ+pi-eta+XY) {	Primary Variables \\
                    							Eq. \eqref{3x3_PQ} $\mathrm{P}^{(0)}_{n}(s), \mathrm{Q}^{(0)}_{n}(t), \mathrm{P}^{(1)}_{n}(-s), \mathrm{Q}^{(1)}_{n}(-t)$	\\ 
                    							Eq. \eqref{pi-eta} $\pi_{n}, \eta_{n}$  \\
                    							Eq. \eqref{XY-elements} $X_{n,n}, Y_{n,n}$};   & \\
				 	& \node [medblock] (ODE) {Coupled Non-linear ODEs\\
    											Prop. \ref{PQ_partial-derivatives}, \ref{PQ_total-derivatives}, \ref{PiEta_derivatives} and \ref{XY_derivatives}};   
    				& \node[narrowblock](Constraints){Constraints Prop. \ref{constraints}}; \\
                    & \node [extrawideblock] (SigmaForm)		{Norm Eq.\eqref{BOPS_norm} $S_{n}(s,t)$, \\ 
                    										Sigma functions Eq.\eqref{sigma} $\sigma_{n}(s,t)$, Eq.\eqref{tau} $\tau_{n}(s,t)$}; & \\
                    & \node [medblock] (Cauchy-Laguerre) 	{Cauchy-Laguerre Gap \\ Eq.\eqref{CL2M_gapGF} $Z^\text{C-L2M}(s,t;m,a,b;\xi,\psi)$};            & \\
                    & \node [medblock] (U-Bures-Hall) 	{Unconstrained Bures-Hall Gap \\ Eq.\eqref{UBH_defn} $Z^\text{UB-H}(s;m,a,\xi)$};  & \\
                    & \node [medblock] (Bures-Hall) 		{Bures-Hall Unit-trace Gap \\ Eq.\eqref{BHFT_gapGF} $Z^\text{B-HFT}(t;m,a,r,\xi)$};   & \\
  };
  \begin{scope} [every path/.style=line]
    \path (Data)        --    (PQ+pi-eta+XY);
    \path (PQ+pi-eta+XY)      --    node [right] {} (ODE);
    \path (Constraints)      --    (ODE);
    \path (ODE)   --    node [right] {Eq.\eqref{K01_Limit@s} and Eq.\eqref{K10_Limit@-t}} (SigmaForm);
    \path (SigmaForm)    --    node [right] {Eq.\eqref{sigma_Defn} or Eq.\eqref{norm+leading-coeff}} (Cauchy-Laguerre);
    \path (Cauchy-Laguerre)    --    node [right] {Eq.\eqref{UBHgf_CL2MM_gf}} (U-Bures-Hall);
    \path (U-Bures-Hall) --    node [right] {Eq.\eqref{FT-U_connect}} (Bures-Hall);
  \end{scope}
\end{tikzpicture}
}
\end{center}
 
The sequential layout of our manuscript is as follows:
\begin{description}
\item[\S \ref{BHFT_to_UCL}]
	We begin by defining the generalised gap probability for the unit-trace Bures-Hall ensemble and employ a Laplace transform relating it to an unconstrained ensemble. 
	In a second step we relate the Pfaffian point process of this unconstrained Bures-Hall ensemble to the determinantal point process of the Cauchy-Laguerre two matrix model via the identity of Forrester-Kieburg \cite{FK_2016}.
	Then we construct the explicit Bi-moment matrices for the Pfaffian and determinant formulae of the generalised gap probabilities for the unconstrained Bures-Hall ensemble 
	and the Cauchy-Laguerre two matrix model.   
\item[\S \ref{Cauchy-BOPS}] 
	Our first part treats the Cauchy Bi-orthogonal system of polynomials $\{P_{n}(x), Q_{n}(y)\}_{n=0}^{\infty}$ for general classes of weights, building on the foundations of \cite{BGS_2010}.
	We recall basic definitions and results such as the rank-one moment relation, the two recurrence relations, the multiplication operators,
	the Stieltjes and associated functions corresponding to the bi-orthogonal polynomials forming a rank-3 matrix system and their transfer matrices.
	We proceed by defining reproducing and related kernels and derive a new and simple Christoffel-Darboux formulae for these.
	From this point we construct spectral differential equations for the polynomials and associated functions, thereby defining spectral Lax matrices.
	Linking relations between the transfer matrices and the Lax matrices of the $P$ and $Q$ systems are studied and their invariants found.
	Using these two novel ingredients we can explicitly connect the Lax matrices for the $P$ and $Q$ sub-systems which hitherto have been parallel yet uncoupled structures.
\item[\S \ref{Cauchy-Laguerre-BOPS}] 
	In the second part we specialise our univariate weight factors to Laguerre densities and to the Cauchy-Laguerre bi-orthogonal system of polynomials,
	which introduces two deformation variables $s,t$.
	In this setting we derive the spectral $x,y$ and deformation $s,t$ derivatives of the polynomials, associated functions and coefficients from first principles,
	thus giving the spectral Lax matrices and the two deformation Lax matrices in explicit rational $x,y$ form for the $P,Q$ system.
	All the Lax matrices are given in terms of polynomials and the first associated functions evaluated at the finite, regular singular points $x=\pm s$, $y=\pm t$.
	In addition we derive a number of relations linking auxiliary coefficients, recurrence relation coefficients and evaluated polynomials and associated functions.  
	We write down all the invariants of the Lax matrices and the closed system of constrained dynamics as coupled, first order nonlinear differential equations.
	In conclusion we verify compatibility of these structures and the integrable nature of the system.
\end{description}

\section{Definitions and Conventions for the Bures-Hall Fixed Trace Ensemble and Unconstrained Ensemble}\label{BHFT_to_UCL}
\setcounter{equation}{0}

\subsection{Bures-Hall Fixed Trace Ensemble}
The joint probability density function (JPDF) of the density matrix eigenvalues or Schmidt values $ \rho_{1},\ldots,\rho_{m} \in \R_{+} $ of the Bures-Hall fixed trace ensemble (B-HFT)
$m\geq 1$, $\Re(a)>-1$, and trace $r>0$ is
\begin{equation}\label{BHFT_jpdf}
	\mathcal{P}^\text{B-HFT}(\rho_{1},\ldots,\rho_{m})
	 = \frac{1}{C^\text{B-HFT}(m,a)m!} \mathbb{1}_{\sum_{j=1}^{m}\rho_{j}-r}
	 \prod_{1\leq j<k \leq m}\frac{(\rho_{k}-\rho_{j})^{2}}{(\rho_{k}+\rho_{j})}
	 \prod_{j=1}^{m} \rho_{j}^{a} .
\end{equation}
The exponent introduced here is related to the model parameters by $a=n-m-\tfrac{1}{2}$.

There are two points to note here. Firstly the ``fixed'' trace will not be unity but taken as variable for reasons that will become clear later. 
Secondly that the support of $ \rho_{j} $ is often taken to be $\R_{+}:=[0,\infty)$ for convenience rather than the strict interval $ 0\leq \rho_{j} \leq r $,
where this is permissible for the existence of the relevant integral.
The normalisation of the B-HFT JPDF is defined as a variation of that given in \cite{SZ_2003}, see Eq. (3.11), $\Re(\alpha)>\frac{1}{2}$, $\Re(\beta)>0$, $a=\alpha-\frac{3}{2}$
\begin{equation}\label{BHFT_norm-defn}
	C^\text{B-HFT}(m,a,\beta)
	= \frac{1}{m!}\int_{0}^{\infty}d\rho_{1} \ldots \int_{0}^{\infty}d\rho_{m}\; \delta\Big( \sum_{j=1}^{m}\rho_{j}-1 \Big)
		\prod_{1\leq j<k \leq m}\frac{(\rho_{k}-\rho_{j})^{\beta}}{(\rho_{k}+\rho_{j})^{\beta/2}}
		\prod_{j=1}^{m} \rho_{j}^{\alpha-\tfrac{3}{2}} ,	 
\end{equation}
and its evaluation given in Eq. (5.14) of \cite{SZ_2003}
\begin{equation}\label{BHFT_norm}
	C^\text{B-HFT}(m,a,\beta)
	= \frac{\pi^{m/2}2^{-m(2\alpha-2+(m-1)\beta/2)}}{m!\Gamma\left(\tfrac{1}{2}m(2\alpha-1+(m-1)\beta/2)\right)}
	  \prod_{j=1}^m \frac{\Gamma(1+j\beta/2)\Gamma(2\alpha-1+(m-j)\beta/2)}{\Gamma(1+\beta/2)\Gamma(\alpha+(m-j)\beta/2)} .
\end{equation}
Consequently in our case the normalisation of \eqref{BHFT_jpdf} for $\beta=2$ is
\begin{equation}\label{FT_norm}
	C^\text{B-HFT}(m,a)
	=  	\frac{\pi^{m/2}2^{-m(2a+m)}}{m!\Gamma\left(\tfrac{1}{2}m(2a+m+1)\right)}
		\prod_{j=1}^m \frac{\Gamma(j+1)\Gamma(2a+j+1)}{\Gamma(a+j+\tfrac{1}{2})} .
\end{equation}
Having introduced some key concepts we are in a position for define a gap probability.
\begin{definition}\label{BHFT_defn}
The generalised gap probability for the B-HFT ensemble, or the generating function thereof, is defined to be 
\begin{multline}\label{BHFT_gapGF}
	Z^\text{B-HFT}(t;m,a,r,\xi)
\\
	 := \frac{1}{C^\text{B-HFT}(m,a)}
	\frac{1}{m!}\left( \int_{0}^{\infty} - \xi\int_{t}^{\infty} \right)d\rho_{1} \ldots \left( \int_{0}^{\infty} - \xi\int_{t}^{\infty} \right)d\rho_{m}\;
	\delta\Big( \sum_{j=1}^{m}\rho_{j}-r \Big)
	\prod_{1\leq j<k \leq m}\frac{(\rho_{k}-\rho_{j})^{2}}{(\rho_{k}+\rho_{j})} \prod_{j=1}^{m} \rho_{j}^{a} ,
\end{multline}
for the cut-off value is $0<t<\infty$ and the generating function variable is $\xi$ and is taken as an indeterminate, possibly complex.
We have written the upper terminals of the integrals as $\infty $ for subsequent convenience but in fact the support of the integrand is bounded due to the trace condition.
The definition \eqref{BHFT_gapGF} corresponds to the normalisation integral formula with the weight function $ w(\rho)=\rho^a $ replaced by the piecewise weight 
$ w(\rho)=(1-\xi \chi_{\rho>t})\rho^a $.
\end{definition}
\begin{remark}
From this definition we note the following.
\mbox{}\\
\noindent (i)
	$ Z^\text{B-HFT} $ is a polynomial of degree $m$ in $\xi$.
\mbox{}\\
\noindent (ii)
	When $\xi=0$ there is no dependence on $t$ and $ Z^\text{B-HFT} = 1 $.
\mbox{}\\
\noindent (iii)
	When $\xi=1$ and $r=1$ then $ Z^\text{B-HFT} $ is the probability that no eigenvalues are contained in the interval $[t,\infty)$
	\begin{equation}\label{BHFT_top-gap}
		Z^\text{B-HFT}(t;m,a,1,1) := \frac{1}{C^\text{B-HFT}(m,a)}
		\frac{1}{m!} \int_{0}^{t}d\rho_{1} \ldots \int_{0}^{t}d\rho_{m}\;
		\delta\Big( \sum_{j=1}^{m}\rho_{j}-1 \Big)
		\prod_{1\leq j<k \leq m}\frac{(\rho_{k}-\rho_{j})^{2}}{(\rho_{k}+\rho_{j})} \prod_{j=1}^{m} \rho_{j}^{a} .
	\end{equation}
	Furthermore as $0<t$ increases $ Z^\text{B-HFT} $ monotonically increases from zero to unity.
	This the reason for choosing the normalisation in \eqref{BHFT_gapGF}. 
\mbox{}\\
\noindent (iv)
	When $r=1$ then coefficient of the leading monomial in $\xi$ of $ Z^\text{B-HFT} $ is the probability that no eigenvalues are contained in the interval $[0,t)$
	\begin{equation}\label{BHFT_bottom-gap}
		[\xi^m]Z^\text{B-HFT}(t;m,a,1,\xi) := \frac{1}{C^\text{B-HFT}(m,a)}
		\frac{(-1)^{m}}{m!} \int_{t}^{\infty}d\rho_{1} \ldots \int_{t}^{\infty}d\rho_{m}\;
		\delta\Big( \sum_{j=1}^{m}\rho_{j}-1 \Big)
		\prod_{1\leq j<k \leq m}\frac{(\rho_{k}-\rho_{j})^{2}}{(\rho_{k}+\rho_{j})} \prod_{j=1}^{m} \rho_{j}^{a} .
	\end{equation}
	In this case as $0<t$ increases $ Z^\text{B-HFT} $ monotonically decreases from unity to zero. 
\mbox{}\\
\noindent (v)
	The coefficient of intermediate monomials in $\xi$ give the conditional probability that a fixed number of eigenvalues $k$ lie in $[0,t)$ and the remainder $m-k$ lie in $[t,\infty)$,
	up to a combinatorial factor of ${m}\choose{k}$.
\end{remark} 
	
\subsection{Unconstrained Bures-Hall Ensemble}\footnote{A better terminology would be "Pfaffian Cauchy-Laguerre".}

The joint probability density function (JPDF) of the  eigenvalues $x_{1},\ldots,x_{m} \in \R_{+}$ of the unconstrained Bures-Hall ensemble (UB-H) $m\geq 1$, $\Re(a)>-1$ is
\begin{equation}\label{UBH_jpdf}
	\mathcal{P}^\text{UB-H}(x_{1},\ldots,x_{m})
	= \frac{1}{C^\text{UB-H}(m,a)m!}
	\prod_{1\leq j<k \leq m}\frac{(x_{k}-x_{j})^{2}}{(x_{k}+x_{j})}
	\prod_{j=1}^{m} x_{j}^{a}e^{-x_{j}} .
\end{equation}
The normalisation for \eqref{UBH_jpdf} is defined by, see Eq. (3.1) of \cite{FK_2016},
\begin{equation}\label{UBH_norm}
	C^\text{UB-H}(m,a)
	= \frac{1}{m!} \int_{0}^{\infty}\;dx_{1} \ldots \int_{0}^{\infty}\;dx_{m} 
	\prod_{1\leq j<k \leq m}\frac{(x_{k}-x_{j})^{2}}{(x_{k}+x_{j})}
	\prod_{j=1}^{m} x_{j}^{a}e^{-x_{j}},
	\quad x_{1},\ldots,x_{m} \in \R_{+} ,
\end{equation}
and has the evaluation
\begin{equation}\label{U_norm}
	C^\text{UB-H}(m,a) = \pi^{m/2} m! 2^{-m(2a+m)} \prod_{j=1}^m \frac{\Gamma(j)\Gamma(2a+j+1)}{\Gamma(a+j+\tfrac{1}{2})} .
\end{equation}

The definition of a gap probability for the UB-H ensemble can be made along similar lines to that of B-HFT.
\begin{definition}\label{UBH_defn}
The generalised gap probability for the UB-H ensemble, or the generating function thereof, is defined to be 
\begin{equation}\label{UBH_gap}
	Z^\text{UB-H}(s;m,a,\xi)
	:= \frac{1}{C^\text{UB-H}(m,a)}
	\frac{1}{m!}\left( \int_{0}^{\infty} - \xi\int_{s}^{\infty} \right)dx_{1}	 \ldots \left( \int_{0}^{\infty} - \xi\int_{s}^{\infty} \right)dx_{m}\;
	\prod_{1\leq j<k \leq m}\frac{(x_{k}-x_{j})^{2}}{(x_{k}+x_{j})} \prod_{j=1}^{m} x_{j}^{a} e^{-x_{j}} ,
\end{equation}
with a cut-off value $0<s<\infty$ and the generating function variable is $\xi$ and is taken as an indeterminate, possibly complex.
The definition \eqref{UBH_gap} corresponds to the normalisation integral formula with the weight function $ w(x)=x^a e^{-x} $ replaced by the piecewise weight 
$ w(x)=(1-\xi \chi_{x>s})x^a e^{-x} $.
\end{definition}

From the above definition we can make essentially the same notes as we did for the B-HFT ensemble.
At this point we can relate the two generating functions (GF) defined in definitions \ref{BHFT_defn} and \ref{UBH_defn}.
We use the notational conventions
\begin{equation}\label{LaplaceXfm}
	F(s) := \mathcal{L}[f(r);s] := \int_{0}^{\infty}\;dr\, e^{-sr}f(r) , \quad \Re(s)>c ,
\end{equation}
for some $c>0$ and $ f(r)={\rm O}(r^{b}) $ as $r\to 0^{+}$ for $\Re(b)>-1$ to ensure the existence of the integral and its inverse
\begin{equation}\label{I-LaplaceXfm}
	f(r) = \frac{1}{2\pi i}\int_{c+i\infty}^{c-i\infty}\;ds\, e^{rs} F(s) =: \mathcal{L}^{-1}[F(s);r] .
\end{equation}
\begin{proposition}
Let $ {\rm Re}(s)>0 $.
The Bures-Hall fixed trace generating function $ Z^\text{\rm B-HFT}(t;m,a,r,\xi) $ 
is given by the inverse Laplace transform of the unconstrained Bures-Hall generating function $ Z^\text{\rm UB-H}(s;m,a,\xi) $ 
\begin{equation}\label{FT-U_connect}
	Z^\text{\rm B-HFT}(t;m,a,r,\xi) = \frac{C^\text{\rm UB-H}(m,a)}{C^\text{\rm B-HFT}(m,a)}\mathcal{L}^{-1}[s^{-m(m+2a+1)/2}Z^\text{\rm UB-H}(st;m,a,\xi);r] ,
\end{equation}
which serves as our principal formula for the fixed trace generating function upon setting $r=1$.
\end{proposition}
\begin{proof}
The Laplace transform of \eqref{BHFT_gapGF} with respect to the trace variable $r$ 
can be exchanged with the $ \rho_1, \ldots, \rho_m $ integrals due to the absolute and uniform convergence of any pair of integrals.
The only integration required here is of the Dirac measure. 
If $ {\rm Re}(s)>0 $ then we make a change of variables $ x_j=s\rho_j $ and the integrand factors into an algebraic factor $ s^{-m(m+2a+1)/2} $ and
the integrand of the unconstrained Bures-Hall ensemble. The non-trivial aspect is that the deformation variable gets mapped $ t \mapsto st $.
\end{proof}

\subsection{Cauchy-Laguerre Two-Matrix Ensemble}

Our second key linkage is to employ the relation of the unconstrained Bures-Hall ensemble, which is a Pfaffian point process, 
with a specialisation of the Cauchy-Laguerre two-matrix ensemble, a determinantal point process.
The first observation of this relation was made in \cite{BGS_2009} and systematically resolved in \cite{FK_2016}.
To begin with we recall some essential knowledge of this ensemble.

The JPDF of the two sets of eigenvalues $x_{1},\ldots,x_{m} \in \R_{+}$ and $y_{1},\ldots,y_{m} \in \R_{+}$ of the Cauchy-Laguerre two-matrix ensemble (C-L2M) $m\geq 1$, $\Re(a),\Re(b)>-1$ is
\begin{multline}\label{CL2M_jpdf}
	\mathcal{P}^\text{C-L2M}(x_{1},\ldots,x_{m};y_{1},\ldots,y_{m})
	= \frac{1}{C^\text{C-L2M}(m,a,b)(m!)^2}
	  \prod_{j=1}^{m} x_{j}^{a}e^{-x_{j}} y_{j}^{b}e^{-y_{j}}
\\	  \times
	  \frac{\prod_{1\leq j<k \leq m}(x_{k}-x_{j})(y_{k}-y_{j})}{\prod_{1\leq j,k \leq m}(x_{j}+y_{k})}
	  \prod_{1\leq j<k \leq m}(x_{k}-x_{j})(y_{k}-y_{j}) .
\end{multline}

The normalisation for \eqref{CL2M_jpdf} is defined by,
\begin{multline}\label{}
	C^\text{C-L2M}(m,a,b)
	= \frac{1}{(m!)^2} \int_{0}^{\infty}\;dx_{1} \ldots \int_{0}^{\infty}\;dx_{m} \int_{0}^{\infty}\;dy_{1} \ldots \int_{0}^{\infty}\;dy_{m} 
	  \prod_{j=1}^{m} x_{j}^{a}e^{-x_{j}} y_{j}^{b}e^{-y_{j}}
\\	  \times
	  \frac{\prod_{1\leq j<k \leq m}(x_{k}-x_{j})(y_{k}-y_{j})}{\prod_{1\leq j,k \leq m}(x_{j}+y_{k})}
	  \prod_{1\leq j<k \leq m}(x_{k}-x_{j})(y_{k}-y_{j}) .
\end{multline}
see \cite{BGS_2009}, \cite{BGS_2010}, and has the evaluation from Eq. (2.5) of \cite{FL_2019} 
\begin{equation}\label{CL2M_norm}
	C^\text{C-L2M}(m,a,b) = \left( \prod_{j=1}^{m-1}j! \right)^2 \prod_{k=0}^{m-1}\Gamma(a+1+k)\Gamma(b+1+k) \prod_{j=1}^{m} \frac{(a+b-1+j)!}{(m+a+b-1+j)!} .
\end{equation} 

The definition of a gap probability for the C-L2M ensemble will be made along similar lines to that of UB-H ensemble.
\begin{definition}\label{CL2M_defn}
The generalised gap probability for the C-L2M ensemble, or the generating function thereof, is defined to be 
\begin{multline}\label{CL2M_gapGF}
	Z^\text{C-L2M}(s,t;m,a,b;\xi,\psi)
	\\
	:= \frac{1}{C^\text{C-L2M}(m,a,b)}\frac{1}{(m!)^2}
	\left( \int_{0}^{\infty} - \xi\int_{s}^{\infty} \right)dx_{1} \ldots \left( \int_{0}^{\infty} - \xi\int_{s}^{\infty} \right)dx_{m}\;
	\left( \int_{0}^{\infty} - \psi\int_{t}^{\infty} \right)dy_{1} \ldots \left( \int_{0}^{\infty} - \psi\int_{t}^{\infty} \right)dy_{m}\;
\\
	\prod_{j=1}^{m} x_{j}^{a} e^{-x_{j}} y_{j}^{b} e^{-y_{j}}
	\frac{\prod_{1\leq j<k \leq m}(x_{k}-x_{j})(y_{k}-y_{j})}{\prod_{1\leq j,k \leq m}(x_{j}+y_{k})}
	\prod_{1\leq j<k \leq m}(x_{k}-x_{j})(y_{k}-y_{j})  ,
\end{multline}
for the cut-off values of $0<s,t<\infty$ and the generating function variables are $\xi, \psi$ and are taken as indeterminates, possibly complex.
The definition \eqref{CL2M_gapGF} corresponds to the normalisation integral formula with the weight function $ w(x,y)=x^a e^{-x}y^a e^{-y} $ replaced by the piecewise weight 
$ w(x,y)=(1-\xi \chi_{x>s})(1-\psi \chi_{y>t})x^a e^{-x}y^b e^{-y} $.
\end{definition}

The final step in our linkage is via the result in \cite{FK_2016}.
\begin{proposition}[\cite{FK_2016}, Eq. (3.11)]
The unconstrained Bures-Hall generating function is related to the Cauchy-Laguerre two matrix model generating function by
\begin{equation}\label{UBHgf_CL2MM_gf}
	\left( Z^\text{\rm UB-H}(s;m,a,\xi) \right)^2 = 2^m Z^\text{\rm C-L2M}(s,s;m,a,a+1,\xi,\xi) .
\end{equation}
\end{proposition}
\begin{proof}
The proof in Prop.~1 and Cor.~2 of \cite{FK_2016} applies in our situation as the partition functions are defined with densities, denoted there as $ \alpha(z) $, $ z\in \R_{+} $,
which can be piece-wise continuous as are our weights $ w_1(x) = (1-\xi \chi_{x>s})x^a e^{-x} $, $ w_2(y) = (1-\psi \chi_{y>t})y^b e^{-y} $.
The logic of their proof follows through unaltered. 
\end{proof}

\subsection{Bi-moments, Determinant and Pfaffian Structures for finite, fixed $m$}

We adopt the standard definition of the incomplete gamma function $ \Gamma(a,z) $ as per Eq.~8.2.E2 of \cite{DLMF}.
In addition we require the definition of a two-variable extension to the upper, incomplete gamma function.
\begin{definition}
Let $x \in \C\backslash (-\infty,0]$, $y \in \C\backslash (-\infty,-x]$ and ${\rm Re}(a)>-1$. Then define $\Gamma_{2}$ by
\begin{equation}\label{Gamma2}
	\Gamma_{2}(a;x,y) := \int_{x}^{\infty} du\; e^{-u}u^{a}(u+y)^{-1} .
\end{equation} 
\end{definition}
The properties of $\Gamma_{2}$ include:
\begin{enumerate}[(i)]
	\item $ \Gamma_{2}(a;x,0) = \Gamma(a,x) $,
	\item $ \Gamma_{2}(a;0,y) = \Gamma(1+a) y^{a}e^{-y}\Gamma(-a,y) $,
	\item Shift up in the exponent $a$
	\begin{equation}
		 \Gamma_{2}(a+1;x,y) + y\Gamma_{2}(a;x,y) = \Gamma(1+a,x) ,
	\end{equation}
	\item Derivative with respect to $x$
	\begin{equation}
		\partial_{x}\Gamma_{2}(a;x,y) = -x^{a}e^{-x}(x+y)^{-1} ,
	\end{equation}
	\item Derivative with respect to $y$
	\begin{equation}
		\partial_{y}\Gamma_{2}(a;x,y) = \frac{a+y}{y}\Gamma_{2}(a;x,y) - \frac{a}{y}\Gamma(a,x) - x^{a}e^{-x}(x+y)^{-1} .
	\end{equation}
\end{enumerate}

The piece-wise discontinuous bi-variate weight is
\begin{equation}\label{piece-weight}
	w(x,y;s,t) = (1-\xi \chi_{x>s})(1-\psi \chi_{y>t})x^a e^{-x}y^b e^{-y} ,
\end{equation}
and the bi-moments $M_{j,k}$ are defined by
\begin{equation}\label{CL2MM_matrix}
	M_{j,k}(s,t;a,b;\xi,\psi) := \left( \int_{0}^{\infty} - \xi\int_{s}^{\infty} \right)dx \left( \int_{0}^{\infty} - \psi\int_{t}^{\infty} \right)dy 
	\frac{1}{x+y} x^{a+j}y^{b+k} e^{-x-y} .
\end{equation}

A key representation of the C-L2M generating function is as a determinant of the foregoing moments.
\begin{proposition}
The Cauchy-Laguerre two matrix model generating function is also given by the bi-moment determinant $m\geq 1$
\begin{equation}\label{CL2MM_det}
	Z^\text{\rm C-L2M}(s,t;m,a,b;\xi,\psi) =  \frac{1}{C^\text{\rm C-L2M}(m,a,b)}\det\left( M_{j,k} \right)^{m-1}_{j,k=0}	,
\end{equation}
where $s,t>0$, $a,b>-1$, $j,k\geq 0$, $\xi,\psi\in\C$.
\end{proposition}

Properties of the bi-moments:
\begin{enumerate}[(i)]
\item
	$x \leftrightarrow y$ species exchange along with their associated parameters
	\begin{equation}
		M_{k,j}(t,s;b,a;\psi,\xi) = M_{j,k}(s,t;a,b;\xi,\psi)
	\end{equation}
\item 
	$ M_{j,k}(s,t;a,b;\xi,\psi) $ is a function of sums $ a+j $, $ b+k $ only.
	Therefore we often discuss the abbreviation
	\begin{equation}
	M_{j,k}(s,t;a,b;\xi,\psi) = M_{0,0}(s,t;a+j,b+k;\xi,\psi) =: M(s,t;a+j,b+k;\xi,\psi) .
	\end{equation}
\item
	The key structural relation which holds for the bi-moment matrix due to the Cauchy kernel independently of the weight itself, 
	so long as all moments exist $ M_{j,k} < \infty $ $j,k = 0, \ldots, \infty $, is
	\begin{equation}\label{Mjk_key}
		M_{j+1,k} + M_{j,k+1} = \left[ \Gamma(a+j+1) - \xi \Gamma(a+j+1,s) \right]\left[ \Gamma(b+k+1) - \psi \Gamma(b+k+1,t) \right] .
	\end{equation}
	i.e. a factorisation into a product of row and column dependent factors.
\item
	A special case is
	\begin{equation}
		M_{j,k}(0,0;a,b;\xi,\psi) = (1-\xi)(1-\psi) M_{j,k}(s,t;a,b;0,0) = (1-\xi)(1-\psi) \frac{\Gamma(a+j+1)\Gamma(b+k+1)}{a+b+j+k+1},
	\end{equation}
\item
	A limiting case, as $s,t \to \infty$, independently of $ \xi, \psi $
	\begin{equation}
		M_{j,k}(s,t;a,b;\xi,\psi) \mathop{\to}\limits_{s,t \to \infty} \frac{\Gamma(a+j+1)\Gamma(b+k+1)}{a+b+j+k+1}
	\end{equation}
\end{enumerate}

A consequence of the semi-classical character of the weight \eqref{piece-weight} is that the spectral $x,y$ derivatives satisfy linear differential relations
\begin{align}
	x\partial_{x} w & = (a-x)w - \xi e^{-x-y}x^{a+1}y^{b} \delta_{x-s}(1-\psi \chi_{y>t}) ,
\label{wCL_spectral-diff:a}
\\
	y\partial_{y} w & = (b-y)w - \psi e^{-x-y}x^{a}y^{b+1} \delta_{y-t}(1-\xi \chi_{x>s}) .
\label{wCL_spectral-diff:b}
\end{align}

\begin{lemma}
The bi-moments $M_{j,k}$ for $ j,k\geq 0 $ satisfy two identities
\begin{equation}\label{M_Id:a}
	\int_{0}^{\infty}\int_{0}^{\infty} \frac{dxdy}{(x+y)^2} w(x,y) x^{j+1}y^{k}
	= (j+1+a)M_{j,k} - M_{j+1,k} -\xi s^{a+j+1}e^{-s}\int_{0}^{\infty}(1-\psi \chi_{y>t}) \frac{dy}{s+y}e^{-y}y^{b+k} ,
\end{equation}
and
\begin{equation}\label{M_Id:b}
	\int_{0}^{\infty}\int_{0}^{\infty} \frac{dxdy}{(x+y)^2} w(x,y) x^{j}y^{k+1}
	= (k+1+b)M_{j,k} - M_{j,k+1} -\psi t^{b+k+1}e^{-t}\int_{0}^{\infty}(1-\xi \chi_{x>s}) \frac{dx}{x+t}e^{-x}x^{a+j} .
\end{equation}
\end{lemma}
\begin{proof}
The first follows from \eqref{wCL_spectral-diff:a} and the second from \eqref{wCL_spectral-diff:b}.
\end{proof}

The following result is the extension of Eq. (4.1) \cite{BGS_2014} for the basic evaluation of the Cauchy-Laguerre moment. 
\begin{proposition}
The parameters $ \alpha_{n} $, $\beta_{n} $ are given by 
\begin{equation}
 	\alpha_{n} = \left[ \Gamma(a+n+1) - \xi \Gamma(a+n+1,s) \right] ,
\qquad
 	\beta_{n} = \left[ \Gamma(b+n+1) - \psi \Gamma(b+n+1,t) \right] .
\end{equation}
The bi-moment matrix elements are given by 
\begin{multline}\label{CL2MM_moment-matrix}
	(a+b+1) M(s,t;a,b;\xi,\psi) =
	\left[ \Gamma(a+1) - \xi \Gamma(a+1,s) \right]\left[ \Gamma(b+1) - \psi \Gamma(b+1,t) \right]
\\
    + \xi s^{a+1}e^{-s} \left[ \Gamma(b+1)e^{s}s^{b}\Gamma(-b,s) - \psi \Gamma_{2}(b;t,s) \right]
    + \psi t^{b+1} e^{-t} \left[ \Gamma(a+1)e^{t}t^{a}\Gamma(-a,t) - \xi \Gamma_{2}(a;s,t) \right] ,
\end{multline}
where $ \Gamma(a,x) $ is the standard upper incomplete Gamma function, see Eq. (8.6.4) of \cite{DLMF}, and $ \Gamma_{2}(a;x,y) $ is a two variable extension defined by \eqref{Gamma2}.
\end{proposition}
\begin{proof}
This follows by combining \eqref{M_Id:a} and \eqref{M_Id:b} with \eqref{Mjk_key}.
\end{proof}

Another consequence of the semi-classical character of the weight \eqref{piece-weight} is that the deformation $s,t$ derivatives also satisfy linear differential relations
\begin{align}
	\partial_{s} w & = \xi \delta_{x-s}(1-\psi \chi_{y>t}) e^{-x-y}x^{a}y^{b} ,
\label{wCL_deform-diff:a}
\\
	\partial_{t} w & = \psi \delta_{y-t}(1-\xi \chi_{x>s}) e^{-x-y}x^{a}y^{b} .
\label{wCL_deform-diff:b}
\end{align}
\begin{proposition}
The bi-moments satisfy the pair of linear third order PDEs in the deformation variables $s,t$
\begin{align}
	\left[ s\partial_{s} + s-a \right](s+t) \partial_{s}\partial_{t} M & = 0 ,
\label{MCL_deform-pde:s}
\\
	\left[ t\partial_{t} + t-b \right](s+t) \partial_{s}\partial_{t} M & = 0 .
\label{MCL_deform-pde:t}
\end{align}
\end{proposition}
\begin{proof}
We note that $ M(s,t;a,b;\xi,\psi) $ satisfies
\begin{equation}
	\partial_{s} M = \xi s^{a}e^{-s} \left( \int_{0}^{\infty} - \psi\int_{t}^{\infty} \right) \frac{dy}{s+y} y^{b} e^{-y} ,
\qquad
	\partial_{t} M = \psi t^{b}e^{-t} \left( \int_{0}^{\infty} - \xi\int_{s}^{\infty} \right) \frac{dx}{x+t} x^{a} e^{-x} ,
\end{equation}
and thus
\begin{equation}
	\partial_{s}\partial_{t} M = \xi\psi s^{a}t^{b}\frac{ e^{-s-t}}{s+t}.
\end{equation}
The left differential factors in \eqref{MCL_deform-pde:s},\eqref{MCL_deform-pde:t} arise as the annihilators of the univariate densities.
\end{proof}

Adapting the results in \cite{FK_2016} we have a Pfaffian representation of the unconstrained B-H generating function.
\begin{proposition}
The unconstrained Bures-Hall generating function can be written as
\begin{equation}\label{UC_Pfaffian-even}
	Z^\text{\rm UB-H}(s;m,a;\xi) =
	\pm \frac{m!}{C^\text{\rm UB-H}(m,a)} {\rm Pf} \left( M^\text{\rm UB-H}_{j,k} \right)_{0\leq j,k \leq m-1} ,
\end{equation}
if $m$ is even and 
\begin{equation}\label{UC_Pfaffian-odd}
	Z^\text{\rm UB-H}(s;m,a;\xi) = \pm \frac{m!}{C^\text{\rm UB-H}(m,a)} {\rm Pf}
 	\left( 	\begin{array}{cc}
				\left( 0 \right) & \left( m^\text{\rm UB-H}_{k} \right)_{0\leq k \leq m-1} \cr
				\left( -m^\text{\rm UB-H}_{j} \right)_{0\leq j \leq m-1} & \left( M^\text{\rm UB-H}_{j,k} \right)_{0\leq j,k \leq m-1} \cr
			\end{array}
	\right) ,
\end{equation}
if $m$ is odd.
The matrix elements are given by
\begin{equation}\label{UC_mj}
	m^\text{\rm UB-H}_{j} = \Gamma(a+1+j)-\xi\Gamma(a+1+j,s) ,
\end{equation}
and
\begin{multline}\label{UC_Mjk}
	M^\text{\rm UB-H}_{j,k} = (j-k) \left[ \Gamma(a+1+j)-\xi\Gamma(a+1+j,s) \right]\left[ \Gamma(a+1+k)-\xi\Gamma(a+1+k,s) \right]
\\
	+ 2\xi s^{2a+2+j+k} \left[ \Gamma(a+2+j)\Gamma(-a-1-j,s) - \Gamma(a+2+k)\Gamma(-a-1-k,s) \right]
\\
	+ 2\xi^2 e^{-s} \left[ s^{a+1+j}\Gamma(a+1+k,s) - s^{a+1+k}\Gamma(a+1+j,s) + s^{a+2+k}\Gamma_{2}(a+j;s,s) - s^{a+2+j}\Gamma_{2}(a+k;s,s)
					\right] .
\end{multline}
\end{proposition}
\begin{proof}
The Pfaffian structures follow immediately from Eq. (A.3) and (A.4) of \cite{FK_2016} and the matrix elements from Eq. (A.5) of this work 
along with the evaluation \eqref{CL2MM_moment-matrix}.
\end{proof}

\section{Bi-orthogonal Polynomials, Recurrence Relations and Christoffel-Darboux formulae in the General Setting}\label{Cauchy-BOPS}
\setcounter{equation}{0}

In this section we are going to lay-out the essential results for the bi-orthogonal polynomial system for the C-L2M ensemble 
but do this in the most general setting possible rather than un-necessarily specialise. 
A lot of the earlier results do carry over and we will briefly recall these in our notation, 
however we are obliged to complete the picture when certain explicit details are missing.
For general densities $ w_1(x), w_2(y) $ with supports $ S_1\times S_2 \subset \R^2_{+} $ we consider the Cauchy ensemble with the product form $ w(x,y)=w_1(x)w_2(y) $
as per (we switch notation henceforth, $m \mapsto n$) 
\begin{multline}\label{Cxy_jpdf}
	\mathcal{P}^\text{C}(x_{1},\ldots,x_{n};y_{1},\ldots,y_{n})
	= \frac{1}{Z^\text{C}_{n}(n!)^2}
	\prod_{j=1}^{n} w_1(x_{j})w_2(y_{j})
\\	  \times
	\frac{\prod_{1\leq j<k \leq n}(x_{k}-x_{j})(y_{k}-y_{j})}{\prod_{1\leq j,k \leq n}(x_{j}+y_{k})}
	\prod_{1\leq j<k \leq n}(x_{k}-x_{j})(y_{k}-y_{j}) .
\end{multline}
Thus the only coupling of the two one-dimensional systems is through the Cauchy kernel.

As before the bi-moments $M_{j,k}$ are defined by, and we require existence of all members,
\begin{equation}\label{Cxy_matrix}
	M_{j,k} = \int_{S_1\times S_2}dxdy \frac{w(x,y)}{x+y} x^{j}y^{k} < \infty ,\quad \forall j,k \geq 0 .
\end{equation}
The Cauchy partition or generating function is the determinant of the truncated Gram matrix $ \mathsf{M} = \left( M_{j,k} \right)_{j,k \geq 0} $ of foregoing moments, $n\geq 0$
\begin{equation}\label{Cxy_det}
	Z^\text{C}_{n} =  \det\left( M_{j,k} \right)^{n-1}_{j,k=0}, \; n\geq 1, \quad Z^\text{C}_{0}=1 .
\end{equation}
One can define an inner product over polynomial spaces $ \cup_{n\geq 0}\Pi_{n}[x] $ using the weight previously defined.
Let $ f,g \in \cup_{n\geq 0}\Pi_{n}[x] $ then
\begin{equation}\label{I-P}
	\langle f,g \rangle := \int_{S_1\times S_2} dxdy \frac{w(x,y)}{x+y}f(x)g(y) .
\end{equation}
The Cauchy bi-orthogonal system with respect to the above weight function 
consists of two sequences of {\it normalised bi-orthogonal polynomials} $\{P_{n}(x),Q_{n}(y)\}^{\infty}_{n=0}$ satisfying the orthogonality relation
\begin{equation}\label{BOPS_norm}
	\langle P_{m},Q_{n} \rangle = \delta_{m,n},\quad P_{n} = S_{n}x^{n}+\Pi_{n-1}[x], \quad Q_{n}(y) = S_{n}y^{n}+\Pi_{n-1}[y] ,
\end{equation}
with leading coefficients indicated. The {\it monic polynomials} $\{\mathpzc{p}_{n}(x),\mathpzc{q}_{n}(y)\}^{\infty}_{n=0}$ of the system are related via
\begin{equation}\label{monic_PQ}
	P_{n}(x) = \frac{1}{\sqrt{h_{n}}} \mathpzc{p}_{n}(x), 
\quad
	Q_{n}(y) = \frac{1}{\sqrt{h_{n}}} \mathpzc{q}_{n}(y),
\end{equation}
where the monic norm and leading coefficient are related by
\begin{equation}
	h_{n} = \frac{1}{S^2_{n}} = \frac{Z^\text{C}_{n+1}}{Z^\text{C}_{n}} .
\label{norm+leading-coeff}
\end{equation}

In \cite{BGS_2010} several kernels are defined, and we will define the first of the analogous reproducing kernels here as
\begin{equation}
	K^{0,0}_{n}(x,y) := \sum_{l=0}^{n} P_l(x)Q_l(y) .
\label{kernel_00}
\end{equation}
Note that this differs slightly from the conventional definitions used in random matrix theory.
The three basic properties that this kernel satisfies are:
the normalisation
\begin{equation}
	\int_{S_1\times S_2}dxdy \frac{w(x,y)}{x+y} K^{0,0}_{n}(x,y) = n+1,
\end{equation}
the reproducing property for any polynomials $ p(x) \in \Pi_{n}[x] $, $ q(y) \in \Pi_{n}[y] $
\begin{equation}
	\int_{S_1\times S_2}dudv \frac{w(u,v)}{u+v} p(u)K^{0,0}_{n}(x,v) = p(x) ,
\qquad
	\int_{S_1\times S_2}dudv \frac{w(u,v)}{u+v} K^{0,0}_{n}(u,y)q(v) = q(y) ,
\end{equation}
and as a projection operator
\begin{equation}
	\int_{S_1\times S_2}dudv \frac{w(u,v)}{u+v} K^{0,0}_{n}(x,v)K^{0,0}_{n}(u,y) = K^{0,0}_{n}(x,y) .
\end{equation}
Subsequent definitions of the other kernels follow at the beginning of \S \ref{C-D_formulae}, see \eqref{kernel_01}, \eqref{kernel_10} and \eqref{kernel_11}.

In addition to the determinantal formula \eqref{Cxy_det} we will find that every aspect of the bi-orthogonal system will have 
determinantal representations involving bordered bi-moment matrices. Our first example is also a well-known result, 
\begin{equation}
	\mathpzc{p}_{n}(x) = \frac{1}{Z^\text{C}_{n}} \det \left( 
	\begin{array}{cc} 
		  \left( M_{j,k} \right)_{0 \leq j \leq n \atop 0 \leq k \leq n-1}
		& \left( x^j \right)_{0 \leq j \leq n}
	\end{array} \right),
\qquad	
	\mathpzc{q}_{n}(y) = \frac{1}{Z^\text{C}_{n}} \det \left( 
	\begin{array}{c}
		  \left( M_{j,k} \right)_{0 \leq j \leq n-1 \atop 0 \leq k \leq n} \cr
		  \left( y^k \right)_{0 \leq k \leq n}
	\end{array} \right) . 
\label{border_PQ} 
\end{equation}
This result demonstrates that a necessary condition for the existence of the bi-orthogonal system is that the norm $ h_{n} $ or the partition functions $ Z^\text{C}_{n} $
should be non-vanishing $ n \geq 0 $ and using Cramer's rules for solving the linear system for the polynomial coefficients shows this to be sufficient. 
The kernel has the bordered bi-moment representation
\begin{equation}
	K_{n}(x,y) = -\frac{1}{Z_{n+1}} \det \left( 
	\begin{array}{cc}
		\left( M_{j,k} \right)_{0 \leq j \leq n \atop 0 \leq k \leq n} & \left( x^{j} \right)_{0 \leq j \leq n} \cr
		 \left( y^k \right)_{0 \leq k \leq n} & 0 
	\end{array} \right) .
\label{border_K} 
\end{equation}

Furthermore we define averages for symmetric functions of the eigenvalues $ f(x_1, \ldots, x_{n};y_1, \ldots, y_{n}) $ via the JPDF \eqref{Cxy_jpdf} by
\begin{equation}
	{\mathcal I}_{n}\left[ f(x,y) \right]
	:= \int_{\mathbb{R}_{+}^{n}\times \mathbb{R}_{+}^{n}}d^nx\,d^ny\; \mathcal{P}^\text{C}(x_1, \ldots, x_{n};y_1, \ldots, y_{n})f(x_1, \ldots, x_{n};y_1, \ldots, y_{n}) .
\end{equation}
All our examples will be products over the eigenvalues $ f(x_1, \ldots, x_{n})g(y_1, \ldots, y_{n}) = \prod_{j,k=1}^{n}f(x_j)g(y_k) $
and the averages possess the normalisation ${\mathcal I}_{n}\left[ 1\right] = 1$.
It is a straightforward task to show that both polynomials can also be expressed as averages of characteristic polynomials within the two matrix model by
\begin{equation}
	\mathpzc{p}_{n}(u) = {\mathcal I}_{n}\left[ (u-x) \right], \qquad
	\mathpzc{q}_{n}(v) = {\mathcal I}_{n}\left[ (v-y) \right] .
\end{equation}
A simple method to derive these latter results will be given in the proof of Prop. \ref{CD_sum}.
Furthermore the kernel can be expressed as the eigenvalue JPDF average of the product of two characteristic polynomials
\begin{equation}
	K^{0,0}_{n}(x,y) = \frac{1}{h_{n}}{\mathcal I}_{n}\left[ (u-x)(v-y) \right] .
\label{charpoly_K}
\end{equation}

\subsection{Recurrence Relations}

We define the semi-infinite vectors of monomials and polynomials by
\begin{align}
	\mathsf{x}^{T} & = (1,x,x^2,\ldots) & \mathsf{y}^{T} & = (1,y,y^2,\ldots) \\
	\mathsf{P}^{T} & = (P_0(x),P_1(x),P_2(x),\ldots) & \mathsf{Q}^{T} & = (Q_0(y),Q_1(y),Q_2(y),\ldots) .
\end{align}
We have the lower triangular matrices $ \mathsf{S}=\left( S_{n,j} \right)_{n,j\geq 0}, \mathsf{T}=\left( T_{n,j} \right)_{n,j\geq 0} $ relating these bases via
\begin{equation}\label{PQ_Expand_Coeff}
	\mathsf{P} = \mathsf{S}\mathsf{x}, \quad \mathsf{Q} = \mathsf{T}\mathsf{y}, 
\end{equation} 
and the factorisation of the inverse Gram matrix $ \mathsf{M}=\left( M_{j,k} \right)_{j,k\geq 0} $ utilising this
\begin{equation}\label{inverse_Gram}
	\mathsf{T}^T \mathsf{S} = \mathsf{M}^{-1} .
\end{equation}
\begin{proposition}[\cite{BGS_2010}]
Let two sequences of univariate moments be defined as $j\geq 0$
\begin{equation}\label{alpha-beta}
	\alpha_{j} := \int_{S_1}dx\; w_1(x)x^{j} ,\qquad 
	\beta_{j} := \int_{S_2}dy\; w_2(y)y^{j} .	
\end{equation}
A key identity, which we call the Cauchy relation because it follows directly from the Cauchy kernel, for the bi-moments is the rank-1 condition
\begin{equation}\label{Cauchy_Id}
	M_{j+1,k}+M_{j,k+1} = \alpha_{j}\beta_{k} .
\end{equation}
Any three of the data sets $ \{ \alpha_{l} \}^{\infty}_{l=0} $, $ \{ \beta_{l} \}^{\infty}_{l=0} $, $ \{ M_{l,0} \}^{\infty}_{l=0} $, $ \{ M_{0,l} \}^{\infty}_{l=0} $
are algebraically independent and uniquely defines the whole moment data set through \eqref{Cauchy_Id}. 
The above four sets are subject to the relation
\begin{equation}
	M_{0,l+1} + (-1)^{l}M_{l+1,0} + \sum_{j=0}^{l}(-1)^{j+1}\alpha_{j}\beta_{l-j} = 0, \quad l\geq 0 .
\end{equation} 
\end{proposition}
The Cauchy identity stamps a structure on the bi-moment matrix which is the extension of Hankel $M_{j+1,k}=M_{j,k+1}$ and Toeplitz $M_{j,k}=M_{j+1,k+1}$ structures in orthogonal or complex bi-orthogonal systems.
Let the semi-infinite shift-matrix be denoted $ \mathsf{\Lambda} = \left( \delta_{j+1,k} \right)_{j,k\geq 0} $.
In matrix form identity \eqref{Cauchy_Id} reads as a rank-one condition
\begin{equation}
	\mathsf{\Lambda}\mathsf{M}+\mathsf{M}\mathsf{\Lambda}^{T} = \mathsf{\alpha}\otimes\mathsf{\beta}^{T} .
\end{equation}
The spectral multiplication for both sets of polynomials define the lower Hessenberg matrices $ \mathsf{X}, \mathsf{Y} $
\begin{equation}\label{Hessenberg}
	x \mathsf{P} = \mathsf{X}\mathsf{P}, \quad
	y \mathsf{Q}^{T} = \mathsf{Q}^{T}\mathsf{Y}^{T} ,
\end{equation}
and which have the factorisations
\begin{equation}
	\mathsf{X} = \mathsf{S}\mathsf{\Lambda}\mathsf{S}^{-1}, \quad
	\mathsf{Y} = \mathsf{T}\mathsf{\Lambda}\mathsf{T}^{-1} . 
\end{equation}
The rank-one relation then can be rewritten as
\begin{equation}
	\mathsf{X}+\mathsf{Y}^{T} = \mathsf{\pi}\otimes\mathsf{\eta}^{T} ,
\end{equation}
where $ \mathsf{\pi} = \mathsf{S}\mathsf{\alpha} $, $ \mathsf{\eta} = \mathsf{T}\mathsf{\beta} $.
Alternatively these new parameters can be defined by, $j\geq 0$
\begin{equation}\label{pi-eta}
	\pi_{j} := \int_{S_1}dx\; w_1(x)P_{j}(x) ,\qquad 
	\eta_{j} := \int_{S_2}dy\; w_2(y)Q_{j}(y) .	
\end{equation}

The three leading diagonals of the multiplication operators $ \mathsf{X}, \mathsf{Y} $, $ j-k \in \{-1,0,1\} $ will feature prominently in what follows.
The remaining components of these matrices are related via row relations linking the sub-diagonal $ j-k=2 $ case to those on the left 
\begin{equation}
	X_{j,k+1} = \frac{\eta_{k+1}}{\eta_{0}} X_{j,0}, \quad j\geq k+3 ; \qquad 
	Y_{j,k+1} = \frac{\pi_{k+1}}{\pi_{0}} Y_{j,0}, \quad j\geq k+3 ,
\end{equation}
or via column relations linking the sub-diagonal $ j-k=2 $ case to those below
\begin{equation}
	X_{j+1,k} = \frac{\pi_{j+1}}{\pi_{k+2}} X_{k+2,k}, \quad j\geq k+2 ; \qquad 
	Y_{j+1,k} = \frac{\eta_{j+1}}{\eta_{k+2}} Y_{k+2,k}, \quad j\geq k+2 . 
\end{equation}
Collecting all these results we can state in summary.
\begin{proposition}
The multiplication operators $ \mathsf{X}, \mathsf{Y} $ have the evaluations along the leading diagonals
\begin{equation}\label{XY-elements}
	\begin{array}{lll}
		j-k=-1 & X_{n,n+1} = \frac{S_{n}}{S_{n+1}} , & Y_{n,n+1} = \frac{S_{n}}{S_{n+1}} ,  \cr
		j-k=0 & X_{n,n} = \pi_{n}\eta_{n}-\langle P_{n}, yQ_{n} \rangle , & Y_{n,n} = \pi_{n}\eta_{n}-\langle xP_{n}, Q_{n} \rangle , \cr
		j-k=1 & X_{n,n-1} = \pi_{n}\eta_{n-1}-\frac{S_{n-1}}{S_{n}} , & Y_{n,n-1} = \pi_{n-1}\eta_{n}-\frac{S_{n-1}}{S_{n}} , \cr
		j-k\geq 2 & X_{j,k} = \pi_{j}\eta_{k} , & Y_{j,k} = \pi_{k}\eta_{j} . \cr
	\end{array} .
\end{equation}
\end{proposition}
\begin{proof}
We express that matrix element $ X_{n,j} = \langle xP_{n}, Q_j \rangle = \pi_{n}\eta_{j} - \langle P_{n}, yQ_j \rangle $ as a determinant using \eqref{border_PQ} 
and note that few, if any, of elements of the bottom row are non-zero.
Those that are non-zero have simple evaluations except for $ j=n $, which we will encounter subsequently.  
\end{proof}

\begin{definition}
We will spell out a set of conditions, which we denote as {\it generic conditions}, to ensure the validity of subsequent statements:
\begin{itemize}
\item	Finite moments, $ M_{j,k} < \infty $ for all $ j,k \geq 0 $,	
\item	Non-vanishing moment determinants or tau-functions, $ Z^\text{C}_{n} \neq 0 $ for all $ n\geq 0 $, and
\item	Non-vanishing auxiliary coefficients, $ \pi_{n}, \eta_{n} \neq 0 $ for all $ n\geq 0 $.
\end{itemize}
\end{definition}

\begin{definition}
Let the univariate moment sequences $ \pi_{n}, \eta_{n} $ be non-vanishing for all $ n\geq 0 $,
and the diagonal matrices be constructed $ \mathsf{D}_{\pi} = \rm{diag}( \mathsf{\pi}) $, $ \mathsf{D}_{\eta} = \rm{diag}( \mathsf{\eta}) $.
Then it is possible to define the following difference operators
\begin{equation}\label{L-operators}
	\check{\mathsf{L}} := (\mathsf{\Lambda}-\mathsf{1})\mathsf{D}_{\pi}^{-1}, \qquad
	\hat{\mathsf{L}} := \mathsf{D}_{\eta}^{-1}(\mathsf{\Lambda}^{T}-\mathsf{1}) ,
\end{equation} 
which are the left and right annihilators of $ \mathsf{\pi}, \mathsf{\eta}^{T} $ respectively.
Written in terms of components $\check{\mathsf{L}}, \hat{\mathsf{L}} $ have upper and lower bi-diagonal structures respectively
\begin{equation}
	\big( \check{\mathsf{L}} \big)_{j,k} = \frac{1}{\pi_{k+1}}\delta_{j+1,k} - \frac{1}{\pi_{k}}\delta_{j,k}, \qquad
	\big( \hat{\mathsf{L}} \big)_{j,k} = \frac{1}{\eta_{k+1}}\delta_{j,k+1} - \frac{1}{\eta_{k}}\delta_{j,k} .
\end{equation}
If we further require finite moments $ \pi_{n}, \eta_{n} < \infty $ for all $ n\geq 0 $ then their inverses exist
which means that $\check{\mathsf{L}}^{-1}, \hat{\mathsf{L}}^{-1} $ are summation operators with upper and lower triangular structures respectively
\begin{equation}
\big( \check{\mathsf{L}}^{-1} \big)_{j,k} = \begin{cases} -\pi_{k} & j\leq k \\ 0 & j> k \end{cases}, \qquad
\big( \hat{\mathsf{L}}^{-1} \big)_{j,k} = \begin{cases} -\eta_{k} & j\geq k \\ 0 & j< k \end{cases} .
\end{equation}
\end{definition}

One can then deduce recurrence relations for the bi-orthogonal polynomials in the form
\begin{equation}
	x \check{\mathsf{L}} \mathsf{P} = \check{\mathsf{A}} \mathsf{P} ,\quad
	y \hat{\mathsf{L}}^{T} \mathsf{Q} = \hat{\mathsf{B}} \mathsf{Q} ,
\end{equation}
where
\begin{align}
	\check{\mathsf{A}} & := \check{\mathsf{L}}\mathsf{X} = \left( \frac{X_{j+1,k}}{\pi_{j+1}} - \frac{X_{j,k}}{\pi_{j}} \right)_{j,k\geq 0} , 
& 	\check{\mathsf{B}} & := \mathsf{Y}\check{\mathsf{L}}^{T} = \left( \frac{Y_{j,k+1}}{\pi_{k+1}} - \frac{Y_{j,k}}{\pi_{k}} \right)_{j,k\geq 0} ,
\label{AB_Down}\\
	\hat{\mathsf{A}} & := \mathsf{X}\hat{\mathsf{L}} = \left( \frac{X_{j,k+1}}{\eta_{k+1}} - \frac{X_{j,k}}{\eta_{k}} \right)_{j,k\geq 0} , 
& 	\hat{\mathsf{B}} & := \hat{\mathsf{L}}^{T}\mathsf{Y} = \left( \frac{Y_{j+1,k}}{\eta_{j+1}} - \frac{Y_{j,k}}{\eta_{j}} \right)_{j,k\geq 0} .
\label{AB_Up}
\end{align}

All of the foregoing developments show that the results of \cite{BGS_2010} continue to hold in our general setting, subject to the stated conditions,
and consequently the bi-orthogonal polynomial system satisfies two parallel recurrence relations.
\begin{proposition}[\cite{BGS_2010}]
Let the coefficients $ \pi_{n} $, $ \eta_{n} $ be non-vanishing for all $ n \geq 0 $, and the norms $ h_{n} $ or $ Z^\text{C}_{n} $ similarly be non-vanishing.
The bi-orthogonal polynomials $ P_{n}(x) $, $ Q_n(y) $ defined by the general orthogonality condition \eqref{BOPS_norm} satisfy uncoupled, third-order scalar recurrence relations of the form
\begin{align}
	x \left( \frac{1}{\pi_{n+1}}P_{n+1} - \frac{1}{\pi_{n}}P_{n} \right)
	& = r_{n,2}P_{n+2} + r_{n,1}P_{n+1} + r_{n,0}P_{n} + r_{n,-1}P_{n-1} ,
\label{P_recur}
\\
	y \left( \frac{1}{\eta_{n+1}}Q_{n+1} - \frac{1}{\eta_{n}}Q_{n} \right)
	& = s_{n,2}Q_{n+2} + s_{n,1}Q_{n+1} + s_{n,0}Q_{n} + s_{n,-1}Q_{n-1} .
\label{Q_recur}
\end{align}
\end{proposition}
\begin{proof}
From the relation $ S_{n+1} = \pi_{n+1}S_{n+2}r_{n,2} $ it is clear that $ \pi_{n+1}, S_{n+1} \neq 0, \infty $ for the recurrence system \eqref{P_recur} 
to be solved in the forward direction.
\end{proof}

\begin{remark}
Explicit evaluations for the coefficients $ \pi_{n}, \eta_{n}, r_{n,j}, s_{n,j} $ are given in Theorem 2.1 of \cite{WW_2021} for the base Cauchy-Laguerre system
when $ s=t=0 $, $ \xi=\psi=0 $ as rational functions of $ n,a,b $ for $ a,b, a+b > -2 $.
We give further details of this case in \S\ref{Undeformed_Case}.
\end{remark}

\begin{remark}
The sequence of four recurrence coefficients are algebraically independent $ \{r_{n,2}, r_{n,1}, r_{n,0}, r_{n,-1} \} $.
From the four leading diagonals of matrix $ \mathsf{X} $ we have the recurrence coefficients
\begin{align}
	r_{n,2} & = \frac{X_{n+1,n+2}}{\pi_{n+1}} = \frac{S_{n+1}}{S_{n+2}\pi_{n+1}} ,
\label{rCFF:+2}\\
	r_{n,1} & = \frac{X_{n+1,n+1}}{\pi_{n+1}} - \frac{X_{n,n+1}}{\pi_{n}} = \frac{Y_{n+1,n}}{\pi_{n}} - \frac{Y_{n+1,n+1}}{\pi_{n+1}} ,
\label{rCFF:+1}\\
	r_{n,0} & = \frac{X_{n+1,n}}{\pi_{n+1}} - \frac{X_{n,n}}{\pi_{n}} = \frac{Y_{n,n}}{\pi_{n}} - \frac{Y_{n,n+1}}{\pi_{n+1}} ,
\label{rCFF:+0}\\
	r_{n,-1} & = \frac{X_{n+1,n-1}}{\pi_{n+1}} - \frac{X_{n,n-1}}{\pi_{n}} = \frac{S_{n-1}}{S_{n}\pi_{n}} .
\label{rCFF:-1}
\end{align} 
Alternatively one has
\begin{gather}
	\frac{X_{n+1,n+2}}{\pi_{n+1}} = r_{n,2},
	\frac{X_{n+1,n+1}}{\pi_{n+1}} = r_{n-1,2}+r_{n,1},
	\frac{X_{n+1,n}}{\pi_{n+1}} = r_{n-2,2}+r_{n-1,1}+r_{n,0},
\\
	\frac{X_{n+1,n-1}}{\pi_{n+1}} = r_{n-3,2}+r_{n-2,1}+r_{n-1,0}+r_{n,-1} .
\label{X_from_r}
\end{gather}
Symmetric statements can be made for the $Q$-recurrence coefficients and matrix $ \mathsf{Y} $
and from this we deduce relations for the $ s_{n,k} $ coefficients of the $Q$ polynomials system in terms of the $ r_{n,k} $ coefficients
\begin{align}
	s_{n,2} & = \frac{Y_{n+1,n+2}}{\eta_{n+1}} = \frac{S_{n+1}}{S_{n+2}\eta_{n+1}} = \frac{\pi_{n+1}}{\eta_{n+1}}r_{n,2} ,
\label{sCFF:+2}\\
	s_{n,1} & = \frac{Y_{n+1,n+1}}{\eta_{n+1}} - \frac{Y_{n,n+1}}{\eta_{n}} = \frac{X_{n+1,n}}{\eta_{n}} - \frac{X_{n+1,n+1}}{\eta_{n+1}} 
				= \frac{\pi_{n+1}}{\eta_{n}}\left[ r_{n-2,2}+r_{n-1,1}+r_{n,0} \right] - \frac{\pi_{n+1}}{\eta_{n+1}}\left[ r_{n-1,2}+r_{n,1} \right] ,
\label{sCFF:+1}\\
	s_{n,0} & = \frac{Y_{n+1,n}}{\eta_{n+1}} - \frac{Y_{n,n}}{\eta_{n}} = \frac{X_{n,n}}{\eta_{n}} - \frac{X_{n,n+1}}{\eta_{n+1}}
				= \frac{\pi_{n}}{\eta_{n}}\left[ r_{n-2,2}+r_{n-1,1} \right] - \frac{\pi_{n}}{\eta_{n+1}} r_{n-1,2} ,
\label{sCFF:+0}\\
	s_{n,-1} & = \frac{Y_{n+1,n-1}}{\eta_{n+1}} - \frac{Y_{n,n-1}}{\eta_{n}} = \frac{S_{n-1}}{S_{n}\eta_{n}} 
				= \frac{\pi_{n}}{\eta_{n}}r_{n,-1} .
\label{sCFF:-1}
\end{align} 
\end{remark}

\begin{remark}
Given knowledge of the recurrence relation and the $ \pi, \eta $ coefficients one can set up the linear system of equations for the polynomial coefficients $ S_{j,k}, 0\leq k\leq j $
in the monomial basis as a recursive system of first inhomogeneous difference equations, all with identical homogeneous parts
\begin{align}
	k=n+2: \qquad
	\frac{1}{\pi_{n+1}} S_{n+1,n+1} - r_{n,2}S_{n+2,n+2} & = 0 ,
\\
	k=n+1: \qquad
	\frac{1}{\pi_{n+1}} S_{n+1,n} - r_{n,2}S_{n+2,n+1} & = \frac{1}{\pi_{n}} S_{n,n} + r_{n,1}S_{n+1,n+1} ,
\\
	k=n: \qquad
	\frac{1}{\pi_{n+1}} S_{n+1,n-1} - r_{n,2}S_{n+2,n} & = \frac{1}{\pi_{n}} S_{n,n-1} + r_{n,1}S_{n+1,n} + r_{n,0}S_{n,n} ,
\\
	1\leq k\leq n-1: \qquad
	\frac{1}{\pi_{n+1}} S_{n+1,k-1} - r_{n,2}S_{n+2,k} & = \frac{1}{\pi_{n}} S_{n,k-1} + r_{n,1}S_{n+1,k} + r_{n,0}S_{n,k} + r_{n,-1}S_{n-1,k} .
\\
	k=0: \qquad
	0 & =   r_{n,2}S_{n+2,0} + r_{n,1}S_{n+1,0} + r_{n,0}S_{n,0} + r_{n,-1}S_{n-1,0} .
\end{align}
Clearly this is only possible if $ \pi_{n+1}r_{n,2} \neq 0, \infty $.
By convention the initial values of the recurrence coefficients are set at 
$ r_{-2,2} = 0 $, $ r_{-1,2} = S_0/(\pi_0 S_{1}) $, $ r_{0,2} = S_1/(\pi_1 S_{2}) $, $ r_{-1,1} = M_{1,0}/(\pi_0 M_{0,0}) $. 
One can readily solve these to generate a finite set of leading polynomial coefficients: $ S_{n+1,n+1} = S_{n+1} $ and
\begin{align}
	\frac{S_{n+1,n}}{S_{n+1}} & = -\sum_{l=0}^{n} X_{l,l} ,
\label{P_subleading}\\
	\frac{S_{n+1,n-1}}{S_{n+1}} & = \sum_{l=1}^{n}\frac{S_{l-1}}{S_{l}}\left( \frac{S_{l-1}}{S_{l}}-\pi_{l}\eta_{l-1} \right) + \sum_{l=1}^{n} X_{l,l}\sum_{m=0}^{l-1}X_{m,m} . 
\label{P_sub-subleading}
\end{align} 
\end{remark}

In \cite{BGS_2010} the notion of a pair of {\it intertwined} polynomials was introduced and here we develop a more complete description of this notion consisting of two pairs of such polynomials.
\begin{definition}
Let the following polynomial pair sequences be defined using the difference operators $ \check{\mathsf{L}}, \hat{\mathsf{L}} $, which we denote as "down" and "up" respectively,
\begin{align}
	\check{\mathsf{P}} & := \check{\mathsf{L}}\mathsf{P}, \qquad
	\check{\mathsf{Q}}^{T} := \mathsf{Q}^{T}\check{\mathsf{L}}^{-1}
\label{PQ_Down}
\\
	\hat{\mathsf{P}} & := \hat{\mathsf{L}}^{-1}\mathsf{P}, \qquad
	\hat{\mathsf{Q}}^{T} := \mathsf{Q}^{T}\hat{\mathsf{L}} .
\label{PQ_Up}
\end{align} 
\end{definition}
\begin{proposition}[\cite{BGS_2010}]\label{PQintertwined}
The two pairs of intertwined polynomials satisfy the following properties:\\
\begin{enumerate}[(i)]
\item
They are orthonormal pairs $ \langle \check{\mathsf{P}}, \check{\mathsf{Q}}^{T} \rangle = \langle \hat{\mathsf{P}}, \hat{\mathsf{Q}}^{T} \rangle = \mathsf{1} $,
\item
They are related with the pair $ \mathsf{P}, \mathsf{Q} $ via the component-wise relations
\begin{align}
	\check{P}_{n}(x) & = \frac{1}{\pi_{n+1}}P_{n+1}(x) - \frac{1}{\pi_{n}}P_{n}(x), & \check{Q}_{n}(y) & = -\sum_{j=0}^{n} \pi_{j}Q_{j}(y) , 
\label{intertwineDN}\\
	\hat{P}_{n}(x) & = -\sum_{k=0}^{n} \eta_{k}P_{k}(x), & \hat{Q}_{n}(y) & = \frac{1}{\eta_{n+1}}Q_{n+1}(y) - \frac{1}{\eta_{n}}Q_{n}(y) ,
\label{intertwineUP}
\end{align}
\item
These can be inverted under generic conditions by
\begin{align}
	\frac{1}{\pi_{n}}P_{n}(x) & = \sum_{j=0}^{n-1}\check{P}_{j}(x) , &  \pi_{n}Q_{n}(y) & = \check{Q}_{n-1}(y)-\check{Q}_{n}(y), 
\label{invertDN}\\
	\eta_{n}P_{n}(x) & = \hat{P}_{n-1}(x)-\hat{P}_{n}(x) , & \frac{1}{\eta_{n}}Q_{n}(y) & = \sum_{k=0}^{n-1}\hat{Q}_{k}(y) ,
\label{invertUP}
\end{align}
\end{enumerate}
\end{proposition}
Alternative characterising relations for $ \hat{P}_{n}(x),\check{Q}_{n}(y) $ are via univariate reproducing relations
\begin{equation}
	\hat{P}_{n}(x) = -\int_{S_2}dy\; w_2(y) K^{0,0}_n(x,y) ,
\qquad
	\check{Q}_{n}(y) = -\int_{S_1}dx\; w_1(x) K^{0,0}_n(x,y) .
\end{equation}

\subsection{Stieltjes Functions}\label{SF}
Our picture of the bi-orthogonal system is not complete, however, without additional essential components to the polynomial sequences themselves.
Firstly we require analogues of the Steiltjes functions and in this case three types will arise.
In our definitions and usage here we will depart from \cite{BGS_2010} however the correspondences between ours and theirs will be simple.

To start with we define the univariate Stieltjes transforms.
\begin{definition}
We define the first univariate transform
\begin{align}
	f_1(z) & := \int_{S_1} dx \frac{w_1(x)}{z-x} ,\quad z\in \C\backslash S_1 ,
\label{Stieltjes_f1}
\\
	f_2(z) & := \int_{S_2} dy \frac{w_2(y)}{z-y} ,\quad z\in \C\backslash S_2 .
\label{Stieltjes_f2}
\end{align}
\end{definition}
\noindent
In addition we require the cross transforms of these Stieltjes transforms themselves, 
which are constructed from the ones above but with their domains reflected to the interval $ (-\infty,0] $ and using the other univariate density.
\begin{definition}
We define the second type of univariate transform
\begin{align}
	g_1(z) & := -\int_{S_1} dx \frac{w_1(x)}{z-x}f_2(-x) = \int_{S_1\times S_2} dxdy\frac{w(x,y)}{x+y}\frac{1}{z-x} ,\quad z\in \C\backslash S_1 ,
\label{Stieltjes_g1}
	\\
	g_2(z) & := -\int_{S_2} dy \frac{w_2(y)}{z-y}f_1(-y) = \int_{S_1\times S_2} dxdy\frac{w(x,y)}{x+y}\frac{1}{z-y} ,\quad z\in \C\backslash S_2 .
\label{Stieltjes_g2}
\end{align}
\end{definition}
These four functions are not independent because of the relation
\begin{align}
	g_2(z)+g_1(-z) = -f_1(-z)f_2(z) ,
\end{align}
as noted in Eq. (5.2) of \cite{BGS_2010}, and which is another consequence of \eqref{Cauchy_Id}.

\begin{remark}
Another important transform of the weight is the symmetric bi-variate one
\begin{equation}
	G(u,v) = \int_{S_1\times S_2} dxdy\frac{w(x,y)}{x+y}\frac{1}{(u-x)(v-y)} ,
\end{equation}
which is a natural object because it is the formal generating function about $ (u,v)=(\infty,\infty) $ for the bi-moments
\begin{equation}
	G(u,v) = \sum_{j\geq 0}\sum_{k\geq 0} u^{-j-1}v^{-k-1}M_{j,k} .
\end{equation}
However this is simply recovered from the ones defined above through the relations
\begin{equation*}
	(u+v)G(u,v) = g_2(v) - g_2(-u) + f_1(u)f_2(v) - f_1(u)f_2(-u) = g_1(u) - g_1(-v) + f_1(u)f_2(v) - f_1(-v)f_2(v) .	
\end{equation*}
\end{remark}

The analytical character of the Stieltjes functions with respect to $ z\in\C $ in their expansions about $ z=\infty $ is important from many perspectives,
including their role as generating functions for quantities already introduced:
the $ f_{1}, f_{2} $ set have the expansions
\begin{equation}
	f_{1}(z) \mathop{=}\limits_{z \to \infty} \sum_{l\geq 0} \alpha_{l} z^{-l-1} ,
\qquad
	f_{2}(z) \mathop{=}\limits_{z \to \infty} \sum_{l\geq 0} \beta_{l} z^{-l-1} , 
\end{equation}
and the $ g_{1}, g_{2} $ set have the expansions
\begin{align}
	g_{1}(z) \mathop{=}\limits_{z \to \infty} \sum_{l\geq 0} M_{l,0} z^{-l-1} ,
\qquad
	g_{2}(z) \mathop{=}\limits_{z \to \infty} \sum_{l\geq 0} M_{0,l} z^{-l-1} . 
\end{align}

\subsection{Associated Polynomials and Functions}\label{APF}
A second set of companion functions to the bi-orthogonal polynomials is required as the third order recurrence relations require another two independent solutions
in order to furnish a full set of three fundamental solutions.
As a stepping stone to this goal one can consider the following definition of $(n-1)$-th degree polynomial solutions to \eqref{P_recur}, \eqref{Q_recur} but displaced one step down, 
$ n \mapsto n-1 $, namely the {\it associated polynomials of the first type}
\begin{equation}
	\tilde{P}^{(1)}_{n-1}(z) := \int_{S_1} dx w_1(x)\frac{P_{n}(z)-P_{n}(x)}{z-x} ,
\qquad
	\tilde{Q}^{(1)}_{n-1}(z) := \int_{S_2} dy w_2(y)\frac{Q_{n}(z)-Q_{n}(y)}{z-y} ,
\label{AssocPoly_P+Q1}
\end{equation}
which parallels the definitions \eqref{Stieltjes_f1}, \eqref{Stieltjes_f2}. 
What is implied by these latter definitions are in fact appropriate ones for {\it associated functions of the first type}, of $n$-th order
\begin{align}
	P^{(1)}_{n}(z) & := \int_{S_1} dx \frac{w_1(x)}{z-x}P_{n}(x) ,
\label{AssocFn_P1}
\\
	Q^{(1)}_{n}(z) & := \int_{S_2} dy \frac{w_2(y)}{z-y}Q_{n}(y) ,
\label{AssocFn_Q1}
\end{align}
because of the key relations
\begin{equation}
	P^{(1)}_{n}(z) = f_1(z)P_{n}(z) - \tilde{P}^{(1)}_{n-1}(z) ,
\qquad
	Q^{(1)}_{n}(z) = f_2(z)Q_{n}(z) - \tilde{Q}^{(1)}_{n-1}(z) .
\label{AssocPoly-AssocFn_P+Q1}
\end{equation}
Finally we need the cross transforms of these functions, again in parallel to the earlier definitions \eqref{Stieltjes_g1} and \eqref{Stieltjes_g2}, 
and we have the {\it associated functions of the second type}
\begin{align}
	P^{(2)}_{n}(z) & := \int_{S_1\times S_2} dxdy\frac{w(x,y)}{x+y}\frac{P_{n}(x)}{z-x} ,
\label{AssocFn_P2}
\\
	Q^{(2)}_{n}(z) & := \int_{S_1\times S_2} dxdy\frac{w(x,y)}{x+y}\frac{Q_{n}(y)}{z-y} .
\label{AssocFn_Q2}
\end{align}
In analogy with the first pair of associated functions we can define their {\it associated polynomials of the second type}
\begin{equation}
	\tilde{P}^{(2)}_{n-1}(z) := \int_{S_1\times S_2} dxdy \frac{w(x,y)}{x+y}\frac{P_{n}(z)-P_{n}(x)}{z-x} ,
\qquad
	\tilde{Q}^{(2)}_{n-1}(z) := \int_{S_1\times S_2} dxdy \frac{w(x,y)}{x+y}\frac{Q_{n}(z)-Q_{n}(y)}{z-y} ,
\label{AssocPoly_P+Q2}
\end{equation}
and their inter-relations
\begin{equation}
	P^{(2)}_{n}(z) = g_1(z)P_{n}(z) - \tilde{P}^{(2)}_{n-1}(z) ,
\qquad
	Q^{(2)}_{n}(z) = g_2(z)Q_{n}(z) - \tilde{Q}^{(2)}_{n-1}(z) .
\label{AssocPoly-AssocFn_P+Q2}
\end{equation}

\begin{proposition}
The associated functions $ P^{(1)}_{n} $, $ P^{(2)}_{n} $ and the $ Q^{(1)}_{n} $, $ Q^{(2)}_{n} $ satisfy the recurrence relations \eqref{P_recur} and \eqref{Q_recur} respectively.
\end{proposition}
\begin{proof}
We provide a proof only for the $P$-system as identical arguments hold for the $Q$-system.
From one copy of \eqref{P_recur} with $x\mapsto z$ we subtract the original equation, multiply by $w_1(x)$ and divide by $z-x$, and integrate the result over the support $S_1$.
Six of the terms are recognised as associated polynomials of the first kind from the definition \eqref{AssocFn_P1} 
while the remaining two terms give $\pi_{n}, \pi_{n+1}$ coefficients from the integral definition \eqref{pi-eta}. The latter terms then mutually cancel.
Now we utilise \eqref{AssocPoly-AssocFn_P+Q1} and write out each associated polynomial in terms of the associated function and polynomial,
and collect on the latter two sets. And we see that the associated function satisfies the recurrence relation given that the corresponding polynomial does in the first instance.
For the second kind associated polynomials we make the same initial subtraction but multiply by $w(x,y)/(x+y)$ and integrate over $S_1\times S_2$ instead.
Again six of the terms are recognised as associated polynomials of the second kind by \eqref{AssocFn_P2} and now the remaining two are of the form
\begin{equation*}
	 \frac{1}{\pi_{n}}\int_{S_1\times S_2} dxdy\frac{w(x,y)}{x+y}P_{n}(x) = \frac{\delta_{n,0}}{S_0\pi_{n}} .
\end{equation*}
If $n\geq 1$ then both terms vanish. 
Applying \eqref{AssocPoly-AssocFn_P+Q2} we can again split the result into two sets and show that the second kind associated function satisfies \eqref{P_recur} 
given that the polynomials so do.
\end{proof}
By introducing the foregoing associated functions we have completed a crucial task - 
constructing a matrix system of bi-orthogonal polynomials and associated functions which satisfy first order matrix recurrence relations.
To this end we define the $ 3\times 3 $ matrix variables for $ n\geq 0 $
\begin{equation}\label{3x3_PQ}
	\mathrm{P}_{n}(x) := 	\left(
	\begin{array}{ccc}
		P_{n+1}(x) 	& P^{(1)}_{n+1}(x) 	& P^{(2)}_{n+1}(x) \cr
		P_{n}(x) 	& P^{(1)}_{n}(x) 	& P^{(2)}_{n}(x) \cr
		P_{n-1}(x) 	& P^{(1)}_{n-1}(x) 	& P^{(2)}_{n-1}(x) \cr
	\end{array}
			\right), 
\qquad
	\mathrm{Q}_{n}(y) := 	\left(
	\begin{array}{ccc}
		Q_{n+1}(y) 	& Q^{(1)}_{n+1}(y) 	& Q^{(2)}_{n+1}(y) \cr
		Q_{n}(y) 	& Q^{(1)}_{n}(y) 	& Q^{(2)}_{n}(y) \cr
		Q_{n-1}(y) 	& Q^{(1)}_{n-1}(y) 	& Q^{(2)}_{n-1}(y) \cr
	\end{array}
			\right) .
\end{equation}
Thus we have a matrix system of first order recurrences
\begin{equation}
	\mathrm{P}_{n+1} = K_{n} \mathrm{P}_{n}, \qquad
	\mathrm{Q}_{n+1} = L_{n} \mathrm{Q}_{n} ,
\label{PQ_recur}
\end{equation}
where the transfer matrices $ K_{n} $ and $ L_{n} $ take the forms
\begin{align}
	K_{n} & =  	\left(
	\begin{array}{ccc}
	\frac{x}{\pi_{n+1}r_{n,2}}-\frac{r_{n,1}}{r_{n,2}} 	& -\frac{x}{\pi_{n}r_{n,2}}-\frac{r_{n,0}}{r_{n,2}}	& -\frac{r_{n,-1}}{r_{n,2}} \cr
	1 	& 0 	& 0 \cr
	0	& 1 	& 0 \cr
	\end{array}
	\right),
\label{3x3_K}\\
	L_{n} & =  	\left(
	\begin{array}{ccc}
	\frac{y}{\eta_{n+1}s_{n,2}}-\frac{s_{n,1}}{s_{n,2}} 	& -\frac{y}{\eta_{n}s_{n,2}}-\frac{s_{n,0}}{s_{n,2}}	& -\frac{s_{n,-1}}{s_{n,2}} \cr
	1 	& 0 	& 0 \cr
	0	& 1 	& 0 \cr
\end{array}
				\right) .
\label{3x3_L}
\end{align}

\begin{proposition}
In addition to the conditions $ S_{n}\neq 0 $, $ \pi_{n}\neq 0 $ for all $ n\geq 0 $, let
$ \alpha_{1}M_{0,0} \neq \alpha_{0}M_{1,0} $ and $ \beta_{1}M_{0,0} \neq \beta_{0}M_{0,1} $.
Then $ P_{n}(x), P^{(1)}_{n}(x), P^{(2)}_{n}(x) $ are a fundamental set of solutions to \eqref{P_recur} and correspondingly  
$ Q_{n}(x), Q^{(1)}_{n}(x), Q^{(2)}_{n}(x) $ are a fundamental set of solutions to \eqref{Q_recur}.
Their respective Casoratians, for $ n\geq 2 $, are given by
\begin{equation}
	\det \mathrm{P}_{n} = (-1)^{n-1} \pi_{n}\frac{S_{n+1}}{S_{n-1}}, 
\qquad
	\det \mathrm{Q}_{n} = (-1)^{n-1} \eta_{n}\frac{S_{n+1}}{S_{n-1}} .
\label{det_PQ}
\end{equation}
\end{proposition}
\begin{proof}
We give the proof for the case of $ \det \mathrm{P}_{n} $ only as the treatment of $ \mathrm{Q}_{n} $ is identical.
Taking the determinant of \eqref{3x3_K} we get the first order recurrence relation
\begin{equation*}
	\det \mathrm{P}_{n+1} = -\frac{\pi_{n+1}}{\pi_{n}}\frac{S_{n+2}S_{n-1}}{S_{n+1}S_{n}}\det \mathrm{P}_{n} ,
\end{equation*}
which is immediately solved. Furthermore we note the initial value is given by
\begin{align*}
	\det \mathrm{P}_{1} & = \int_{S_1} dx_1 \frac{w_1(x_1)}{x-x_1} \int_{S_1\times S_2} dx_2dy_2 \frac{w(x_2,y_2)}{(x_2+y_2)(x-x_2)}
	\left(
	\begin{array}{ccc}
	P_{2}(x) 	& P_{2}(x_1) 	& P_{2}(x_2) \cr
	P_{1}(x) 	& P_{1}(x_1) 	& P_{1}(x_2) \cr
	P_{0}(x) 	& P_{0}(x_1) 	& P_{0}(x_2) \cr
	\end{array}
	\right),
\\
	& = S_{2}S_{1}S_{0}\int_{S_1} dx_1 w_1(x_1) \int_{S_1\times S_2} dx_2dy_2 \frac{w(x_2,y_2)}{(x_2+y_2)}(x_1-x_2) ,
\\
	& = S_{2}S_{1}S_{0} \left[ \alpha_{1}M_{0,0} - \alpha_{0}M_{1,0} \right] .
\end{align*}
Finally the initial factors simplify because $ S_{1}S_{0}^2 \left[ \alpha_{1}M_{0,0} - \alpha_{0}M_{1,0} \right]=\pi_{1} $.
\end{proof}

At this point in our development we need to examine the $ z\to\infty $ character of the associated functions defined so far and we observe that
\begin{align}
	P^{(1)}_{n}(z) & \mathop{=}\limits_{z \to \infty} \pi_{n}z^{-1} + {\rm O}(z^{-2}) ,
\\
	Q^{(1)}_{n}(z) & \mathop{=}\limits_{z \to \infty} \eta_{n}z^{-1} + {\rm O}(z^{-2}) ,
\\
	P^{(2)}_{n}(z) & \mathop{=}\limits_{z \to \infty} S_{0}^{-1}\delta_{n,0}z^{-1} + {\rm O}(z^{-2}) ,
\\
	Q^{(2)}_{n}(z) & \mathop{=}\limits_{z \to \infty} S_{0}^{-1}\delta_{n,0}z^{-1} + {\rm O}(z^{-2}) .
\end{align}
This reveals a deficiency in that the analysis of Riemann-Hilbert and related problems requires partner functions to the $n$-degree polynomials for which
their expansions are missing the first $ n $ terms, i.e their large-$z$ expansions start with the $ z^{-n-1} $ term.
To remedy this we need to identify {\it associated functions of the third type}, which we do with the following definitions
\begin{align}
	P^{(3)}_{n}(z) & := \int_{S_1\times S_2} dxdy \frac{w(x,y)}{x+y}\frac{P_{n}(x)}{z-y} = -\int_{S_2} dy \frac{w_2(y)}{z-y}P^{(1)}_{n}(-y) ,
\\
	Q^{(3)}_{n}(z) & := \int_{S_1\times S_2} dxdy \frac{w(x,y)}{x+y}\frac{Q_{n}(y)}{z-x} = -\int_{S_1} dx \frac{w_1(x)}{z-x}Q^{(1)}_{n}(-x) .
\end{align}
We can then see that this set satisfies the above requirement because of
\begin{align}
	P^{(3)}_{n}(z) & \mathop{=}\limits_{z \to \infty} \sum_{l \geq n} z^{-l-1}\int_{S_1\times S_2} dxdy \frac{w(x,y)}{x+y}P_{n}(x)y^{l} = S_{n}^{-1}z^{-n-1} + {\rm O}(z^{-n-2}) ,
\label{AssocFn_P3}\\
	Q^{(3)}_{n}(z) & \mathop{=}\limits_{z \to \infty} \sum_{l \geq n} z^{-l-1}\int_{S_1\times S_2} dxdy \frac{w(x,y)}{x+y}x^{l}Q_{n}(y) = S_{n}^{-1}z^{-n-1} + {\rm O}(z^{-n-2}) ,	
\label{AssocFn_Q3}
\end{align}
by virtue of the orthogonality conditions $ \langle P_{n}, y^{l} \rangle = \langle x^{l}, Q_{n} \rangle = 0 $ for $ 0 \leq l \leq n-1 $.
Furthermore this third set has linear relations with the first two via the relations
\begin{equation}
	P^{(3)}_{n}(z) = - P^{(2)}_{n}(-z) - f_2(z)P^{(1)}_{n}(-z) ,
\qquad
	Q^{(3)}_{n}(z) = - Q^{(2)}_{n}(-z) - f_1(z)Q^{(1)}_{n}(-z) .
\end{equation}
  
\subsection{Christoffel-Darboux Formulae}\label{C-D_formulae}
In addition to the definition of the kernel \eqref{kernel_00} we give a full set of related sums, 
of which the first three will assume importance in subsequent developments, given here by
\begin{align}
	K^{0,1}_{n}(x,y) & := \sum_{l=0}^{n} P_l(x)Q^{(1)}_l(y) ,
\label{kernel_01}\\
	K^{1,0}_{n}(x,y) & := \sum_{l=0}^{n} P^{(1)}_l(x)Q_l(y) ,
\label{kernel_10}\\
	K^{1,1}_{n}(x,y) & := \sum_{l=0}^{n} P^{(1)}_l(x)Q^{(1)}_l(y) {+} \frac{1}{x+y},
\label{kernel_11}\\
	K^{0,2}_{n}(x,y) & := \sum_{l=0}^{n} P_l(x)Q^{(2)}_l(y) {-} \frac{1}{x+y},
\label{kernel_02}\\
	K^{1,2}_{n}(x,y) & := \sum_{l=0}^{n} P^{(1)}_l(x)Q^{(2)}_l(y) {-} \frac{f_1(x)}{x+y},
\label{kernel_12}\\
	K^{2,0}_{n}(x,y) & := \sum_{l=0}^{n} P^{(2)}_l(x)Q_l(y) {-} \frac{1}{x+y}.
\label{kernel_20}\\
	K^{2,1}_{n}(x,y) & := \sum_{l=0}^{n} P^{(2)}_l(x)Q^{(1)}_l(y) {-} \frac{f_2(y)}{x+y},
\label{kernel_21}\\
	K^{2,2}_{n}(x,y) & := \sum_{l=0}^{n} P^{(2)}_l(x)Q^{(2)}_l(y) {-} \frac{g_1(x)+g_2(y)}{x+y} .
\label{kernel_22}
\end{align}

However to begin with we are going re-derive the Christoffel-Darboux summation formula of \eqref{kernel_00} for a number of reasons elaborated on subsequently.
\begin{proposition}\label{CD_sum}
Let generic conditions apply.  
The Christoffel-Darboux sum for the reproducing kernel \eqref{kernel_00} has the evaluation
\begin{equation}
	(x+y)K^{0,0}_{n}(x,y) = \hat{P}_{n}(x)\check{Q}_{n}(y) + \frac{S_{n}}{S_{n+1}}\left[ P_{n}(x)Q_{n+1}(y)+P_{n+1}(x)Q_{n}(y) \right] .
\label{C-Dsum}
\end{equation}
\end{proposition}
\begin{proof}
Our proof of this evaluation is simple and can be reduced to seven steps, as follows:
\begin{align*}
	(x+y)\sum_{j=0}^{n} P_{j}(x)Q_{j}(y)
	& = \sum_{j=0}^{n}\sum_{l=0}^{j+1} X_{j,l}P_{l}(x)Q_{j}(y) + \sum_{j=0}^{n}\sum_{l=0}^{j+1} Y_{j,l}P_{j}(x)Q_{l}(y) , \tag*{(1)}
\\
	& = \sum_{j=0}^{n}\sum_{l=0}^{j+1} X_{j,l}P_{l}(x)Q_{j}(y) + \sum_{l=0}^{n}\sum_{j=0}^{l+1} Y_{l,j}P_{l}(x)Q_{j}(y) , \tag*{(2)}
\\
	& = \sum_{j=0}^{n}\sum_{l=0}^{j+1} X_{j,l}P_{l}(x)Q_{j}(y) + \sum_{j=0}^{n+1}\sum_{l=j-1}^{n} Y_{l,j}P_{l}(x)Q_{j}(y) , \tag*{(3)}
\\
	& = \sum_{j=0}^{n}\sum_{l=0}^{n+1} X_{j,l}P_{l}(x)Q_{j}(y) + \sum_{j=0}^{n+1}\sum_{l=0}^{n} Y_{l,j}P_{l}(x)Q_{j}(y) , \tag*{(4)}
\\
	& = \sum_{j=0}^{n}\sum_{l=0}^{n} \left( X_{j,l}+Y_{l,j} \right)P_{l}(x)Q_{j}(y) + \sum_{j=0}^{n} X_{j,n+1}P_{n+1}(x)Q_{j}(y) + \sum_{l=0}^{n} Y_{l,n+1}P_{l}(x)Q_{n+1}(y) , \tag*{(5)}
\\
	& = \sum_{j=0}^{n}\sum_{l=0}^{n} \pi_{j}\eta_{l}P_{l}(x)Q_{j}(y) + P_{n+1}(x)\sum_{j=0}^{n} X_{j,n+1}Q_{j}(y) + Q_{n+1}(y)\sum_{l=0}^{n} Y_{l,n+1}P_{l}(x) , \tag*{(6)}
\\
	& = \hat{P}_{n}(x)\check{Q}_{n}(y) + X_{n,n+1}P_{n+1}(x)Q_{n}(y) + Y_{n,n+1}P_{n}(x)Q_{n+1}(y) , \tag*{(7)}
\end{align*}
where the reasoning is:
\begin{enumerate}
	\item using the Hessenberg matrices \eqref{Hessenberg},
	\item exchange the summation labels $ j \leftrightarrow l $ in the second term,
	\item exchange the summation order in the second term,
	\item in the first term extend the upper terminal in the $l$-sum from $j+1$ to the maximum $n+1$ because $ X_{j,l}=0 $ if $l>j+1$, 
	      likewise in the second term extend the lower terminal in the $l$-sum from $j-1$ to the minimum $0$ because $ Y_{l,j}=0 $ if $l<j-1$,   
	\item combine the common summands in both terms - the remainder giving two single sums,
	\item use the Cauchy identity \eqref{Cauchy_Id} in the double sum,
	\item in each of the two single sums only one term contributes,
	\item and finally use the results \eqref{intertwineDN},\eqref{intertwineUP} and evaluations \eqref{rCFF:+2},\eqref{sCFF:+2}, 
\end{enumerate}
to arrive at \eqref{C-Dsum}.
Clearly the key ingredient is \eqref{Cauchy_Id}.
\end{proof}

\begin{remark}
An evaluation for the Christoffel-Darboux sum was given in Thm. (4.1), Eqs. (4.2,4.3) of \cite{BGS_2010}, which we can confirm is correct, 
and when re-written in our notations and further simplified would appear as
\begin{multline*}
	- \frac{\eta_{n}}{\eta_{n+1}}s_{n-1,2}Q_{n}(y)\hat{P}_{n+1}(x) + s_{n,2}Q_{n+1}(y)\hat{P}_{n-1}(x)
\\
	 - \frac{y}{\eta_{n+1}}Q_{n+1}(y)\hat{P}_{n}(x) + s_{n,1}Q_{n+1}(y)\hat{P}_{n}(x) + s_{n,2}Q_{n+2}(y)\hat{P}_{n}(x) . 
\end{multline*}
This expression is more complicated than \eqref{C-Dsum} and involves additional non-trivial data from the recurrence coefficients such as $s_{n,1}$.
Furthermore it is manifestly unsymmetrical in roles the two polynomials play $ \{x, P_{n}(x) \} \leftrightarrow \{y, Q_{n}(y) \} $ due to specific choices made in the derivation
and obscures the underlying symmetry of the kernels, as we will shortly see.
\end{remark}

\begin{corollary}
Let us define functions $ \hat{P}^{(1)}_{n}(x), \check{Q}^{(1)}_{n}(y) $ and  $ \hat{P}^{(2)}_{n}(x), \check{Q}^{(2)}_{n}(y) $ 
associated with the intertwined polynomials $ \hat{P}_{n}(x), \check{Q}_{n}(y) $ as per \eqref{AssocFn_P1},\eqref{AssocFn_Q1} and \eqref{AssocFn_P2},\eqref{AssocFn_Q2} respectively. 
The Christoffel-Darboux sums for the mixed kernels \eqref{kernel_01} -- \eqref{kernel_22} have the evaluations
\begin{equation}
	(x+y)K^{0,1}_{n}(x,y) = \hat{P}_{n}(x) \left( \check{Q}^{(1)}_{n}(y)-1 \right)
	 + \frac{S_{n}}{S_{n+1}}\left[ P_{n}(x)Q^{(1)}_{n+1}(y)+P_{n+1}(x)Q^{(1)}_{n}(y) \right] ,
\label{C-Dsum_01}
\end{equation}
\begin{equation}
	(x+y)K^{1,0}_{n}(x,y) = \left( \hat{P}^{(1)}_{n}(x)-1 \right)\check{Q}_{n}(y)
	+ \frac{S_{n}}{S_{n+1}}\left[ P^{(1)}_{n}(x)Q_{n+1}(y)+P^{(1)}_{n+1}(x)Q_{n}(y) \right] ,
\label{C-Dsum_10}
\end{equation}
\begin{equation}
	(x+y)K^{1,1}_{n}(x,y) = \left( \hat{P}^{(1)}_{n}(x)-1 \right) \left( \check{Q}^{(1)}_{n}(y)-1 \right)
	+ \frac{S_{n}}{S_{n+1}}\left[ P^{(1)}_{n}(x)Q^{(1)}_{n+1}(y)+P^{(1)}_{n+1}(x)Q^{(1)}_{n}(y) \right] ,
\label{C-Dsum_11}
\end{equation}
\begin{equation}
	(x+y)K^{0,2}_{n}(x,y) = \hat{P}_{n}(x)\check{Q}^{(2)}_{n}(y)
	+ \frac{S_{n}}{S_{n+1}}\left[ P_{n}(x)Q^{(2)}_{n+1}(y)+P_{n+1}(x)Q^{(2)}_{n}(y) \right] ,
\label{C-Dsum_02}
\end{equation}
\begin{equation}
	(x+y)K^{1,2}_{n}(x,y) = \left( \hat{P}^{(1)}_{n}(x)-1 \right) \check{Q}^{(2)}_{n}(y)
	+ \frac{S_{n}}{S_{n+1}}\left[ P^{(1)}_{n}(x)Q^{(2)}_{n+1}(y)+P^{(1)}_{n+1}(x)Q^{(2)}_{n}(y) \right] ,
\label{C-Dsum_12}
\end{equation}
\begin{equation}
	(x+y)K^{2,0}_{n}(x,y) = \hat{P}^{(2)}_{n}(x)\check{Q}_{n}(y)
	+ \frac{S_{n}}{S_{n+1}}\left[ P^{(2)}_{n}(x)Q_{n+1}(y)+P^{(2)}_{n+1}(x)Q_{n}(y) \right] ,
\label{C-Dsum_20}
\end{equation}
\begin{equation}
	(x+y)K^{2,1}_{n}(x,y) = \hat{P}^{(2)}_{n}(x) \left( \check{Q}^{(1)}_{n}(y)-1 \right)
	+ \frac{S_{n}}{S_{n+1}}\left[ P^{(2)}_{n}(x)Q^{(1)}_{n+1}(y)+P^{(2)}_{n+1}(x)Q^{(1)}_{n}(y) \right] ,
\label{C-Dsum_21}
\end{equation}
\begin{equation}
	(x+y)K^{2,2}_{n}(x,y) = \hat{P}^{(2)}_{n}(x)\check{Q}^{(2)}_{n}(y)
	+ \frac{S_{n}}{S_{n+1}}\left[ P^{(2)}_{n}(x)Q^{(2)}_{n+1}(y)+P^{(2)}_{n+1}(x)Q^{(2)}_{n}(y) \right] .
\label{C-Dsum_22}
\end{equation}
\end{corollary}
\begin{proof}
The calculation follows the method of proof in Prop. \ref{CD_sum}, but facilitated by the additional identities
\begin{alignat*}{3}
	x P^{(1)}_{n}(x) & =  \pi_{n} + \sum_{j=0}^{n+1} X_{n,j}P^{(1)}_{j}(x) , \qquad
	&& y Q^{(1)}_{n}(y) = \eta_{n} + \sum_{j=0}^{n+1} Y_{n,j}Q^{(1)}_{j}(y) ,
\\
	x P^{(2)}_{n}(x) & =  S^{-1}_{0}\delta_{n,0} + \sum_{j=0}^{n+1} X_{n,j}P^{(2)}_{j}(x) , \qquad
	&& y Q^{(2)}_{n}(y) = S^{-1}_{0}\delta_{n,0} + \sum_{j=0}^{n+1} Y_{n,j}Q^{(2)}_{j}(y) .
\end{alignat*}
Repeating those steps outlined earlier we have \eqref{C-Dsum_01}--\eqref{C-Dsum_22}.
\end{proof}

Having found this Christoffel-Darboux sum a number of remarkable implications follow, the first of which now unfolds.
\begin{proposition}
The polynomials $ P_{n}(x), Q_{n}(y) $ and the intertwined polynomials $ \hat{P}_{n}(x), \check{Q}_{n}(y) $ satisfy a pair of second order difference equations
\begin{align}
	(x-X_{n,n})P_{n} - \frac{S_{n}}{S_{n+1}}P_{n+1} + \frac{S_{n-1}}{S_{n}}P_{n-1} + \pi_{n}\hat{P}_{n-1} & = 0 ,
\label{PUp_recur}\\ 
	(y-Y_{n,n})Q_{n} - \frac{S_{n}}{S_{n+1}}Q_{n+1} + \frac{S_{n-1}}{S_{n}}Q_{n-1} + \eta_{n}\check{Q}_{n-1} & = 0 .
\label{QDown_recur}
\end{align}
\end{proposition} 
\begin{proof}
The Christoffel-Darboux sum \ref{C-Dsum} implies the above identities directly.
By differencing this formula and making the substitutions $ \hat{P}_{n} = \hat{P}_{n-1}-\eta_{n}P_{n} $, $ \check{Q}_{n} = \check{Q}_{n-1}-\pi_{n}Q_{n} $
we find after simplification that the result can be factorised in manner utilising $ P_{n}(x) $ and $ Q_{n}(y) $ thus
\begin{multline*}
	\pi_{n}\eta_{n} =
	\frac{1}{P_{n}(x)} \left[ xP_{n} - \frac{S_{n}}{S_{n+1}}P_{n+1} + \frac{S_{n-1}}{S_{n}}P_{n-1} + \pi_{n}\hat{P}_{n-1}  \right]
\\
	+ \frac{1}{Q_{n}(y)} \left[ yQ_{n} - \frac{S_{n}}{S_{n+1}}Q_{n+1} + \frac{S_{n-1}}{S_{n}}Q_{n-1} + \eta_{n}\check{Q}_{n-1}  \right] .
\end{multline*} 
Clearly the left-hand side is a constant independent of $ x $ and $ y $, whilst on the right-hand the first term is only dependent on $x$ and the second term only on $y$.
Therefore the two ratios are in fact constants independent of either $ x $ or $ y $. 
We can determine the two unknown constants by employing the Hessenberg formula \eqref{Hessenberg} and find they are $ X_{n,n} $ and $ Y_{n,n} $ for the $ P, Q $ equations respectively.
\end{proof}

\begin{remark}
From \eqref{PUp_recur}, \eqref{QDown_recur} we can confirm the relations 
\begin{gather}
	X_{n,n-1} = \pi_{n}\eta_{n-1} - \frac{S_{n-1}}{S_{n}} ,\qquad
	Y_{n,n-1} = \pi_{n-1}\eta_{n} - \frac{S_{n-1}}{S_{n}} ,
\\
	X_{n,k} = \pi_{n}\eta_{k} , \qquad Y_{n,k} = \pi_{k}\eta_{n}, \quad k\leq n-2 , 	
\end{gather}
first given in \eqref{XY-elements}.
Two particularly useful relations arising from \eqref{PUp_recur}, \eqref{QDown_recur} are
\begin{equation}
	X_{n,n} = \frac{\pi_{n}S_{n-1}}{\pi_{n-1}S_{n}} + \pi_{n}r_{n-1,1} , \qquad
	Y_{n,n} = \frac{\eta_{n}S_{n-1}}{\eta_{n-1}S_{n}} + \eta_{n}s_{n-1,1} ,
\end{equation}
which have appeared before in \eqref{X_from_r}.
\end{remark}

\begin{remark}
The pair of recurrences \eqref{PUp_recur}, \eqref{QDown_recur} represent a lowering of the order of the difference equations \eqref{P_recur}, \eqref{Q_recur} by one unit
through a summation, but of course at the expense of introducing other variables in the process.
Another way to view this pair is as the simplest relations to express $ \hat{P}_{n}(x), \check{Q}_{n}(y) $ in terms of the basic polynomial set of $ P_{n}(x), Q_{n}(y) $,
which are independent of the properties listed in Prop. \ref{PQintertwined}.
There are alternative forms of this relationship which we record here for future use
\begin{gather}
	\pi_{n}\hat{P}_{n}(x) = \frac{S_{n}}{S_{n+1}}P_{n+1}(x) - \left( Y_{n,n}+x \right)P_{n}(x) - \frac{S_{n-1}}{S_{n}}P_{n-1}(x) ,
\label{hatPsubs}
\\
	\eta_{n}\check{Q}_{n}(y) = \frac{S_{n}}{S_{n+1}}Q_{n+1}(y) - \left( X_{n,n}+y \right)Q_{n}(y) - \frac{S_{n-1}}{S_{n}}Q_{n-1}(y) ,
\label{checkQsubs}
\\
	\pi_{n}\left( \hat{P}^{(1)}_{n}(x)-1 \right) = \frac{S_{n}}{S_{n+1}}P^{(1)}_{n+1}(x) - \left( Y_{n,n}+x \right)P^{(1)}_{n}(x) - \frac{S_{n-1}}{S_{n}}P^{(1)}_{n-1}(x) ,
\label{hatP1subs}
\\
	\eta_{n}\left( \check{Q}^{(1)}_{n}(y)-1 \right) = \frac{S_{n}}{S_{n+1}}Q^{(1)}_{n+1}(y) - \left( X_{n,n}+y \right)Q^{(1)}_{n}(y) - \frac{S_{n-1}}{S_{n}}Q^{(1)}_{n-1}(y) ,
\label{checkQ1subs}
\\
	\pi_{n}\hat{P}^{(2)}_{n}(x)-S_{0}^{-1}\delta_{n,0} = \frac{S_{n}}{S_{n+1}}P^{(2)}_{n+1}(x) - \left( Y_{n,n}+x \right)P^{(2)}_{n}(x) - \frac{S_{n-1}}{S_{n}}P^{(2)}_{n-1}(x) ,
\label{hatP2subs}
\\
	\eta_{n}\check{Q}^{(2)}_{n}(y)-S_{0}^{-1}\delta_{n,0} = \frac{S_{n}}{S_{n+1}}Q^{(2)}_{n+1}(y) - \left( X_{n,n}+y \right)Q^{(2)}_{n}(y) - \frac{S_{n-1}}{S_{n}}Q^{(2)}_{n-1}(y) .
\label{checkQ2subs}
\end{gather} 
We remind the reader that the coefficients $ X_{n,n} $ and $ Y_{n,n} $ are directly related to the recurrence coefficients $ r_{n-1,1} $ and $ s_{n-1,1} $,
as well as the norms and the coefficients $ \pi_{j}, \eta_{k} $. In addition $ X_{n,n}+Y_{n,n}=\pi_{n}\eta_{n} $.
\end{remark}

The evaluations given of the Christoffel-Darboux sums also imply further interesting identities, which we refer to as anti-incidence identities,
and involve both the $P,Q$ sides of the bi-orthogonal system.
Hitherto all of our results have involved relations of the $P$ side amongst themselves, and parallel relations solely involving the $Q$ side, 
without any mutual coupling. 
\begin{corollary}\label{anti-incidence}
For all $x$, $n\geq 0$, the bi-orthogonal polynomials satisfy the coupled bi-linear relations
\begin{equation}
	\hat{P}_{n}(x)\check{Q}_{n}(-x) + \frac{S_{n}}{S_{n+1}}\left[ P_{n}(x)Q_{n+1}(-x)+P_{n+1}(x)Q_{n}(-x) \right] = 0,
\label{anti-CDsum-00}
\end{equation} 
whilst their associated functions satisfy the coupled bi-linear relations
\begin{align}
	\hat{P}_{n}(x)\check{Q}^{(1)}_{n}(-x) + \frac{S_{n}}{S_{n+1}}\left[ P_{n}(x)Q^{(1)}_{n+1}(-x)+P_{n+1}(x)Q^{(1)}_{n}(-x) \right] & = \hat{P}_{n}(x) ,
\label{anti-CDsum-01}\\
	\hat{P}^{(1)}_{n}(x)\check{Q}_{n}(-x) + \frac{S_{n}}{S_{n+1}}\left[ P^{(1)}_{n}(x)Q_{n+1}(-x)+P^{(1)}_{n+1}(x)Q_{n}(-x) \right] & = \check{Q}_{n}(-x) ,
\label{anti-CDsum-10}\\
	\hat{P}^{(1)}_{n}(x)\check{Q}^{(1)}_{n}(-x) + \frac{S_{n}}{S_{n+1}}\left[ P^{(1)}_{n}(x)Q^{(1)}_{n+1}(-x)+P^{(1)}_{n+1}(x)Q^{(1)}_{n}(-x) \right] 
	& = \hat{P}^{(1)}_{n}(x) + \check{Q}^{(1)}_{n}(-x) ,
\label{anti-CDsum-11}\\
	\hat{P}_{n}(x)\check{Q}^{(2)}_{n}(-x) + \frac{S_{n}}{S_{n+1}}\left[ P_{n}(x)Q^{(2)}_{n+1}(-x)+P_{n+1}(x)Q^{(2)}_{n}(-x) \right] & = -1 ,
\label{anti-CDsum-02}\\
	\hat{P}^{(1)}_{n}(x)\check{Q}^{(2)}_{n}(-x) + \frac{S_{n}}{S_{n+1}}\left[ P^{(1)}_{n}(x)Q^{(2)}_{n+1}(-x)+P^{(1)}_{n+1}(x)Q^{(2)}_{n}(-x) \right] & = \check{Q}^{(2)}_{n}(-x)-f_1(x) ,
\label{anti-CDsum-12}\\
	\hat{P}^{(2)}_{n}(x)\check{Q}_{n}(-x) + \frac{S_{n}}{S_{n+1}}\left[ P^{(2)}_{n}(x)Q_{n+1}(-x)+P^{(2)}_{n+1}(x)Q_{n}(-x) \right] & = -1 ,
\label{anti-CDsum-20}\\
	\hat{P}^{(2)}_{n}(x)\check{Q}^{(1)}_{n}(-x) + \frac{S_{n}}{S_{n+1}}\left[ P^{(2)}_{n}(x)Q^{(1)}_{n+1}(-x)+P^{(2)}_{n+1}(x)Q^{(1)}_{n}(-x) \right] & = \hat{P}^{(2)}_{n}(x)-f_2(-x) ,
\label{anti-CDsum-21}\\
	\hat{P}^{(2)}_{n}(x)\check{Q}^{(2)}_{n}(-x) + \frac{S_{n}}{S_{n+1}}\left[ P^{(2)}_{n}(x)Q^{(2)}_{n+1}(-x)+P^{(2)}_{n+1}(x)Q^{(2)}_{n}(-x) \right] 
	& = f_1(x)f_2(-x) .
\label{anti-CDsum-22}
\end{align}
\end{corollary}
\begin{proof}
These relations follow immediately from \eqref{C-Dsum}, \eqref{C-Dsum_01}, \eqref{C-Dsum_10} and \eqref{C-Dsum_11}. 
\end{proof}

Taken together all the preceding developments reveal a fundamental structure for the kernels.
Let us define the column vectors $ \mathrm{P}_{n}^{(\mu)} $, $ \mu=0,1,2 $ of $ \mathrm{P}_{n} = \left( \mathrm{P}_{n}^{(0)},\mathrm{P}_{n}^{(1)},\mathrm{P}_{n}^{(2)} \right) $ and similarly for $ \mathrm{Q}_{n} $.
We construct forms of the kernels whereby the intertwined polynomials have been eliminated.
\begin{proposition}
Assume generic conditions hold.
The kernels $ K^{\mu,\nu}_{n}(x,y) $ for $ \mu, \nu = 0,1,2 $ have the bilinear forms
\begin{equation}
	\pi_{n}\eta_{n}(x+y) K^{\mu,\nu}_{n}(x,y) = \mathrm{P}_{n}^{(\mu)}(x)^{T} G_{n}(x,y) \mathrm{Q}_{n}^{(\nu)}(y) ,
\label{symmetric_kernel}
\end{equation}
where $ G(x,y) $ is the $ 3\times 3 $ matrix
\begin{equation}
	G_{n}(x,y) =     
	\begin{pmatrix}
	\frac{S_{n}^2}{S_{n+1}^2}						& \frac{S_{n}}{S_{n+1}}\left[Y_{n,n}-y\right]	& 	-\frac{S_{n-1}}{S_{n+1}} 					\cr
	\frac{S_{n}}{S_{n+1}}\left[X_{n,n}-x\right]		& \left[X_{n,n}+y\right]\left[Y_{n,n}+x\right]	&  \frac{S_{n-1}}{S_{n}}\left[Y_{n,n}+x\right]	\cr
	-\frac{S_{n-1}}{S_{n+1}}						& \frac{S_{n-1}}{S_{n}}\left[X_{n,n}+y\right]	&  \frac{S_{n-1}^2}{S_{n}^2}					\cr
	\end{pmatrix} .
\label{Gmatrix}
\end{equation}
The last column of $ G $ is proportional to the components of $ \hat{P}_{n}(x) $ and the last row is proportional to the components of $ \check{Q}_{n}(y) $.
The characteristic polynomial of $G$ gives the cubic spectral curve
\begin{multline}
	\det(G_{n}(x,y) - \mathbb{I}_3\lambda) = -\lambda^3 + \left( \frac{S_{n}^2}{S_{n+1}^2}+\frac{S_{n-1}^2}{S_{n}^2}+\left[Y_{n,n}+x\right]\left[X_{n,n}+y\right] \right)\lambda^2
\\
	- \pi_{n}\eta_{n} \frac{S_{n}^2}{S_{n+1}^2}(x+y)\lambda - \pi_{n}^2\eta_{n}^2\frac{S_{n-1}^2}{S_{n+1}^2} . 
\label{cubic}
\end{multline}
The anti-incidence formulae in Cor. \ref{anti-incidence} are then expressions of orthogonality between 
$ \mathrm{P}_{n}^{(0)}(x), \mathrm{P}_{n}^{(1)}(x) $ on the one hand and $ \mathrm{Q}_{n}^{(0)}(-x), \mathrm{Q}_{n}^{(1)}(-x) $ on the other. 
Furthermore the $G_{n}$ matrix is invertible if $ n\geq 1 $
\begin{equation}
	\pi_{n}\eta_{n}G_{n}(x,y)^{-1} =     
	\begin{pmatrix}
		0												& \frac{S_{n+1}}{S_{n}}	&  -\frac{S_{n+1}}{S_{n-1}}\left[Y_{n,n}+x\right] 	\cr
		\frac{S_{n+1}}{S_{n}}							& 0						&   \frac{S_{n}}{S_{n-1}}							\cr
		-\frac{S_{n+1}}{S_{n-1}}\left[X_{n,n}+y\right]	& \frac{S_{n}}{S_{n-1}}	&  -\frac{S_{n}^2}{S_{n-1}^2}(x+y)					\cr
	\end{pmatrix} .
\label{inverse_Gmatrix}
\end{equation}
\end{proposition}

\subsection{Spectral Derivatives}\label{SpectralStructures}

Although we will address the question of derivatives of our bi-orthogonal system of polynomials $ P_{n}(x), Q_{n}(y) $ 
with respect to $x,y$ in \S \ref{Cauchy-Laguerre-BOPS} in systematic detail,
we need to treat this problem in the general setting, i.e for more general classes of weights whereby all moments exist and bi-moment determinants are non-vanishing.
In fact it is more convenient and transparent to treat the generic situation.
 
For general classes of weights the bi-orthogonal polynomials possess derivatives with respect to the spectral variables which can be parametrised as
\begin{align}
	W_1(x) \partial_{x} P_{n}(x) & = \Theta^{+}_{n}(x) P_{n+1}(x) + \Omega_{n}(x) P_{n}(x) + \Theta^{-}_{n}(x) P_{n-1}(x) ,
\label{SpectralDiff_P}\\
	W_2(y) \partial_{y} Q_{n}(y) & = \Delta^{+}_{n}(y) Q_{n+1}(y) + \Gamma_{n}(y) Q_{n}(y) + \Delta^{-}_{n}(y) Q_{n-1}(y) ,
\label{SpectralDiff_Q}
\end{align}
where $ W_{1,2} $, $ \Theta^{\pm}_{n} $, $ \Delta^{\pm}_{n} $, $ \Omega_{n} $ and $ \Gamma_{n} $ are polynomials in either $ x,y $ as the case may be and their degrees are fixed, 
i.e. independent of $n$. 
That this is true follows from the expansion of arbitrary monomials into either of the bi-orthogonal $P,Q$ bases, 
and the fact that the polynomials satisfy third-order or four-term recurrence relations so that any sequence of polynomials can be folded into a set of three with consecutive degrees.
In addition it is straight forward to prove this statement in an explicit and constructive manner by adapting the techniques found in Bauldry \cite{Bau_1990}, 
and Bonan and Clark \cite{BC_1990}, however we hesitate to do this here due to the fact that such an exercise would distract us from our primary task.

Utilising the recurrence relations one can then show that the matrix formulations of \eqref{SpectralDiff_P}, \eqref{SpectralDiff_Q} satisfy the relations
\begin{equation}
	W_1(x)\partial_{x} 
	\begin{pmatrix}
		P_{n+1} \\ P_{n} \\ P_{n-1}
	\end{pmatrix}
	= A_{n}(x)
	\begin{pmatrix}
		P_{n+1} \\ P_{n} \\ P_{n-1}
	\end{pmatrix} ,
\label{Dx_P-vector}
\end{equation}
and
\begin{equation}
	W_2(y)\partial_{y} 
	\begin{pmatrix}
		Q_{n+1} \\ Q_{n} \\ Q_{n-1}
	\end{pmatrix}
	= D_{n}(y)
	\begin{pmatrix}
		Q_{n+1} \\ Q_{n} \\ Q_{n-1}
	\end{pmatrix} ,
\label{Dy_Q-vector}
\end{equation}
where the first Lax matrix $ A_{n} $ is
\begin{multline}\label{P_Diff}
	A_{n}(x) =
\\
	\begin{pmatrix}
		\Omega_{n+1}+\frac{1}{r_{n,2}}\left(\frac{x}{\pi_{n+1}}-r_{n,1} \right)\Theta^{+}_{n+1} & 
		\Theta^{-}_{n+1}-\frac{1}{r_{n,2}}\left(\frac{x}{\pi_{n}}+r_{n,0} \right)\Theta^{+}_{n+1} &
		-\frac{r_{n,-1}}{r_{n,2}}\Theta^{+}_{n+1} \\
		\Theta^{+}_{n} & \Omega_{n} & \Theta^{-}_{n} \\
		-\frac{r_{n-1,2}}{r_{n-1,-1}}\Theta^{-}_{n-1}  &
		\Theta^{+}_{n-1}+\frac{1}{r_{n-1,-1}}\left(\frac{x}{\pi_{n}}-r_{n-1,1} \right)\Theta^{-}_{n-1} &
		\Omega_{n-1}-\frac{1}{r_{n-1,-1}}\left(\frac{x}{\pi_{n-1}}+r_{n-1,0} \right)\Theta^{-}_{n-1}
	\end{pmatrix} ,
\end{multline}
and the second Lax matrix $ D_{n} $ is
\begin{multline}\label{Q_Diff}
	D_{n}(y) =
\\
	\begin{pmatrix}
		\Gamma_{n+1}+\frac{1}{s_{n,2}}\left(\frac{y}{\eta_{n+1}}-s_{n,1} \right)\Delta^{+}_{n+1} & 
		\Delta^{-}_{n+1}-\frac{1}{s_{n,2}}\left(\frac{y}{\eta_{n}}+s_{n,0} \right)\Delta^{+}_{n+1} &
		-\frac{s_{n,-1}}{s_{n,2}}\Delta^{+}_{n+1} \\
		\Delta^{+}_{n} & \Gamma_{n} & \Delta^{-}_{n} \\
		-\frac{s_{n-1,2}}{s_{n-1,-1}}\Delta^{-}_{n-1}  &
		\Delta^{+}_{n-1}+\frac{1}{s_{n-1,-1}}\left(\frac{y}{\eta_{n}}-s_{n-1,1} \right)\Delta^{-}_{n-1} &
		\Gamma_{n-1}-\frac{1}{s_{n-1,-1}}\left(\frac{y}{\eta_{n-1}}+s_{n-1,0} \right)\Delta^{-}_{n-1}
	\end{pmatrix} .
\end{multline}
However the nine polynomials appearing in each of these matrices satisfy inter-relations, allowing us to reduce the number of independent entries.
These relations follow by performing exact summation of the compatibility relations
or by employing explicit representations for the polynomials constructed using the techniques of \cite{Bau_1990}, \cite{BC_1990}.
To this end we define the matrix invariants $ \Lambda_l(A) $, $ l\in\N $ of an arbitrary square matrix $A$ from the generating function expanded about $t=0$, 
$ \det({\mathbb{1}}-t A) = 1 - \Lambda_1(A) t + \Lambda_2(A) t^2 - \Lambda_3(A) t^3 + \ldots $.
In our case of $ 3\times 3$ matrices the expansion terminates after the cubic term.
\begin{proposition}\label{Spectral_Invariants}
The Lax matrices in \eqref{SpectralDiff_P}, \eqref{SpectralDiff_Q} satisfy the following relations:\\
The first invariant or trace conditions are
\begin{equation}\label{trace-A}
	\Lambda_1(A_{n}) =
	\Omega_{n+1}+\Omega_{n}+\Omega_{n-1}
	+\frac{1}{r_{n,2}}\left(\frac{x}{\pi_{n+1}}-r_{n,1} \right)\Theta^{+}_{n+1}-\frac{1}{r_{n-1,-1}}\left(\frac{x}{\pi_{n-1}}+r_{n-1,0} \right)\Theta^{-}_{n-1} = V_1 ,
\end{equation}
and
\begin{equation}\label{trace-D}
	\Lambda_1(D_{n}) =
	\Gamma_{n+1}+\Gamma_{n}+\Gamma_{n-1}
	+\frac{1}{s_{n,2}}\left(\frac{y}{\eta_{n+1}}-s_{n,1} \right)\Delta^{+}_{n+1}-\frac{1}{s_{n-1,-1}}\left(\frac{y}{\eta_{n-1}}+s_{n-1,0} \right)\Delta^{-}_{n-1} = V_2 .
\end{equation}
Note that the right-hand sides of these latter two equations are independent of $n$.\\
The second invariants are
\begin{equation}\label{2ndInvariant-A}
	\Lambda_2(A_{n+1}) = \Lambda_2(A_{n})
	-\frac{W_1}{r_{n,2}\pi_{n+1}}\left[ \Theta^{+}_{n+1}+\frac{S_{n+1}S_{n}}{S_{n+2}S_{n-1}}\Theta^{-}_{n} \right] ,
\end{equation}
and
\begin{equation}\label{2ndInvariant-D}
	\Lambda_2(D_{n+1}) = \Lambda_2(D_{n})
	-\frac{W_2}{s_{n,2}\eta_{n+1}}\left[ \Delta^{+}_{n+1}+\frac{S_{n+1}S_{n}}{S_{n+2}S_{n-1}}\Delta^{-}_{n} \right] .
\end{equation}
\\
The third invariant or determinant relations are
\begin{multline}\label{3rdInvariant-A}
	\Lambda_3(A_{n+1}) = \Lambda_3(A_{n})
	-\frac{W_1}{r_{n,2}\pi_{n+1}}
		\Bigg\{ 
			 \Theta^{+}_{n+1}\Omega_{n} + \frac{S_{n+1}S_{n}}{S_{n+2}S_{n-1}}\Theta^{-}_{n}\Omega_{n+1} + \frac{\pi_{n+1}}{\pi_{n}}\Theta^{-}_{n+1}\Theta^{+}_{n}
		\\
				- \left[ \frac{1}{r_{n,-1}}\left( \frac{x}{\pi_{n}}+r_{n,0} \right)
						-\frac{S_{n+1}S_{n}}{S_{n+2}S_{n-1}}\frac{1}{r_{n,2}}\left( \frac{x}{\pi_{n+1}}-r_{n,1} \right)
						-\frac{\pi_{n}S_{n+1}S_{n}}{\pi_{n+1}S_{n+2}S_{n-1}} 
				  \right]\Theta^{+}_{n+1}\Theta^{-}_{n}
		\Bigg\} ,
\end{multline}
and
\begin{multline}\label{3rdInvariant-D}
	\Lambda_3(D_{n+1}) = \Lambda_3(D_{n})
	-\frac{W_2}{s_{n,2}\eta_{n+1}}
		\Bigg\{ 
			\Delta^{+}_{n+1}\Gamma_{n} + \frac{S_{n+1}S_{n}}{S_{n+2}S_{n-1}}\Delta^{-}_{n}\Gamma_{n+1} + \frac{\eta_{n+1}}{\eta_{n}}\Delta^{-}_{n+1}\Delta^{+}_{n}
\\
			- \left[ \frac{1}{s_{n,-1}}\left( \frac{y}{\eta_{n}}+s_{n,0} \right)
					-\frac{S_{n+1}S_{n}}{S_{n+2}S_{n-1}}\frac{1}{s_{n,2}}\left( \frac{y}{\eta_{n+1}}-s_{n,1} \right)
					-\frac{\eta_{n}S_{n+1}S_{n}}{\eta_{n+1}S_{n+2}S_{n-1}} 
			  \right]\Delta^{+}_{n+1}\Delta^{-}_{n}
		\Bigg\} .
\end{multline}
\end{proposition}
\begin{proof}
We only give details for the polynomial $P$ case as the $Q$ case is virtually identical after making a transcriptions of variables.
Let us define the abbreviations
\begin{equation*}
	\omega_{n} = \pi_{n}\frac{S_{n+1}}{S_{n-1}} ,  \qquad
	\zeta_{n} = \frac{1}{r_{n,2}}\left( \frac{x}{\pi_{n+1}}-r_{n,1} \right) , \qquad
	\gamma_{n} = \frac{1}{r_{n,-1}}\left( \frac{x}{\pi_{n}}+r_{n,0} \right) .
\end{equation*} 
The compatibility of the $P$-recurrence system \eqref{PQ_recur} and the spectral differential system \eqref{P_Diff} is expressed by the matrix equation
\begin{equation}
	A_{n+1} = W_1 \partial_{x}K_{n}\cdot K_{n}^{-1} + K_{n}\cdot A_{n}\cdot K_{n}^{-1} ,
\end{equation}
or in a scalar form between \eqref{SpectralDiff_P} and \eqref{P_recur}.
Choosing initially the scalar form we find three relations:
\begin{equation*}
	-\Omega_{n+2}+\Omega_{n-1}-\zeta_{n+1} \Theta^{+}_{n+2}+\gamma_{n} \Theta^{-}_{n}+\zeta_{n} \Theta^{+}_{n+1}-\gamma_{n-1} \Theta^{-}_{n-1} = 0,
\end{equation*}
\begin{equation*}
	\zeta_{n} \left(\Omega_{n+2}-\Omega_{n+1}\right)+\left(\zeta_{n} \zeta_{n+1}-\gamma_{n+1}\right)\Theta^{+}_{n+2}-\zeta_{n}^2\Theta^{+}_{n+1}
	+\frac{\omega_{n+1}}{\omega_{n}}\gamma_{n} \Theta^{+}_{n}+\Theta^{-}_{n+2}-\frac{\omega_{n+1}\omega_{n-1}}{\omega_{n}^2}\Theta^{-}_{n-1} = \frac{W_1}{r_{n,2}\pi_{n+1}},
\end{equation*}
and
\begin{equation*}
	\gamma_{n} \left(\Omega_{n+2}-\Omega_{n}\right)+\left(\gamma_{n} \zeta_{n+1}+\frac{\omega_{n+2} \omega_{n}}{\omega_{n+1}^2}\right)\Theta^{+}_{n+2}
	-\gamma_{n} \zeta_{n} \Theta^{+}_{n+1}-\Theta^{+}_{n-1}+\frac{\omega_{n}}{\omega_{n+1}}\zeta_{n}\Theta^{-}_{n+1}-\frac{\omega_{n-1}}{\omega_{n}}\zeta_{n-1}\Theta^{-}_{n-1} =
	 \frac{W_1}{r_{n,-1}\pi_{n}} .
\end{equation*}
The first relation is an obvious perfect difference and implies the trace constancy, \eqref{trace-A}. 
The matrix form of the compatibility gives only three new relations (the other six being trivial) along the first row, 
of which the $(1,1)$ component is the trace relation and the $(1,2)$, $(1,3)$ components are linear combinations of the first and second, 
and of the first and third relations given above respectively.
The $(1,1)$ component is the only one which allows an exact summation of a linear combination of the spectral polynomials -
the others can only be constructed out of quadratic and cubic combinations of these polynomials 

To construct these non-linear invariants it is more transparent to use the matrix form.
We employ identities for the invariants of sums of matrices \cite{Ami_1979-80}, \cite{RS_1987}, which in our context of a sum of two arbitrary $ 3\times 3 $ matrices $A,B$ are
\begin{gather*}
	\Lambda_1(A+B) = \Lambda_1(A)+\Lambda_1(B) ,
\\
	\Lambda_2(A+B) = \Lambda_2(A)+\Lambda_2(B) + \Lambda_1(A)\Lambda_1(B) -\Lambda_1(A\cdot B) ,
\end{gather*}
\begin{multline*}
	\Lambda_3(A+B) = \Lambda_3(A)+\Lambda_3(B)
	- \Lambda_1(A)\Lambda_1(A\cdot B) - \Lambda_1(B)\Lambda_1(A\cdot B)
\\
	+ \Lambda_1(A)\Lambda_2(B) + \Lambda_1(B)\Lambda_2(A) + \Lambda_1(A\cdot A\cdot B) + \Lambda_1(A\cdot B\cdot B) .
\end{multline*}
The evaluations of the right-hand sides is simplified by the fact that $\partial_{x}K_{n}\cdot K_{n}^{-1}$ is traceless, nilpotent and also has vanishing second and third invariants.
\end{proof}

The matrix formulation \eqref{Dx_P-vector}, \eqref{Dy_Q-vector} only contains the polynomial components and we need to incorporate the associated functions,
which will entail a $x,y$-dependent yet $n$-independent gauge transformation of the $ 3\times 3 $ matrix variables \eqref{3x3_PQ}.
Thus we define new matrix variables by inserting univariate weight factors in the second and third columns thus
\begin{gather}
	\mathcal{P}_{n}(x) := 	\left(
	\begin{array}{ccc}
		P_{n+1}(x) 	& \frac{1}{\displaystyle w_1(x)}P^{(1)}_{n+1}(x) 	& \frac{1}{\displaystyle w_1(x)w_2(-x)}P^{(2)}_{n+1}(x) \cr
		P_{n}(x) 	& \frac{1}{\displaystyle w_1(x)}P^{(1)}_{n}(x) 		& \frac{1}{\displaystyle w_1(x)w_2(-x)}P^{(2)}_{n}(x) \cr
		P_{n-1}(x) 	& \frac{1}{\displaystyle w_1(x)}P^{(1)}_{n-1}(x) 	& \frac{1}{\displaystyle w_1(x)w_2(-x)}P^{(2)}_{n-1}(x) \cr
	\end{array}
			\right) = \left( \mathcal{P}_{n}^{(0)}(x),\mathcal{P}_{n}^{(1)}(x),\mathcal{P}_{n}^{(2)}(x) \right), 
\label{new_3x3:P}\\
	\mathcal{Q}_{n}(y) := 	\left(
	\begin{array}{ccc}
		Q_{n+1}(y) 	& \frac{1}{\displaystyle w_2(y)}Q^{(1)}_{n+1}(y) 	& \frac{1}{\displaystyle w_1(-y)w_2(y)}Q^{(2)}_{n+1}(y) \cr
		Q_{n}(y) 	& \frac{1}{\displaystyle w_2(y)}Q^{(1)}_{n}(y) 		& \frac{1}{\displaystyle w_1(-y)w_2(y)}Q^{(2)}_{n}(y) \cr
		Q_{n-1}(y) 	& \frac{1}{\displaystyle w_2(y)}Q^{(1)}_{n-1}(y) 	& \frac{1}{\displaystyle w_1(-y)w_2(y)}Q^{(2)}_{n-1}(y) \cr
	\end{array}
			\right) = \left( \mathcal{Q}_{n}^{(0)}(y),\mathcal{Q}_{n}^{(1)}(y),\mathcal{Q}_{n}^{(2)}(y) \right).
\label{new_3x3:Q}
\end{gather}

Our last task is to compute the confluence of the kernels at anti-incidence points for several reasons.
The first is that we will see subsequently, in \eqref{sigma-form}, 
the $\sigma$-functions of the deformed Laguerre densities will be anti-incidence confluences of the two mixed kernels $ K^{0,1}_{n} $ and $ K^{1,0}_{n} $.
The second reason is that they furnish an additional identity linking the Lax matrices $ A_{n} $ and $ D_{n} $ pairwise, 
and thus linking the bi-orthogonal pair $P$ and $Q$.
\begin{proposition}
The confluent kernels $ K^{\mu,\nu}_{n}(x,-x) $ for $ \mu, \nu = 0,1 $ at anti-incidence are the bilinear forms
\begin{align}
	\pi_{n}\eta_{n}K^{0,0}_{n}(x,-x) & =
	\begin{cases}
		\mathcal{P}_{n}^{(0)}(x)^{T} \left( \frac{1}{W_1(x)}A_{n}(x)^{T} G_{n}(x,-x)+\partial_{x} G_{n}(x,-x) \right) \mathcal{Q}_{n}^{(0)}(-x)  ,& \\
		\mathcal{P}_{n}^{(0)}(x)^{T} \left( \frac{1}{W_2(-x)} G_{n}(x,-x)D_{n}(-x)+\partial_{y} G_{n}(x,-x) \right) \mathcal{Q}_{n}^{(0)}(-x)  ,&
	\end{cases}
\label{anti-incident_K00}
\end{align}
\begin{align}
	\pi_{n}\eta_{n}K^{1,0}_{n}(x,-x) & = w_1(x)
	\begin{cases}
		\mathcal{P}_{n}^{(1)}(x)^{T} \left( \frac{1}{W_1(x)}A_{n}(x)^{T} G_{n}(x,-x)+\partial_{x} G_{n}(x,-x) \right) \mathcal{Q}_{n}^{(0)}(-x)  ,& \\
		\mathcal{P}_{n}^{(1)}(x)^{T} \left( \frac{1}{W_2(-x)} G_{n}(x,-x)D_{n}(-x)+\partial_{y} G_{n}(x,-x) \right) \mathcal{Q}_{n}^{(0)}(-x)  ,&
	\end{cases}
\label{anti-incident_K10}
\end{align}
\begin{align}
	\pi_{n}\eta_{n}K^{0,1}_{n}(x,-x) & = w_2(-x)
	\begin{cases}
		\mathcal{P}_{n}^{(0)}(x)^{T} \left( \frac{1}{W_1(x)}A_{n}(x)^{T} G_{n}(x,-x)+\partial_{x} G_{n}(x,-x) \right) \mathcal{Q}_{n}^{(1)}(-x)  ,	& \\
		\mathcal{P}_{n}^{(0)}(x)^{T} \left( \frac{1}{W_2(-x)} G_{n}(x,-x)D_{n}(-x)+\partial_{y} G_{n}(x,-x) \right) \mathcal{Q}_{n}^{(1)}(-x)  .&
	\end{cases}
\label{anti-incident_K01}
\end{align}
\end{proposition}
An important consequence of the above relations is that the spectral matrix for the $ Q $ polynomials is fully determined by that of the $ P $ polynomials and vice-versa.
\begin{corollary}
Assume generic conditions hold. In addition assume $ W_1(x) \neq 0 $, $ W_2(x) \neq 0 $ and $ n\geq 1$.
The spectral matrices $ A_{n} $, $ D_{n} $ are related by the identities
\begin{multline}
	\frac{1}{W_2(-x)}D_{n}(-x) = \frac{1}{W_1(x)}G_{n}(x,-x)^{-1}A_{n}(x)^{T} G_{n}(x,-x)
\\
	+ \frac{1}{3} \left[ \frac{1}{W_2(-x)}{\rm Tr}D_{n}(-x) - \frac{1}{W_1(x)}{\rm Tr}A_{n}(x) \right]{\mathbb{1}}_{3}
\\
	+ \frac{1}{\pi_{n}\eta_{n}}
	\begin{pmatrix}
	-1								& \frac{S_{n+1}}{S_{n}}\left(X_{n,n}-x\right)	&  \frac{S_{n+1}S_{n-1}}{S_{n}^2} 	\cr
	 0								& 0												&  0								\cr
	-\frac{S_{n}^2}{S_{n+1}S_{n-1}}	& -\frac{S_{n}}{S_{n-1}}\left(Y_{n,n}+x\right)	&  1								\cr
	\end{pmatrix} .
\label{D-A-relation}
\end{multline} 
\end{corollary}
\begin{proof}
This follows from equating the pair of evaluations \eqref{anti-incident_K00}--\eqref{anti-incident_K01} for arbitrary $ \mathcal{P}_{n}^{(0)}(x), \mathcal{Q}_{n}^{(0)}(-x) $.
However equality of the matrices holds only modulo a term proportional to $ G_{n}(x,-x) $ or
\begin{equation*}
	\frac{1}{W_2(-x)}D_{n}(-x) 
	= \frac{1}{W_1(x)}G_{n}(x,-x)^{-1}A_{n}(x)^{T} G_{n}(x,-x) + G_{n}(x,-x)^{-1}\left( \partial_{x} G_{n}(x,-x)-\partial_{y} G_{n}(x,-x) \right) + m_{n}(x) {\mathbb{1}}_{3} .
\end{equation*}
The proportionality factor $ m_{n}(x) $ is most easily found using the trace of this relation.
\end{proof}
\begin{remark}
In practice $ W_2(-x) $ is equal to $ W_1(x) $ modulo a sign factor, therefore simplifying \eqref{D-A-relation} further
and the traces are also simple polynomials due to \eqref{trace-A} and \eqref{trace-D}, as will be seen subsequently.
\end{remark}

\subsection{The Undeformed Cauchy-Laguerre Weight and Boundary Condition Data}\label{Undeformed_Case}

Putting $ \xi, \psi = 0 $, or setting $ s, t \to 0^{+} $ and extracting the trivial factors of $ (1-\xi) $, $ (1-\psi) $ from the partition function, 
we turn off the deformation and are back to the undeformed Cauchy-Laguerre system.
We examine this in detail here for three reasons - firstly it serves as a check of the foregoing theory, 
secondly it provides some insight into these results for our subsequent development,
and lastly it furnishes us with some of the boundary condition data as $ s,t \to 0^{+} $ which we require for solving for the differential equations in $s,t$.

The spectral data consists of
$ W_1(x) = x $, $ W_2(y) = y $
and the moment data
\begin{gather}
	\alpha_{j} = \Gamma(a+1+j) ,
\qquad
	\beta_{k} = \Gamma(b+1+k) ,
\\
	\pi_{n} = S_{n} \frac{n!\Gamma(a+1+n)\Gamma(a+b+1+n)}{\Gamma(a+b+1+2n)} ,
\qquad
	\eta_{n} = S_{n} \frac{n!\Gamma(b+1+n)\Gamma(a+b+1+n)}{\Gamma(a+b+1+2n)} ,
\end{gather}
along with $ \pi_{n}\eta_{n} = a+b+1+2n $. This latter result is deduced from
\begin{equation}
	\frac{S_{n+1}^2}{S_{n}^2}
	= \frac{(2n+a+b+3)(2n+a+b+2)^2(2n+a+b+1)}{(n+1)^2(n+a+b+1)^2(n+1+a)(n+1+b)} .
\end{equation}
The four-term recurrence relation is determined by $ X_{n,n}, Y_{n,n} $ which have the rational forms
\begin{align}
	X_{n,n} & = \frac{(n+1)(n+1+a)(n+a+b+1)}{2n+a+b+2} - \frac{n(n+a)(n+a+b)}{2n+a+b} ,
\\
	Y_{n,n} & = \frac{(n+1)(n+1+b)(n+a+b+1)}{2n+a+b+2} - \frac{n(n+b)(n+a+b)}{2n+a+b} .
\end{align} 

The spectral polynomials are
\begin{equation}
	\Theta^{+}_{n} = -\frac{S_{n}}{S_{n+1}} ,
\qquad
	\Theta^{-}_{n} = +\frac{S_{n-1}}{S_{n}} ,
\qquad
	\Omega_{n} = x-n-a-b-1+Y_{n,n} ,
\end{equation}
and the Lax matrix has entries
\begin{align}
	A_{n,11} & = n+1 - \frac{(n+1)(n+1+a)(n+a+b+1)}{(2n+a+b+2)(2n+a+b+1)} ,
\\
	A_{n,12} & = \frac{S_{n+1}}{S_{n}}\frac{(n+1)(n+1+a)(n+a+b+1)}{(2n+a+b+2)(2n+a+b+1)}\left(x+Y_{n,n}\right) ,
\\
	A_{n,13} & = \frac{S_{n+1}S_{n-1}}{S_{n}^2}\frac{(n+1)(n+1+a)(n+a+b+1)}{(2n+a+b+2)(2n+a+b+1)} ,
\\
	A_{n,21} & = -\frac{S_{n}}{S_{n+1}} ,
\\
	A_{n,22} & =  x-n-a-b-1+Y_{n,n} ,
\\
	A_{n,23} & = \frac{S_{n-1}}{S_{n}} ,
\\
	A_{n,31} & = -\frac{S_{n-1}}{S_{n+1}} \frac{(2n+a+b)(2n+a+b-1)}{n(n+a)(n+a+b)} ,
\\
	A_{n,32} & = \frac{S_{n-1}}{S_{n}} \frac{(2n+a+b)(2n+a+b-1)}{n(n+a)(n+a+b)}\left(x-X_{n,n}\right)  ,
\\
	A_{n,33} & = \frac{n(n+b)(n+a+b)}{(2n+a+b+1)(2n+a+b)} - n-a-b .
\end{align}
The invariants are
$ {\rm Tr}A_{n} = x-2a-b $,
$ \Lambda_2(A_{n}) = (a+b) (a-x)+x $, and
$ {\rm Det}A_{n} = -(n+1)(a+b+n)x $.

This concludes our discussion of the generic fundamental structures of bi-orthogonal systems with the Cauchy kernel for broad classes of absolutely continuous weights which factorise
and the bi-moment matrices containing finite moments which are always non-singular.

\section{Bi-orthogonal Polynomials and Lax Pairs of Spectral and Deformation Derivatives in the Cauchy-Laguerre Setting}\label{Cauchy-Laguerre-BOPS}
\setcounter{equation}{0}

In the second phase of our study we will restrict ourselves to weights that satisfy the semi-classical criteria of satisfying a Pearson type equation,
and in order not to deviate from our original aims we will resume the Laguerre form as per our definitions and conventions \S \ref{BHFT_to_UCL}.

\subsection{Spectral and Deformation Derivatives of the Bi-orthogonal Polynomials and Associated Functions}

Our intention here is to derive all the required derivatives, both the spectral $x,y$ and deformation $s,t$, 
from first principles utilising the explicit semi-classical character of the univariate densities. 
\begin{proposition}
The Stieltjes functions $f_{1}(z,s), f_{2}(z,t)$ satisfy the following linear, inhomogeneous partial differential equations with respect to the spectral variable $z$ (either $x$ or $y$)
\begin{align}
	\partial_{z}f_{1}(z,s) & = \frac{a-z}{z}f_{1}(z,s) - \frac{a}{z}f_{1}(0,s) - \xi s^{a}e^{-s}\frac{1}{z-s} ,
\label{Dz_f1}
\\
	\partial_{z}f_{2}(z,t) & = \frac{b-z}{z}f_{2}(z,t) - \frac{b}{z}f_{2}(0,t) - \psi t^{b}e^{-t}\frac{1}{z-t} ,
\label{Dz_f2}
\end{align}
where we observe regular singularities at $z=0,s$ and $z=0,t$ respectively and both have an irregular singularity of Poincare rank unity at $ z=\infty$.
In addition they satisfy linear partial differential equations with respect to the deformation variables $s,t$
\begin{equation}
	\partial_{s}f_{1}(z,s) = \xi s^{a}e^{-s}\frac{1}{z-s} ,
\qquad
	\partial_{t}f_{2}(z,t) = \psi t^{b}e^{-t}\frac{1}{z-t} .
\label{DeformDiff_f1+f2}
\end{equation}
\end{proposition}
\begin{proof}
The spectral derivatives \eqref{Dz_f1},\eqref{Dz_f2} follow from differentiating the definitions \eqref{Stieltjes_f1},\eqref{Stieltjes_f2} under the integration, 
employing the identity $ \partial_{z}(z-x)^{-1} = -\partial_{x}(z-x)^{-1} $ in the case of $f_{1}$ say, and integrating by parts. 
At this juncture we then utilise the univariate forms of either the weight derivatives \eqref{wCL_spectral-diff:a} or \eqref{wCL_deform-diff:b}.
One then makes suitable partial fraction expansions of the resulting integrands into rational forms with respect to $x$ with simple poles 
and identifies these integrals with the definitions of the Stieltjes functions.
The deformation derivatives are found by a simple differentiation of the definitions \eqref{Stieltjes_f1},\eqref{Stieltjes_f2}.
\end{proof}

Before deriving the spectral derivatives of the bi-orthogonal polynomials we require the calculation of a preliminary quantity.
Identities \eqref{M_Id:a} and \eqref{M_Id:b} need to be re-written in the polynomial basis, which is
\begin{multline}
	\int_{0}^{\infty}\int_{0}^{\infty} \frac{dxdy}{(x+y)^2} w(x,y) xP_{n}(x)Q_{m}(y)
\\
	= \langle x\partial_x P_{n}, Q_{m} \rangle + (a+1)\delta_{n,m} - \langle xP_{n}, Q_{m} \rangle + \xi s^{a+1}e^{-s}P_{n}(s)Q^{(1)}_{m}(-s) ,
\label{SQkernel_x}
\end{multline}
and
\begin{multline}
	\int_{0}^{\infty}\int_{0}^{\infty} \frac{dxdy}{(x+y)^2} w(x,y) yP_{n}(x)Q_{m}(y)
\\
	= \langle P_{n}, y\partial_y Q_{m} \rangle + (b+1)\delta_{n,m} - \langle P_{n}, yQ_{m} \rangle + \psi t^{b+1}e^{-t}P^{(1)}_{n}(-t)Q_{m}(t) .
\label{SQkernel_y}
\end{multline}
In addition we have the following sum identity.
\begin{lemma}
The sum of inner products $ \langle x\partial_x P_{n}, Q_{m} \rangle $ and $ \langle P_{n}, y\partial_y Q_{m} \rangle $ has the evaluation
\begin{multline}
	\langle x\partial_x P_{n}, Q_{m} \rangle + \langle P_{n}, y\partial_y Q_{m} \rangle
	= \pi_{n}\eta_{m} - (a+b+1)\delta_{n,m} - \xi s^{a+1}e^{-s}P_{n}(s)Q^{(1)}_{m}(-s) - \psi t^{b+1}e^{-t}P^{(1)}_{n}(-t)Q_{m}(t) .
\label{sum_BOP-derivatives}
\end{multline}
\end{lemma}
Having the above results we can evaluate certain integrals with the squared Cauchy kernel.
\begin{lemma}\label{inverseSQ}
The integrals on the left-hand sides of \eqref{SQkernel_x} and \eqref{SQkernel_y} have the case-by-case evaluations:
\begin{align}
	\int_{0}^{\infty}\int_{0}^{\infty} \frac{dxdy}{(x+y)^2} w(x,y) xP_{n}(x)Q_{m}(y) 
	& = \left\{
		\begin{alignedat}{3}
			\xi s^{a+1}e^{-s}P_{n}(s)Q^{(1)}_{m}(-s) ,&&\quad m \geq n+2 ,\\
			-\frac{S_{n}}{S_{n+1}} + \xi s^{a+1}e^{-s}P_{n}(s)Q^{(1)}_{n+1}(-s) ,&&\quad m = n+1 ,\\
			n+a+1 - X_{n,n} + \xi s^{a+1}e^{-s}P_{n}(s)Q^{(1)}_{n}(-s) ,&&\quad m=n , \\
			\frac{S_{n-1}}{S_{n}} - \psi t^{b+1}e^{-t}P^{(1)}_{n}(-t)Q_{n-1}(t) ,&&\quad m = n-1 ,\\
			- \psi t^{b+1}e^{-t}P^{(1)}_{n}(-t)Q_{m}(t) ,&&\quad m \leq n-2 ,\\
		\end{alignedat}
		\right.
\\
	\int_{0}^{\infty}\int_{0}^{\infty} \frac{dxdy}{(x+y)^2} w(x,y) P_{n}(x)yQ_{m}(y) 
	& = \left\{
		\begin{alignedat}{3}
			-\xi s^{a+1}e^{-s}P_{n}(s)Q^{(1)}_{m}(-s) ,&&\quad m \geq n+2 ,\\
			 \frac{S_{n}}{S_{n+1}} - \xi s^{a+1}e^{-s}P_{n}(s)Q^{(1)}_{n+1}(-s) ,&&\quad m = n+1 ,\\
			n+b+1 - Y_{n,n} + \psi t^{b+1}e^{-t}P^{(1)}_{n}(-t)Q_{n}(t) ,&&\quad m=n , \\
			-\frac{S_{n-1}}{S_{n}} + \psi t^{b+1}e^{-t}P^{(1)}_{n}(-t)Q_{n-1}(t) ,&&\quad m = n-1 ,\\
			 \psi t^{b+1}e^{-t}P^{(1)}_{n}(-t)Q_{m}(t) ,&&\quad m \leq n-2 ,\\
		\end{alignedat}	
		\right.
\end{align} 
\end{lemma}

We now have made all the preparations in order to deduce the spectral derivatives of the polynomial pair $P_{n}$, $Q_{n}$.
\begin{proposition}
The polynomials $ P_{n}, Q_{n} $ satisfy the following linear, first order homogeneous partial differential equations with respect to the spectral variables $x$ or $y$
\begin{multline}
	x\partial_{x} P_{n}(x) =
	-(n+a+b+1)P_{n}(x) - \pi_{n}\hat{P}_{n}(x) 
\\
	- \xi s^{a+1}e^{-s}\frac{P_{n}(s)}{x-s}
	\left\{
	\left[ \check{Q}^{(1)}_{n}(-s)-1 \right]\hat{P}_{n}(x) 
	+ \frac{S_{n}}{S_{n+1}} Q^{(1)}_{n+1}(-s)P_{n}(x)
	+ \frac{S_{n}}{S_{n+1}} Q^{(1)}_{n}(-s)P_{n+1}(x)
	\right\}
\\
	- \psi t^{b+1}e^{-t}\frac{P^{(1)}_{n}(-t)}{x+t}
	\left\{
	\check{Q}_{n}(t)\hat{P}_{n}(x) + \frac{S_{n}}{S_{n+1}} Q_{n+1}(t)P_{n}(x) + \frac{S_{n}}{S_{n+1}} Q_{n}(t)P_{n+1}(x)
	\right\} ,
\\
	= -(n+a+b+1)P_{n}(x) - \pi_{n}\hat{P}_{n}(x) - \xi s^{a+1}e^{-s}P_{n}(s) K^{0,1}_{n}(x,-s) - \psi t^{b+1}e^{-t}P^{(1)}_{n}(-t) K^{0,0}_{n}(x,t) ,
\label{P_diff-x}
\end{multline}
and
\begin{multline}
	y\partial_{y} Q_{n}(y) =
	-(n+a+b+1)Q_{n}(y) - \eta_{n}\check{Q}_{n}(y) 
\\
	- \xi s^{a+1}e^{-s}\frac{Q^{(1)}_{n}(-s)}{s+y}
	\left\{
	\hat{P}_{n}(s)\check{Q}_{n}(y) + \frac{S_{n}}{S_{n+1}}P_{n+1}(s) Q_{n}(y) + \frac{S_{n}}{S_{n+1}} P_{n}(s)Q_{n+1}(y)
	\right\}
\\
	- \psi t^{b+1}e^{-t}\frac{Q_{n}(t)}{y-t}
	\left\{
	\left[ \hat{P}^{(1)}_{n}(-t)-1 \right]\check{Q}_{n}(y) 
	+ \frac{S_{n}}{S_{n+1}} P^{(1)}_{n+1}(-t)Q_{n}(y)
	+ \frac{S_{n}}{S_{n+1}} P^{(1)}_{n}(-t)Q_{n+1}(y) 
	\right\} ,
\\
	= -(n+a+b+1)Q_{n}(y) - \eta_{n}\check{Q}_{n}(y) - \xi s^{a+1}e^{-s}Q^{(1)}_{n}(-s) K^{0,0}_{n}(s,y) - \psi t^{b+1}e^{-t}Q_{n}(t) K^{1,0}_{n}(-t,y) .
\label{Q_diff-y}
\end{multline}
\end{proposition}
\begin{proof}
We consider only the $P_{n}$ polynomials as the method applies identically for the $Q_{n}$ case.
Let $ x\partial_{x} P_{n}(x) = \sum_{j=0}^{n} a_{n,j} P_{j}(x) $ and consequently $ a_{n,j} = \langle x\partial_{x} P_{n}, Q_{j} \rangle $.
Integrating by parts and assuming $ {\rm Re}(a) > -1 $ we find 
\begin{equation*}
	\langle x\partial_{x} P_{n}, Q_{j} \rangle = -\int_{} dxdy P_{n}Q_{j} \left[ \frac{y}{(x+y)^2}w + \frac{1}{x+y}x\partial_x w \right] .
\end{equation*}
Utilising the semi-classical Pearson equation \eqref{wCL_spectral-diff:a} we deduce
\begin{equation*}
	a_{n,j} = - \int_{} dxdy\; wP_{n}Q_{j} \left[ \frac{y}{(x+y)^2} + \frac{a-x}{x+y} \right] - \xi sw_{1}(s)P_{n}(s)Q^{(1)}_{j}(-s) .
\end{equation*}
For the $Q$ case we define $ y\partial_{y} Q_{n}(y) = \sum_{k=0}^{n} b_{n,k} Q_{k}(y) $ and in the same manner find
\begin{equation*}
	b_{n,k} = - \int_{} dxdy\; wP_{k}Q_{n} \left[ \frac{x}{(x+y)^2} + \frac{b-y}{x+y} \right] - \psi tw_{2}(t)P^{(1)}_{k}(-t)Q_{n}(t) .
\end{equation*}
The last step is to evaluate the inverse square kernels using the results of Lemma \ref{inverseSQ} and the second term of the integrands using Prop. \ref{XY-elements}.
\end{proof}

We continue the process for the pair of first associated functions, using the results of the previous proposition.
\begin{proposition}
The associated functions $ P^{(1)}_{n}, Q^{(1)}_{n} $ satisfy the following linear, first order homogeneous partial differential equations with respect to the spectral variables $x$ or $y$
\begin{multline}
	x\partial_{x} P^{(1)}_{n}(x) =
	-(n+b+1+x)P^{(1)}_{n}(x) + \pi_{n}\left( 1-\hat{P}^{(1)}_{n}(x) \right)
\\
	- \xi s^{a+1}e^{-s}\frac{P_{n}(s)}{x-s}
	\left\{
	\left[ \check{Q}^{(1)}_{n}(-s)-1 \right]\left[ \hat{P}^{(1)}_{n}(x)-1 \right] 
	+ \frac{S_{n}}{S_{n+1}} Q^{(1)}_{n+1}(-s)P^{(1)}_{n}(x)
	+ \frac{S_{n}}{S_{n+1}} Q^{(1)}_{n}(-s)P^{(1)}_{n+1}(x)
	\right\}
\\
	- \psi t^{b+1}e^{-t}\frac{P^{(1)}_{n}(-t)}{x+t}
	\left\{
	\check{Q}_{n}(t)\left[ \hat{P}^{(1)}_{n}(x)-1 \right] + \frac{S_{n}}{S_{n+1}} Q_{n+1}(t)P^{(1)}_{n}(x) + \frac{S_{n}}{S_{n+1}} Q_{n}(t)P^{(1)}_{n+1}(x)
	\right\} ,
\\
	= -(n+b+1+x)P^{(1)}_{n}(x) + \pi_{n}\left( 1-\hat{P}^{(1)}_{n}(x) \right) - \xi s^{a+1}e^{-s}P_{n}(s) K^{1,1}_{n}(x,-s) - \psi t^{b+1}e^{-t}P^{(1)}_{n}(-t) K^{0,1}_{n}(x,t) ,
\label{P1_diff-x}
\end{multline}
and
\begin{multline}
	y\partial_{y} Q^{(1)}_{n}(y) =
	-(n+a+1+y)Q^{(1)}_{n}(y) + \eta_{n}\left[ 1-\check{Q}^{(1)}_{n}(y) \right] 
\\
	- \xi s^{a+1}e^{-s}\frac{Q^{(1)}_{n}(-s)}{s+y}
	\left\{
	\hat{P}_{n}(s)\left[ \check{Q}^{(1)}_{n}(y)-1 \right] + \frac{S_{n}}{S_{n+1}}P_{n+1}(s)Q^{(1)}_{n}(y) + \frac{S_{n}}{S_{n+1}}P_{n}(s)Q^{(1)}_{n+1}(y)
	\right\}
\\
	- \psi t^{b+1}e^{-t}\frac{Q_{n}(t)}{y-t}
	\left\{
	\left[ \hat{P}^{(1)}_{n}(-t)-1 \right]\left[ \check{Q}^{(1)}_{n}(y)-1 \right] 
	+ \frac{S_{n}}{S_{n+1}} P^{(1)}_{n+1}(-t)Q^{(1)}_{n}(y)
	+ \frac{S_{n}}{S_{n+1}} P^{(1)}_{n}(-t)Q^{(1)}_{n+1}(y) 
	\right\} , 
\\
	= -(n+a+1+y)Q^{(1)}_{n}(y) + \eta_{n}\left[ 1-\check{Q}^{(1)}_{n}(y) \right] - \xi s^{a+1}e^{-s}Q^{(1)}_{n}(-s) K^{0,1}_{n}(s,y) - \psi t^{b+1}e^{-t}Q_{n}(t) K^{1,1}_{n}(-t,y) .
\label{Q1_diff-y}
\end{multline}
\end{proposition}
\begin{proof}
Again we only provide details for the $ P^{(1)}_{n}(z) $ function as the methods are the same for the $ Q^{(1)}_{n}(z) $ function.
One differentiates the definition \eqref{AssocFn_P1} with respect to $z$ and interchanges the differentiation and integration over $x$.
The identity $ \partial_{z}(z-x)^{-1} = -\partial_{x}(z-x)^{-1} $ is employed again and we integrate by parts with respect to $x$. 
For the term with $ \partial_{x}w_1(x) $ we use the explicit weight formula and for this term and the others, do a partial fraction decomposition of the rational functions of $x$.
The resulting integrals have a term with $ \partial_{x}P_{n}(x) $ whilst the others do not and these latter integrals can be evaluated via the previous definitions.
Now we insert the results of the previous proposition in the form of \eqref{P_diff-x} which has a rational dependence on $x$ with simple poles, 
namely the first of the two formulae.
Now all the remaining integrals can be evaluated in terms of associated functions with arguments of $z$, $s$ or $-t$ as the case maybe.
Our final step is to simplify the resulting expression by using the identities \eqref{anti-CDsum-01} with $x=s$ and \eqref{anti-CDsum-10} with $x=-t$.
Restoring $ z\mapsto x $ yields \eqref{P1_diff-x}. 
\end{proof}

We now turn our attention to the deformation derivatives of the polynomial pair.
\begin{proposition}
The polynomials $ P_{n}, Q_{n} $ satisfy the following linear, first order homogeneous partial differential equations with respect to the deformation variables $s,t$
\begin{align}
	\partial_{s} P_{n}(x) & =
	\xi s^{a}e^{-s} P_{n}(s) \left[ \tfrac{1}{2}P_{n}(x)Q^{(1)}_{n}(-s) + K^{0,1}_{n-1}(x,-s) \right] ,
\label{Pdiff-s}\\
	\partial_{s} Q_{n}(y) & =
	\xi s^{a}e^{-s} Q^{(1)}_{n}(-s) \left[ \tfrac{1}{2}P_{n}(s)Q_{n}(y) + K^{0,0}_{n-1}(s,y) \right] ,
\label{Qdiff-s}\\
	\partial_{t} P_{n}(x) & =
	\psi t^{b}e^{-t} P^{(1)}_{n}(-t) \left[ \tfrac{1}{2}P_{n}(x)Q_{n}(t) + K^{0,0}_{n-1}(x,t) \right] ,
\label{Pdiff-t}\\
	\partial_{t} Q_{n}(y) & =
	\psi t^{b}e^{-t} Q_{n}(t) \left[ \tfrac{1}{2}P^{(1)}_{n}(-t)Q_{n}(y) + K^{1,0}_{n-1}(-t,y) \right] .
\label{Qdiff-t}
\end{align}
\end{proposition}
\begin{proof}
We will demonstrate the proof only for case of $ \partial_{s} P_{n} $.
Writing the expansion formula $ \partial_{s}P_{n}(x) = \sum_{j=0}^{n}c_{n,j}P_{j}(x) $ we note that $ c_{n,j} = \langle \partial_{s}P_{n},Q_{j} \rangle $, $ 0\leq j\leq n $.
For $ 0\leq i \leq n $ the orthonormality condition is expressible as
\begin{equation}
	\delta_{i,0} = \langle P_{n},Q_{n-i} \rangle = \int\, dxdy\frac{w(x,y)}{x+y}P_{n}(x)Q_{n-i}(y) .
\end{equation}
Differentiating this with respect to $s$, employing \eqref{wCL_deform-diff:a} and the definition \eqref{AssocFn_Q1} we deduce
\begin{equation}
	0 = -\xi s^{a}e^{-s}P_{n}(s)Q^{(1)}_{n-i}(-s) + c_{n,n-i} + \delta_{i,0}\langle P_{n},\partial_{s}Q_{n} \rangle .
\end{equation}
Thus for $ 0< i\leq n-1 $ this furnishes us with the solution $ c_{n,n-i} = \xi s^{a}e^{-s}P_{n}(s)Q^{(1)}_{n-i}(-s) $ however the case $ i=0 $ will have to be treated separately.
For this case we have
\begin{equation}
	 \langle P_{n},\partial_{s}Q_{n} \rangle = - c_{n,n} + \xi s^{a}e^{-s}P_{n}(s)Q^{(1)}_{n}(-s) .
\end{equation}
Also note that $ c_{n,n} = S_{n}^{-1}\partial_{s}S_{n} $, as can be found by comparing the leading terms of the expansion formula.
However in the $ Q_{n} $ expansion formula $ \partial_{s}Q_{n}(x) = \sum_{j=0}^{n}d_{n,j}Q_{j}(x) $ we observe that $ d_{n,n} = S_{n}^{-1}\partial_{s}S_{n} $ by the same argument.
So $ \langle P_{n},\partial_{s}Q_{n} \rangle =: d_{n,n} = c_{n,n} $ and thus
\begin{equation}
	c_{n,n} = \tfrac{1}{2}\xi s^{a}e^{-s}P_{n}(s)Q^{(1)}_{n}(-s) .
\end{equation}
In conclusion we have \eqref{Pdiff-s}.
\end{proof}

Our final task in this subsection is now to compute the deformation derivatives of the first associated functions.
\begin{proposition}
The associated functions $ P^{(1)}_{n}, Q^{(1)}_{n} $ satisfy the following linear, first order homogeneous partial differential equations with respect to the deformation variables $s,t$
\begin{align}
	\partial_{s} P^{(1)}_{n}(x) & =
	\xi s^{a}e^{-s} P_{n}(s) \left[ \tfrac{1}{2}P^{(1)}_{n}(x)Q^{(1)}_{n}(-s) + K^{1,1}_{n-1}(x,-s) \right] ,
\label{P1diff-s}\\
	\partial_{s} Q^{(1)}_{n}(y) & =
	\xi s^{a}e^{-s} Q^{(1)}_{n}(-s) \left[ \tfrac{1}{2}P_{n}(s)Q^{(1)}_{n}(y) + K^{0,1}_{n-1}(s,y) \right] ,
\label{Q1diff-s}\\
	\partial_{t} P^{(1)}_{n}(x) & =
	\psi t^{b}e^{-t} P^{(1)}_{n}(-t) \left[ \tfrac{1}{2}P^{(1)}_{n}(x)Q_{n}(t) + K^{1,0}_{n-1}(x,t) \right] ,
\label{P1diff-t}\\
	\partial_{t} Q^{(1)}_{n}(y) & =
	\psi t^{b}e^{-t} Q_{n}(t) \left[ \tfrac{1}{2}P^{(1)}_{n}(-t)Q^{(1)}_{n}(y) + K^{1,1}_{n-1}(-t,y) \right] .
\label{Q1diff-t}
\end{align}
\end{proposition}
\begin{proof}
We only consider $ P^{(1)}_{n}(x) $ and derivation with respect to $s$ as the other cases are amenable to an identical treatment.
Firstly we differentiate the definition for $ P^{(1)}_{n}(x) $ \eqref{AssocFn_P1} with respect to $s$ and find two terms arise: 
one containing $ \partial_{s}w_1(x,s) $ and one containing $ \partial_{s}P_{n}(x) $.
For the former we use $ \xi\delta_{x-s}x^{a}e^{-x} $ and for the latter we use formula \eqref{Pdiff-s} from the previous proposition.
The integrals can be immediately evaluated using the definition \eqref{AssocFn_P1}, and in addition we recognise the definition of the
$ K^{1,1}_{n-1}(x,-s) $ kernel.
\end{proof}

The reader should note that the Stieltjes functions, bi-orthogonal polynomials and associated functions are fully described through the complete notations
$ f_1(z;s) $, $ P_{n}(x;s,t) $, $ P^{(1)}_{n}(x;s,t) $, etc even though some of these dependencies, usually the deformation variables,
will be suppressed from time to time. 
We will find in the ensuing developments that the spectral variable of these objects can be tied to a deformation variable, either $ \pm s $ or $ \pm t $,
and so the total dependence on the deformation variable comes from two sources. 
We will denote the pure deformation derivative by $ \partial_s, \partial_t $ whereas a total derivative is by $ d/ds, d/dt $. 

\begin{definition}
Let the pair of $\sigma$-functions $ \sigma_{n}(s,t;\xi,\psi), \tau_{n}(s,t;\xi,\psi) $ be defined by
\begin{equation}
	\sigma_{n}(s,t) := s\partial_{s} \log Z^\text{\rm C-L2M}_{n}, \qquad
	\tau_{n}(s,t) := t\partial_{t} \log Z^\text{\rm C-L2M}_{n} .
\label{sigma_Defn}
\end{equation}
\end{definition}

\begin{proposition}\label{sigma-form}
The $\sigma$-functions have evaluations in terms of reproducing kernels
\begin{align}
	\sigma_{n}(s,t) & = -\xi s^{a+1}e^{-s} K_{n-1}^{0,1}(s,-s;s,t) ,
\label{sigma}
\\
	\tau_{n}(s,t) & = - \psi t^{b+1}e^{-t}  K_{n-1}^{1,0}(-t,t;s,t) .
\label{tau}
\end{align}
They satisfy the compatibility relation $ t\partial_{t} \sigma_{n}(s,t) = s\partial_{s} \tau_{n}(s,t) $.
\end{proposition}
\begin{proof}
The inverse norm of the bi-orthogonal system satisfies the first order partial differential equations with respect to the deformation variables $s,t$
\begin{align}
	\frac{\partial_{s} S^2_{n}}{S^2_{n}} & = \xi s^{a}e^{-s} P_{n}(s)Q^{(1)}_{n}(-s) ,
\label{PQNormdiff-s}\\
	\frac{\partial_{t} S^2_{n}}{S^2_{n}} & = \psi t^{b}e^{-t} P^{(1)}_{n}(-t)Q_{n}(t) .
\label{PQNormdiff-t}
\end{align}
From \eqref{norm+leading-coeff} we have $ \log Z_{n}^\text{C-L2M} = -\sum_{l=0}^{n-1}\log S^2_{l} $ 
and the resulting sum of derivatives using the two relations above is recognised as the definitions of the $0,1$ and $1,0$ kernels.
\end{proof}

\subsection{Lax Matrices}
In this subsection we will consolidate the results of the previous subsection into a standard matrix formulation which furnishes greater economy and transparency.
In addition we will exclusively report only those results applying to the $P$-system for the sake of brevity and in the knowledge that through the symmetry exchanges 
$ x \leftrightarrow y,\,s \leftrightarrow t,\,a \leftrightarrow b,\,\xi \leftrightarrow \psi,\,P_{n} \leftrightarrow Q_{n} $ 
one can write down the corresponding results for the $Q$-system by inspection.
The polynomial form of the spectral derivative of the bi-orthogonal $P$-polynomial is
\begin{equation}
	x(x-s)(x+t) \partial_{x} P_{n}(x) = \Theta^{+}_{n}(x) P_{n+1}(x) + \Omega_{n}(x) P_{n}(x) + \Theta^{-}_{n}(x) P_{n-1}(x) ,
\label{Diff-x_P}
\end{equation}
while the polynomial forms of the deformation derivatives of this polynomial are given by
\begin{equation}
	(x-s) \partial_{s} P_{n}(x) = \Xi^{+}_{n}(x) P_{n+1}(x) + \Upsilon_{n}(x) P_{n}(x) + \Xi^{-}_{n}(x) P_{n-1}(x) ,
\label{Diff-s_P}
\end{equation}
and
\begin{equation}
	(x+t) \partial_{t} P_{n}(x) = \Phi^{+}_{n}(x) P_{n+1}(x) + \Psi_{n}(x) P_{n}(x) + \Phi^{-}_{n}(x) P_{n-1}(x) .
\label{Diff-t_P}
\end{equation}
Our next result gives explicit forms for the coefficient polynomials appearing in the forgoing spectral derivatives.
\begin{proposition}
The spectral polynomials $\Theta^{\pm}_{n}(x)$, ${\rm deg}_{x}\Theta^{\pm}_{n} = 2$, have the evaluations
\begin{multline}
	\frac{S_{n+1}}{S_{n}}\Theta^{+}_{n}(x) = -(x-s)(x+t) 
\\
	- \xi s^{a+1}e^{-s}\frac{(x+t)}{\pi_{n}\eta_{n}}P_{n}(s) 
	\left[ \frac{S_{n}}{S_{n+1}}Q^{(1)}_{n+1}(-s)+(Y_{n,n}+s)Q^{(1)}_{n}(-s)-\frac{S_{n-1}}{S_{n}}Q^{(1)}_{n-1}(-s) \right]
\\
	- \psi t^{b+1}e^{-t}\frac{(x-s)}{\pi_{n}\eta_{n}}P^{(1)}_{n}(-t) 
	\left[ \frac{S_{n}}{S_{n+1}}Q_{n+1}(t)+(Y_{n,n}-t)Q_{n}(t)-\frac{S_{n-1}}{S_{n}}Q_{n-1}(t) \right] ,
\label{Theta+}
\end{multline}
and
\begin{multline}
	\frac{S_{n}}{S_{n-1}}\Theta^{-}_{n}(x) = (x-s)(x+t) 
\\
	+ \xi s^{a+1}e^{-s}\frac{(x+t)}{\pi_{n}\eta_{n}}P_{n}(s) 
	\left[ \frac{S_{n}}{S_{n+1}}Q^{(1)}_{n+1}(-s)-(X_{n,n}-s)Q^{(1)}_{n}(-s)-\frac{S_{n-1}}{S_{n}}Q^{(1)}_{n-1}(-s) \right]
\\
	+ \psi t^{b+1}e^{-t}\frac{(x-s)}{\pi_{n}\eta_{n}}P^{(1)}_{n}(-t) 
	\left[ \frac{S_{n}}{S_{n+1}}Q_{n+1}(t)-(X_{n,n}+t)Q_{n}(t)-\frac{S_{n-1}}{S_{n}}Q_{n-1}(t) \right] ,
\label{Theta-}
\end{multline}
while the spectral polynomial $\Omega_{n}(x)$, ${\rm deg}_{x}\Omega_{n} = 3$, has the evaluation
\begin{multline}
	\Omega_{n}(x) = (x-s)(x+t)\left[ x+Y_{n,n}-n-a-b-1 \right] 
\\
	- \xi s^{a+1}e^{-s}\frac{(x+t)}{\pi_{n}\eta_{n}}P_{n}(s) 
	\left[ \frac{S_{n}}{S_{n+1}}(X_{n,n}-x)Q^{(1)}_{n+1}(-s)+(Y_{n,n}+x)(X_{n,n}-s)Q^{(1)}_{n}(-s)+\frac{S_{n-1}}{S_{n}}(Y_{n,n}+x)Q^{(1)}_{n-1}(-s) \right]
\\
	- \psi t^{b+1}e^{-t}\frac{(x-s)}{\pi_{n}\eta_{n}}P^{(1)}_{n}(-t) 
	\left[ \frac{S_{n}}{S_{n+1}}(X_{n,n}-x)Q_{n+1}(t)+(Y_{n,n}+x)(X_{n,n}+t)Q_{n}(t)+\frac{S_{n-1}}{S_{n}}(Y_{n,n}+x)Q_{n-1}(t) \right] .
\label{Omega}
\end{multline}
\end{proposition}

The spectral Lax matrices $ \left( A^{(\nu)}_{l,m}(s,t) \right)_{l,m\geq 0} $ corresponding to the singular points $ \nu \in \{0,s,-t,\infty\} $
appear in the rational form of the spectral derivative of the bi-orthogonal polynomial
\begin{multline}
	\partial_{x} P_{n}(x) =
	\left( \frac{A^{(0)}_{n,n+1}}{x} + \frac{A^{(s)}_{n,n+1}}{x-s} + \frac{A^{(-t)}_{n,n+1}}{x+t} + A^{(\infty)}_{n,n+1} \right) P_{n+1}(x)
\\
	+ \left( \frac{A^{(0)}_{n,n}}{x} + \frac{A^{(s)}_{n,n}}{x-s} + \frac{A^{(-t)}_{n,n}}{x+t} + A^{(\infty)}_{n,n} \right) P_{n}(x)
\\
	+ \left( \frac{A^{(0)}_{n,n-1}}{x} + \frac{A^{(s)}_{n,n-1}}{x-s} + \frac{A^{(-t)}_{n,n-1}}{x+t} + A^{(\infty)}_{n,n-1} \right) P_{n-1}(x) .
\label{LaxA-P}
\end{multline}
The deformation Lax matrices $ \left( B^{(\nu)}_{l,m}(s,t) \right)_{l,m\geq 0} $, $ \left( C^{(\nu)}_{l,m}(s,t) \right)_{l,m\geq 0} $ 
corresponding to the singular points $ \nu \in \{s,\infty\} $ and $ \nu \in \{-t,\infty\} $ respectively
appear in the rational form of the deformation derivatives of the bi-orthogonal polynomial
\begin{equation}
	\partial_{s} P_{n}(x)
	= \left( \frac{B^{(s)}_{n,n+1}}{x-s} + B^{(\infty)}_{n,n+1} \right) P_{n+1}(x)
	+ \left( \frac{B^{(s)}_{n,n}}{x-s} + B^{(\infty)}_{n,n} \right) P_{n}(x)
	+ \left( \frac{B^{(s)}_{n,n-1}}{x-s} + B^{(\infty)}_{n,n-1} \right) P_{n-1}(x) ,
\label{LaxB-P}
\end{equation}
and
\begin{equation}
	\partial_{t} P_{n}(x)
	= \left( \frac{C^{(-t)}_{n,n+1}}{x+t} + C^{(\infty)}_{n,n+1} \right) P_{n+1}(x)
	+ \left( \frac{C^{(-t)}_{n,n}}{x+t} + C^{(\infty)}_{n,n} \right) P_{n}(x)
	+ \left( \frac{C^{(-t)}_{n,n-1}}{x+t} + C^{(\infty)}_{n,n-1} \right) P_{n-1}(x) .
\label{LaxC-P}
\end{equation}
Assembling these results into modified - redefined - $3\times 3$ matrix variables, which are related to \eqref{3x3_PQ} by a gauge transformation,
\begin{equation}\label{new_3x3_PQ}
	\mathcal{P}_{n} := 	\left(
	\begin{array}{ccc}
	P_{n+1}(x) 	& e^{x}x^{-a}P^{(1)}_{n+1}(x) 	& x^{-a-b}P^{(2)}_{n+1}(x) \cr
	P_{n}(x) 	& e^{x}x^{-a}P^{(1)}_{n}(x) 	& x^{-a-b}P^{(2)}_{n}(x) \cr
	P_{n-1}(x) 	& e^{x}x^{-a}P^{(1)}_{n-1}(x) 	& x^{-a-b}P^{(2)}_{n-1}(x) \cr
	\end{array}
	\right) , \quad
	\mathcal{Q}_{n} := 	\left(
	\begin{array}{ccc}
	Q_{n+1}(y) 	& e^{y}y^{-b}Q^{(1)}_{n+1}(y) 	& y^{-a-b}Q^{(2)}_{n+1}(y) \cr
	Q_{n}(y) 	& e^{y}y^{-b}Q^{(1)}_{n}(y) 	& y^{-a-b}Q^{(2)}_{n}(y) \cr
	Q_{n-1}(y) 	& e^{y}y^{-b}Q^{(1)}_{n-1}(y) 	& y^{-a-b}Q^{(2)}_{n-1}(y) \cr
	\end{array}
\right) , 
\end{equation}
and the spectral and deformation Lax matrices 
$	\mathcal{A}_{n} = \left( A_{n+l,n+m} \right)_{l,m=+1,0,-1} $, 
$	\mathcal{B}_{n} = \left( B_{n+l,n+m} \right)_{l,m=+1,0,-1} $,
$	\mathcal{C}_{n} = \left( C_{n+l,n+m} \right)_{l,m=+1,0,-1} $
we can express the differential and recurrence relations as
\begin{equation}\label{xst-LaxPairs}
	\partial_{x} \mathcal{P}_{n} = \mathcal{A}_{n} \mathcal{P}_{n} ,
\qquad
	\partial_{s} \mathcal{P}_{n} = \mathcal{B}_{n} \mathcal{P}_{n} ,
\qquad
	\partial_{t} \mathcal{P}_{n} = \mathcal{C}_{n} \mathcal{P}_{n} ,
\qquad
	\mathcal{P}_{n+1} = \mathcal{K}_{n} \mathcal{P}_{n} .
\end{equation}
The elements of the coefficient elements are rational in $x$ with the following structure
\begin{align}
	A_{n+l,n+m} & = \sum_{\nu\in \{0,s,-t\}} \frac{1}{x-\nu} A^{(\nu)}_{n+l,n+m} + A^{(\infty)}_{n+l,n+m} ,
\\
	B_{n+l,n+m} & = \frac{1}{x-s} B^{(s)}_{n+l,n+m} + B^{(\infty)}_{n+l,n+m} ,
\\
	C_{n+l,n+m} & = \frac{1}{x+t} C^{(-t)}_{n+l,n+m} + C^{(\infty)}_{n+l,n+m} ,
\end{align}
or in $3\times 3$ matrix form defining the residue matrices $ \mathcal{A}^{(\cdot)}_{n} $, $ \mathcal{B}^{(\cdot)}_{n} $, $ \mathcal{C}^{(\cdot)}_{n} $
\begin{gather}
	\mathcal{A}_{n}(x) = \frac{1}{x} \mathcal{A}^{(0)}_{n} + \frac{1}{x-s} \mathcal{A}^{(s)}_{n} + \frac{1}{x+t} \mathcal{A}^{(-t)}_{n} + \mathcal{A}^{(\infty)}_{n} ,
\label{Rational_A}\\
	\mathcal{B}_{n}(x) = \frac{1}{x-s} \mathcal{B}^{(s)}_{n} + \mathcal{B}^{(\infty)}_{n} ,
\\
	\mathcal{C}_{n}(x) = \frac{1}{x+t} \mathcal{C}^{(-t)}_{n} + \mathcal{C}^{(\infty)}_{n} .
\end{gather}
It is evident that the spectral differential system \eqref{xst-LaxPairs},\eqref{Rational_A} has regular singularities at $ x=0,s,-t $ 
and an irregular one of Poincar{\'e} rank unity at $x=\infty$.

Summarising the results of our previous calculations, we have the explicit forms of the spectral residue matrices $ \mathcal{A}^{(\cdot)}_{n} $ featuring in \eqref{Rational_A}
\begin{equation}
	\mathcal{A}^{(\infty)}_{n} = 
	\begin{pmatrix}
	0	& \frac{\pi_{n+1}}{\pi_{n}}	& 0	\\
	0	& 1							& 0	\\
	0	& \frac{\pi_{n-1}}{\pi_{n}}	& 0
	\end{pmatrix} ,	
\label{Aresidue@Infty}	
\end{equation}
\begin{equation}
	\mathcal{A}^{(s)}_{n} = 
		- \xi s^{a}e^{-s}\frac{1}{\pi_{n}\eta_{n}} \mathrm{P}_{n}(s) \otimes\left( G_{n}(s,-s)\mathrm{Q}^{(1)}_{n}(-s) \right)^{T} ,
\label{AresidueSing@s}
\end{equation}
\begin{equation}
	\mathcal{A}^{(-t)}_{n} =
		+ \psi t^{b}e^{-t}\frac{1}{\pi_{n}\eta_{n}} \mathrm{P}^{(1)}_{n}(-t) \otimes\left( G_{n}(-t,t)\mathrm{Q}_{n}(t) \right)^{T} ,
\label{AresidueSing@-t}
\end{equation}
and
\begin{equation}
	\mathcal{A}^{(\Sigma)}_{n} := \mathcal{A}^{(0)}_{n} + \mathcal{A}^{(s)}_{n} + \mathcal{A}^{(-t)}_{n} =
\\
	\begin{pmatrix}
		n+1-\frac{S_{n}\pi_{n+1}}{S_{n+1}\pi_{n}}	& \mathcal{A}^{(0)}_{1,0}															& \frac{S_{n-1}\pi_{n+1}}{S_{n}\pi_{n}}	\\
		-\frac{S_{n}}{S_{n+1}}						& -a-1+\frac{S_{n}\pi_{n+1}}{S_{n+1}\pi_{n}}-\frac{S_{n-1}\pi_{n-1}}{S_{n}\pi_{n}}	& \frac{S_{n-1}}{S_{n}}	\\
		-\frac{S_{n}\pi_{n-1}}{S_{n+1}\pi_{n}}		& \mathcal{A}^{(0)}_{-1,0}															& -n-a-b+\frac{S_{n-1}\pi_{n-1}}{S_{n}\pi_{n}}
	\end{pmatrix} ,
\label{AresidueSing@0}
\end{equation}
where
\begin{multline}
	\mathcal{A}^{(0)}_{1,0} = \frac{\pi_{n+1}}{\pi_{n}}Y_{n,n}
\\
	+ \xi s^{a+1}e^{-s}\frac{1}{\pi_{n}\eta_{n}}P_{n+1}(s) 
	\left[ \frac{S_{n}}{S_{n+1}}Q^{(1)}_{n+1}(-s)-(X_{n,n}-s)Q^{(1)}_{n}(-s)-\frac{S_{n-1}}{S_{n}}Q^{(1)}_{n-1}(-s) \right]
\\
	+ \psi t^{b+1}e^{-t}\frac{1}{\pi_{n}\eta_{n}}P^{(1)}_{n+1}(-t) 
	\left[ \frac{S_{n}}{S_{n+1}}Q_{n+1}(t)-(X_{n,n}+t)Q_{n}(t)-\frac{S_{n-1}}{S_{n}}Q_{n-1}(t) \right] ,	
\label{AresidueSing@0:a}
\end{multline}
and
\begin{multline}
	\mathcal{A}^{(0)}_{-1,0} = -\frac{\pi_{n-1}}{\pi_{n}}X_{n,n}
\\
	+ \xi s^{a+1}e^{-s}\frac{1}{\pi_{n}\eta_{n}}P_{n-1}(s) 
	\left[ \frac{S_{n}}{S_{n+1}}Q^{(1)}_{n+1}(-s)+(Y_{n,n}+s)Q^{(1)}_{n}(-s)-\frac{S_{n-1}}{S_{n}}Q^{(1)}_{n-1}(-s) \right]
\\
	+ \psi t^{b+1}e^{-t}\frac{1}{\pi_{n}\eta_{n}}P^{(1)}_{n-1}(-t) 
	\left[ \frac{S_{n}}{S_{n+1}}Q_{n+1}(t)+(Y_{n,n}-t)Q_{n}(t)-\frac{S_{n-1}}{S_{n}}Q_{n-1}(t) \right] .	
\label{AresidueSing@0:b}
\end{multline}

Furthermore define certain auxiliary matrices
\begin{equation}
	\mathcal{B}^{(\infty,0)}_{n} := \tfrac{1}{2}\xi s^{a}e^{-s}
	\begin{pmatrix}
		P_{n+1}(s)Q^{(1)}_{n+1}(-s)	& 0								& 0								\\
		0							& -P_{n}(s)Q^{(1)}_{n}(-s)		& 0								\\
		0							& -2P_{n-1}(s)Q^{(1)}_{n}(-s)	& -P_{n-1}(s)Q^{(1)}_{n-1}(-s)	\\
	\end{pmatrix} ,
\end{equation}
\begin{equation}
	\mathcal{B}^{(\infty,\bar{0})}_{n} := \tfrac{1}{2}\xi s^{a}e^{-s}
	\begin{pmatrix}
		P_{n+1}(s)Q^{(1)}_{n+1}(-s)	& 0								& 0								\\
		0							& -P_{n}(s)Q^{(1)}_{n}(-s)		& 0								\\
		0							& -2P_{n}(s)Q^{(1)}_{n-1}(-s)	& -P_{n-1}(s)Q^{(1)}_{n-1}(-s)	\\
	\end{pmatrix} ,
\end{equation}
\begin{equation}
	\mathcal{C}^{(\infty,0)}_{n} := \tfrac{1}{2}\psi t^{b}e^{-t} 
				\begin{pmatrix}
					P^{(1)}_{n+1}(-t)Q_{n+1}(t)	& 0								& 0								\\
					0							& -P^{(1)}_{n}(-t)Q_{n}(t)		& 0								\\
					0							& -2P^{(1)}_{n-1}(-t)Q_{n}(t)	& -P^{(1)}_{n-1}(-t)Q_{n-1}(t)	\\
				\end{pmatrix} ,
\end{equation}
\begin{equation}
	\mathcal{C}^{(\infty,\bar{0})}_{n} := \tfrac{1}{2}\psi t^{b}e^{-t} 
				\begin{pmatrix}
					P^{(1)}_{n+1}(-t)Q_{n+1}(t)	& 0								& 0								\\
					0							& -P^{(1)}_{n}(-t)Q_{n}(t)		& 0								\\
					0							& -2P^{(1)}_{n}(-t)Q_{n-1}(t)	& -P^{(1)}_{n-1}(-t)Q_{n-1}(t)	\\
				\end{pmatrix} .
\end{equation}

The explicit forms of the $s$-deformation residue matrices $ \mathcal{B}^{(\cdot)}_{n} $ are:
$ \mathcal{B}^{(s)}_{n} = -\mathcal{A}^{(s)}_{n} $,
and
\begin{equation}
	\mathcal{B}^{(\infty)}_{n} = \mathcal{B}^{(\infty,0)}_{n}
	- \frac{1}{\pi_{n}\eta_{n}}\xi s^{a}e^{-s}\left[ \frac{S_{n}}{S_{n+1}}Q^{(1)}_{n+1}(-s)-(X_{n,n}-s)Q^{(1)}_{n}(-s)-\frac{S_{n-1}}{S_{n}}Q^{(1)}_{n-1}(-s) \right]
	\begin{pmatrix}
		0	& P_{n+1}(s)	& 0	\\
		0	& P_{n}(s)		& 0	\\
		0	& P_{n-1}(s)	& 0
	\end{pmatrix} .
\label{BresidueSing@Infty}		
\end{equation}
Similarly we find the explicit forms of the $t$-deformation residue matrices $ \mathcal{C}^{(\cdot)}_{n} $ are:
$ \mathcal{C}^{(-t)}_{n} = \mathcal{A}^{(-t)}_{n} $,
and
\begin{equation}
	\mathcal{C}^{(\infty)}_{n} = \mathcal{C}^{(\infty,0)}_{n}
	- \frac{1}{\pi_{n}\eta_{n}}\psi t^{b}e^{-t}\left[ \frac{S_{n}}{S_{n+1}}Q_{n+1}(t)-(X_{n,n}+t)Q_{n}(t)-\frac{S_{n-1}}{S_{n}}Q_{n-1}(t) \right]
	\begin{pmatrix}
		0	& P^{(1)}_{n+1}(-t)	& 0	\\
		0	& P^{(1)}_{n}(-t)	& 0	\\
		0	& P^{(1)}_{n-1}(-t)	& 0
	\end{pmatrix} .
\label{CresidueSing@Infty}		
\end{equation}

We record here the invariants of the $\mathcal{A}^{(\cdot)}_{n}$ residue matrices given above, as these exhibit the spectral data for the integrable system.
\begin{proposition}
Let us recall the definition $ \mathcal{A}^{\Sigma}_{n} := \mathcal{A}^{(0)}_{n}+\mathcal{A}^{(s)}_{n}+\mathcal{A}^{(-t)}_{n} $.
The traces of single residue matrices have the evaluations
\begin{equation*}
	\Tr(\mathcal{A}^{(0)}_{n}) = -2a-b,\quad 
	\Tr(\mathcal{A}^{(s)}_{n}) = 0 ,\quad 
	\Tr(\mathcal{A}^{(-t)}_{n}) = 0 ,\quad 
	\Tr(\mathcal{A}^{(\infty)}_{n}) = 1 ,\quad
	\Tr(\mathcal{A}_{n}(x)) = 1-\frac{2a+b}{x} ,
\end{equation*}
while their second invariants have the evaluations
\begin{equation*}
	\Lambda_2(\mathcal{A}^{(0)}_{n}) = a(a+b),\quad 
	\Lambda_2(\mathcal{A}^{(s)}_{n}) = 0 ,\quad 
	\Lambda_2(\mathcal{A}^{(-t)}_{n}) = 0 ,\quad 
	\Lambda_2(\mathcal{A}^{(\infty)}_{n}) = 0 ,
\end{equation*}
and their determinants have the evaluations
\begin{equation*}
	\det(\mathcal{A}^{(0)}_{n}) = 0,\quad 
	\det(\mathcal{A}^{(s)}_{n}) = 0 ,\quad 
	\det(\mathcal{A}^{(-t)}_{n}) = 0 ,\quad 
	\det(\mathcal{A}^{(\infty)}_{n}) = 0 .
\end{equation*}
In addition we record the second invariant
\begin{multline*}
	\Lambda_2(\mathcal{A}^{(\Sigma)}_{n}) = 
		- \left(n (a+b)-a (a+b-1)+n^2+n+1\right)
\\
	 	+ \frac{S_{n-1}}{S_{n}}\frac{\pi_{n-1}}{\pi_{n}} \left(X_{n,n}+b+n-1\right) + \frac{S_{n}}{S_{n+1}}\frac{\pi_{n+1}}{\pi_{n}} \left(Y_{n,n}+a+n+2\right) 	 	
\\
		- \frac{S_{n-1}^2}{S_{n}^2}\frac{\pi_{n-1}^2}{\pi_{n}^2}
		+ 2\frac{S_{n-1}}{S_{n+1}}\frac{\pi_{n-1}\pi_{n+1}}{\pi_{n}^2}  
		- \frac{S_{n}^2}{S_{n+1}^2}\frac{\pi_{n+1}^2}{\pi_{n}^2} 
\\
	+\frac{\xi s^{a+1}e^{-s}}{\pi_{n}\eta_{n}}P_{n+1}(s)
	 \left[
	 	\frac{S_{n}^2}{S_{n+1}^2}Q^{(1)}_{n+1}(-s) - \frac{S_{n}}{S_{n+1}}(X_{n,n}-s)Q^{(1)}_{n}(-s) - \frac{S_{n-1}}{S_{n+1} }Q^{(1)}_{n-1}(-s) 
	 \right]
\\
	+\frac{\xi s^{a+1}e^{-s}}{\pi_{n}\eta_{n}} P_{n-1}(s)
	 \left[
	 	- \frac{S_{n-1}}{S_{n+1}}Q^{(1)}_{n+1}(-s) - \frac{S_{n-1}}{S_{n}}(Y_{n,n}+s)Q^{(1)}_{n}(-s) + \frac{S_{n-1}^2}{S_{n}^2}Q^{(1)}_{n-1}(-s)
	 \right]
\\
	+\frac{\psi t^{b+1}e^{-t}}{\pi_{n}\eta_{n}}P^{(1)}_{n+1}(-t)
	 \left[
	 	\frac{S_{n}^2}{S_{n+1}^2}Q_{n+1}(t) - \frac{S_{n}}{S_{n+1}}(X_{n,n}+t)Q_{n}(t) - \frac{S_{n-1}}{S_{n+1} }Q_{n-1}(t) 
	 \right]
\\
	+\frac{\psi t^{b+1}e^{-t}}{\pi_{n}\eta_{n}} P^{(1)}_{n-1}(-t)
	 \left[
	 	- \frac{S_{n-1}}{S_{n+1}}Q_{n+1}(t) - \frac{S_{n-1}}{S_{n}}(Y_{n,n}-t)Q_{n}(t) + \frac{S_{n-1}^2}{S_{n}^2}Q_{n-1}(t)
	 \right] ,
\end{multline*}
and its determinant
\begin{multline*}
	\det(\mathcal{A}^{(\Sigma)}_{n}) = (a+1) (n+1) (a+b+n)
\\	
		- (a+b+n)\frac{S_{n}}{S_{n+1}}\frac{\pi_{n+1}}{\pi_{n}} (Y_{n,n}+a+n+2) + (n+1)\frac{S_{n-1}}{S_{n}}\frac{\pi_{n-1}}{\pi_{n}} (X_{n,n}+b+n-1)
\\
	 + (a+b+n)\frac{S_{n}^2}{S_{n+1}^2}\frac{\pi_{n+1}^2}{\pi_{n}^2}
	 - (a+b-1)\frac{S_{n-1}}{S_{n+1}}\frac{\pi_{n+1}\pi_{n-1}}{\pi_{n}^2} 
	 - (n+1)\frac{S_{n-1}^2}{S_{n}^2} \frac{\pi_{n-1}^2}{\pi_{n}^2} 
\\
	- (a+b+n) \frac{\xi s^{a+1}e^{-s}}{\pi_{n}\eta_{n}} P_{n+1}(s)
		 \left[ \frac{S_{n}^2}{ S_{n+1}^2 }Q^{(1)}_{n+1}(-s) - \frac{S_{n}}{S_{n+1}}(X_{n,n}-s)Q^{(1)}_{n}(-s) - \frac{S_{n-1}}{S_{n+1}}Q^{(1)}_{n-1}(-s) \right]
\\
	+ (n+1) \frac{\xi s^{a+1}e^{-s}}{\pi_{n}\eta_{n}} P_{n-1}(s)
		 \left[ - \frac{S_{n-1}}{S_{n+1}}Q^{(1)}_{n+1}(-s) - \frac{S_{n-1}}{S_{n}}(Y_{n,n}+s)Q^{(1)}_{n}(-s) + \frac{S_{n-1}^2}{S_{n}^2}Q^{(1)}_{n-1}(-s) \right]
\\
	- (a+b+n) \frac{\psi t^{b+1}e^{-t}}{\pi_{n}\eta_{n}} P^{(1)}_{n+1}(-t)
		 \left[ \frac{ S_{n}^2}{ S_{n+1}^2 }Q_{n+1}(t) - \frac{S_{n}}{S_{n+1}}(X_{n,n}+t)Q_{n}(t) - \frac{S_{n-1}}{S_{n+1}}Q_{n-1}(t) \right]
\\
	+ (n+1) \frac{\psi t^{b+1}e^{-t}}{\pi_{n}\eta_{n}} P^{(1)}_{n-1}(-t)
		 \left[ - \frac{S_{n-1}}{S_{n+1}}Q_{n+1}(t) - \frac{S_{n-1}}{S_{n}}(Y_{n,n}-t)Q_{n}(t) + \frac{S_{n-1}^2}{S_{n}^2}Q_{n-1}(t) \right] .
\\	
\end{multline*}
\end{proposition}
Following up on the invariants of the $\mathcal{A}^{(\cdot)}_{n}$ residue matrices we report evaluations of the traces of all distinct pair-wise products of these.
\begin{proposition}
The traces of pairwise products of residue matrices are given by
\begin{gather*}
	\Tr(\mathcal{A}^{(\infty)}_{n}\mathcal{A}^{(s)}_{n}) = -\xi s^{a}e^{-s}\frac{P_{n}(s)}{\pi^2_{n}\eta_{n}}(\pi_{n+1},\pi_{n},\pi_{n-1})G_{n}(s,-s)\mathrm{Q}^{(1)}_{n}(-s) ,
\\
	\Tr(\mathcal{A}^{(\infty)}_{n}\mathcal{A}^{(-t)}_{n}) = \psi t^{b}e^{-t}\frac{P^{(1)}_{n}(-t)}{\pi_{n}^2\eta_{n}}(\pi_{n+1},\pi_{n},\pi_{n-1})G_{n}(-t,t)\mathrm{Q}_{n}(t) ,
\\
	\Tr(\mathcal{A}^{(\infty)}_{n}\mathcal{A}^{(\Sigma)}_{n}) = -a-1,
\end{gather*}
along with
\begin{equation*}
	\Tr(\mathcal{A}^{(s)}_{n}\mathcal{A}^{(-t)}_{n}) 
	= -\xi\psi s^{a}e^{-s}t^{b}e^{-t}\frac{1}{\pi^2_{n}\eta^2_{n}} 
		\mathrm{P}^{(1)}_{n}(-t)^{T}G_{n}(s,-s)\mathrm{Q}^{(1)}_{n}(-s) \times \mathrm{P}_{n}(s)^{T}G_{n}(-t,t)\mathrm{Q}_{n}(t) .
\end{equation*}
The remaining pair have the evaluations of
\begin{equation}
	\Tr(\mathcal{A}^{(\Sigma)}_{n}\mathcal{A}^{(s)}_{n}) = 
	-\xi s^{a}e^{-s}\frac{1}{\pi_{n}\eta_{n}} \mathrm{P}_{n}(s)^{T} \left(\mathcal{A}^{(\Sigma)}_{n}\right)^{T}G_{n}(s,-s) \mathrm{Q}^{(1)}_{n}(-s),
\end{equation}
and
\begin{equation}
	\Tr(\mathcal{A}^{(\Sigma)}_{n}\mathcal{A}^{(-t)}_{n}) = 
	\psi t^{b}e^{-t}\frac{1}{\pi_{n}\eta_{n}}\mathrm{P}^{(1)}_{n}(-t)^{T} \left(\mathcal{A}^{(\Sigma)}_{n}\right)^{T}G_{n}(-t,t) \mathrm{Q}_{n}(t).
\end{equation}
\end{proposition}

We note that the $\sigma$-functions of Prop. \ref{sigma-form} have been evaluated as the $ 0,1 $ and $1,0 $ kernels with anti-incidence arguments.
However these arguments are set at the singular points $ x=s,-t $ and as one can see from \eqref{anti-incident_K01} and \eqref{anti-incident_K10} that a division
is carried out by $ W_1(x) $ with a simple zero at any singularity. Thus we need to carry out a limiting procedure in order to evaluate Prop. \ref{sigma-form}.
The first step in this direction is to evaluate derivatives at $ x=s, -t $.
\begin{lemma}
Let us assume $ \mathcal{P}^{(\nu)}_{n} \in C^{2}(0,\infty) $ for $ \nu \in \{0,1,2\} $.
The spectral derivative $ \partial_{x}\mathcal{P}^{(\nu)}_{n}(s) $ at the singular point $ x=s $ is given by
\begin{equation}
	\partial_{x}\mathcal{P}^{(\nu)}_{n}(s) = 
	\left[ \mathbb{1}_3 + \mathcal{A}^{(s)}_{n} \right]
	\left( \mathcal{A}^{(\infty)}_{n}+s^{-1}\mathcal{A}^{(0)}_{n}+(s+t)^{-1}\mathcal{A}^{(-t)}_{n} \right) \mathcal{P}^{(\nu)}_{n}(s) ,
\label{DerivPmatrix@s}
\end{equation}
and the corresponding derivative at $ x=-t $ is
\begin{equation}
	\partial_{x}\mathcal{P}^{(\nu)}_{n}(-t) = 
	\left[ \mathbb{1}_3 + \mathcal{A}^{(-t)}_{n} \right]
	\left( \mathcal{A}^{(\infty)}_{n}-t^{-1}\mathcal{A}^{(0)}_{n}-(s+t)^{-1}\mathcal{A}^{(s)}_{n} \right) \mathcal{P}^{(\nu)}_{n}(-t) .
\label{DerivPmatrix@-t}
\end{equation}
One can interpret the coefficient of $ \mathcal{P}^{(\nu)}_{n}(s) $ on the right-hand side of \eqref{DerivPmatrix@s} as $ \mathcal{A}_{n}(s) $, 
and similarly for \eqref{DerivPmatrix@-t}. 
Corresponding formulae hold for $ \mathcal{Q}^{(\nu)}_{n} $.
\end{lemma}
\begin{proof}
We give the proof for the first case only as the treatment of the second is the same.
Separating the partial fraction form of \eqref{xst-LaxPairs} and utilising \eqref{AresidueSing@s} we can write
\begin{multline*}
	\partial_{x}\mathcal{P}^{(\nu)}_{n}(x) = 
	\left( \mathcal{A}^{(\infty)}_{n}+x^{-1}\mathcal{A}^{(0)}_{n}+(x+t)^{-1}\mathcal{A}^{(-t)}_{n} \right)\mathcal{P}^{(\nu)}_{n}(x)
	+ (x-s)^{-1}a^{(s)}_{n}\mathcal{P}^{(0)}_{n}(s)\otimes \left(G_{n}(s,-s)\mathcal{Q}^{(1)}_{n}(-s)\right)^{T}\mathcal{P}^{(\nu)}_{n}(x) ,
\end{multline*}
where we abbreviate $ a^{(s)} = -(-)^{b}\xi s^{a+b}/\pi_{n}\eta_{n} $.
Expanding $ x=s+\epsilon $ we find that in the second term of the right-hand side there is a within it a single term of order $ \epsilon^{-1} $ having coefficient 
\begin{equation*}
	a^{(s)}_{n}\mathcal{P}^{(0)}_{n}(s)\otimes \left( G_{n}(s,-s)\mathcal{Q}^{(1)}_{n}(-s) \right)^{T}\mathcal{P}^{(\nu)}_{n}(s)
	= a^{(s)}_{n}\mathcal{P}^{(0)}_{n}(s)\cdot \mathcal{P}^{(\nu)}_{n}(s)^{T}G_{n}(s,-s)\mathcal{Q}^{(1)}_{n}(-s)
	= 0 ,
\end{equation*}
and vanishing by orthogonality. Thus in the limit $ \epsilon\to 0 $ we deduce
\begin{equation*}
	\left[ \mathbb{1}_3 - a^{(s)}_{n}\mathcal{P}^{(0)}_{n}(s)\otimes \left( G_{n}(s,-s)\mathcal{Q}^{(1)}_{n}(-s) \right)^{T} \right]\mathcal{P}^{(\nu)\prime}_{n}(s)
	= \left( \mathcal{A}^{(\infty)}_{n}+s^{-1}\mathcal{A}^{(0)}_{n}+(s+t)^{-1}\mathcal{A}^{(-t)}_{n} \right)\mathcal{P}^{(\nu)}_{n}(s) .
\end{equation*}
For arbitrary column vectors $A,B$ subject to $ A^{T}B=0 $ we note that the matrix on the left-hand side has determinant unity and the inverse can be performed using
$ \left[ \mathbb{1}_3 - A\otimes B^{T} \right]^{-1} = \mathbb{1}_3 + A\otimes B^{T} $, and consequently \eqref{DerivPmatrix@s} follows.
\end{proof}

\begin{proposition}
Let generic conditions apply and $n\geq 1$.
The anti-incidence limit of the kernel $K^{0,1}_{n}$, at the singular point $ x=s $ is
\begin{equation}
	\pi_{n}\eta_{n} K^{0,1}_{n}(s,-s)
\\
	=\mathrm{P}^{(0)}_{n}(s)^{T} 
	\left\{
	\begin{pmatrix}
		0 & 0 & 0 \\ \frac{\pi_{n+1}}{\pi_{n}} & 1 & \frac{\pi_{n-1}}{\pi_{n}}+\frac{S_{n}}{S_{n-1}} \\ 0 & 0 & 0
	\end{pmatrix}
	+ s^{-1}\mathcal{A}^{(\Sigma)}_{n}{}^{T}
	- \frac{t}{s(s+t)}\mathcal{A}^{(-t)}_{n}{}^{T}
	\right\} G_{n}(s,-s)\mathrm{Q}^{(1)}_{n}(-s) ,
\label{K01_Limit@s}
\end{equation}
and that of the kernel $K^{1,0}_{n}$, at the singular point $ x=-t $ is
\begin{equation}
	\pi_{n}\eta_{n} K^{1,0}_{n}(-t,t) 
\\
	=\mathrm{P}^{(1)}_{n}(-t)^{T} 
	\left\{
	\begin{pmatrix}
		0 & 0 & 0 \\ \frac{\pi_{n+1}}{\pi_{n}} & 1 & \frac{\pi_{n-1}}{\pi_{n}}+\frac{S_{n}}{S_{n-1}} \\ 0 & 0 & 0
	\end{pmatrix}
	- t^{-1}\mathcal{A}^{(\Sigma)}_{n}{}^{T}
	+ \frac{s}{t(s+t)}\mathcal{A}^{(s)}_{n}{}^{T}
	\right\} G_{n}(-t,t)\mathrm{Q}^{(0)}_{n}(t) .
\label{K10_Limit@-t}
\end{equation}
\end{proposition}
\begin{proof}
We start with the first form of \eqref{anti-incident_K01} and writing the spectral matrix in partial fraction form we can split off the term
with the pole at $x=s$ from the remainder thus,
\begin{multline*}
	\pi_{n}\eta_{n} K^{0,1}_{n}(x,-x) 
\\
	= (-x)^{b}e^{x}
	\mathcal{P}^{(0)}_{n}(x){}^{T}\left\{ \left( \mathcal{A}^{(\infty)}_{n} +x^{-1}\mathcal{A}^{(0)}_{n}+(x+t)^{-1}\mathcal{A}^{(-t)}_{n} \right)^{T}G_{n}(x,-x)
	+ \partial_{x}G_{n}(x,-x) \right\} \mathcal{Q}^{(1)}_{n}(-x)
\\
	+ (x-s)^{-1}(-x)^{b}e^{x}\mathcal{P}^{(0)}_{n}(x){}^{T} \mathcal{A}^{(s)}_{n}{}^{T}G_{n}(x,-x)\mathcal{Q}^{(1)}_{n}(-x) .
\end{multline*}
The first term on the right-hand side can be directly evaluated as $ x\to s $ along with
\begin{equation*}
	\partial_{x}G_{n}(s,-s) =
	\begin{pmatrix}
	0 & 0 & 0 \\ -\frac{S_{n}}{S_{n+1}} & X_{n,n}-s & \frac{S_{n-1}}{S_{n}} \\ 0 & 0 & 0 
	\end{pmatrix} .
\end{equation*}
In the second term we set $ x=s+\epsilon $ and expand up to order $ \epsilon $ which yields
\begin{multline*}
	\epsilon^{-1} (-s)^{b}e^{s}\mathcal{P}^{(0)}_{n}(s){}^{T} \mathcal{A}^{(s)}_{n}{}^{T}G_{n}(s,-s)\mathcal{Q}^{(1)}_{n}(-s)
	+ (-s)^{b}e^{s}\frac{b-s}{s}\mathcal{P}^{(0)}_{n}(s){}^{T} \mathcal{A}^{(s)}_{n}{}^{T}G_{n}(s,-s)\mathcal{Q}^{(1)}_{n}(-s) 
\\
	+ (-s)^{b}e^{s}\mathcal{P}^{(0)\prime}_{n}(s){}^{T} \mathcal{A}^{(s)}_{n}{}^{T}G_{n}(s,-s)\mathcal{Q}^{(1)}_{n}(-s)
	+ (-s)^{b}e^{s}\mathcal{P}^{(0)}_{n}(s){}^{T} \mathcal{A}^{(s)}_{n}{}^{T}(\partial_{x}G_{n}-\partial_{y}G_{n})(s,-s)\mathcal{Q}^{(1)}_{n}(-s)
\\
	- (-s)^{b}e^{s}\mathcal{P}^{(0)}_{n}(s){}^{T} \mathcal{A}^{(s)}_{n}{}^{T}G_{n}(s,-s)\mathcal{Q}^{(1)\prime}_{n}(-s)
	+ {\rm O}(\epsilon) .
\end{multline*}
From the direct product structure of \eqref{AresidueSing@s} we see that
\begin{equation*}
	\mathcal{P}^{(0)}_{n}(s){}^{T} \mathcal{A}^{(s)}_{n}{}^{T}G_{n}(s,-s)\mathcal{Q}^{(1)}_{n}(-s) = 
	a^{(s)}_{n} \left( \mathcal{P}^{(0)}_{n}(s){}^{T}G_{n}(s,-s)\mathcal{Q}^{(1)}_{n}(-s) \right)^2 = 0 ,
\end{equation*}
due to \eqref{symmetric_kernel}. Thus we can take the limit $ \epsilon\to 0 $.
Combining all the ${\rm O}(1)$ terms, employing the identity \eqref{D-A-relation} (relating $ \mathcal{D}_{n}(-s) $ to $ \mathcal{A}_{n}(s) $) we can effect some cancellations
and finally substituting for $ \mathcal{A}_{n}(s) $ using \eqref{DerivPmatrix@s} we arrive at
\begin{multline*}
	\pi_{n}\eta_{n} K^{0,1}_{n}(s,-s) 
\\
	= (-s)^{b}e^{s}
	\mathcal{P}^{(0)}_{n}(s){}^{T} \Bigg\{
	\partial_{x}G_{n}(s,-s)\cdot G_{n}(s,-s)^{-1} - \frac{a+s}{s}\mathcal{A}^{(s)T}_{n} + \mathcal{A}^{(\infty)T}_{n}+s^{-1}\mathcal{A}^{(0)T}_{n}+(s+t)^{-1}\mathcal{A}^{(-t)T}_{n}
\\
	+ \left[ \mathcal{A}^{(\infty)T}_{n}+s^{-1}\mathcal{A}^{(0)T}_{n}+(s+t)^{-1}\mathcal{A}^{(-t)T}_{n}, \mathcal{A}^{(s)T}_{n} \right] 
	\left( \mathbb{1}_3 + \mathcal{A}^{(s)T}_{n}\right)
	\Bigg\} G_{n}(s,-s)\mathcal{Q}^{(1)}_{n}(-s) .
\end{multline*}
Next we use the identities noted above, and find that a number of terms drop out.
\end{proof}

\subsection{Constrained Dynamical System}\label{Constrained_Dynamics}

Our dynamical system will consist of $12$ independent, primary polynomial and associated function variables
$\{P_{n+\delta}(s)\}_{\delta=1,0,-1}$, 
$\{P^{(1)}_{n+\delta}(-t)\}_{\delta=1,0,-1}$, 
$\{Q_{n+\delta}(t)\}_{\delta=1,0,-1}$, 
$\{Q^{(1)}_{n+\delta}(-s)\}_{\delta=1,0,-1}$,  
along with $6$ parameters
$\{\pi_{n+\delta}\}_{\delta=1,0,-1}$, $\{\eta_{n+\delta}\}_{\delta=1,0,-1}$, 
the recurrence coefficients $X_{n,n}$, $Y_{n,n}$ 
and the norming coefficients $\{S_{n+\delta}\}_{\delta=1,0,-1}$.
Out of the total of $23$ variables we have eight constraints, the first two of which follow from \eqref{symmetric_kernel}
\begin{align}
	\mathrm{P}_{n}^{(0)}(s)^{T} G_{n}(s,-s) \mathrm{Q}_{n}^{(1)}(-s) & = 0 ,
\label{Bilinear_01}
\\
	\mathrm{P}_{n}^{(1)}(-t)^{T} G_{n}(-t,t) \mathrm{Q}_{n}^{(0)}(t) & = 0 ,	
\label{Bilinear_10}
\end{align}
with a further six given below.  
\begin{proposition}\label{constraints}
For all $ n\geq 0 $, $ s,t $, $ \xi,\psi $ and $ a,b $ the following constraints apply:
\begin{equation}
	X_{n,n}+Y_{n,n} = \pi_{n}\eta_{n} = 2n+a+b+1 + \xi s^{a+1}e^{-s}P_{n}(s)Q^{(1)}_{n}(-s) + \psi t^{b+1}e^{-t}P^{(1)}_{n}(-t)Q_{n}(t) ,
\label{00_constraint} 
\end{equation}
\begin{multline}
	X_{n,n} = \frac{S_{n}}{S_{n+1}}\pi_{n+1}\eta_{n} - \frac{S_{n-1}}{S_{n}}\pi_{n}\eta_{n-1}
\\
	- \xi s^{a+1}e^{-s}
	\left[ \frac{S_{n}}{S_{n+1}}P_{n+1}(s)Q^{(1)}_{n}(-s) - \frac{S_{n-1}}{S_{n}}P_{n}(s)Q^{(1)}_{n-1}(-s)  \right]
\\
	- \psi t^{b+1}e^{-t}
	\left[ \frac{S_{n}}{S_{n+1}}P^{(1)}_{n+1}(-t)Q_{n}(t) - \frac{S_{n-1}}{S_{n}}P^{(1)}_{n}(-t)Q_{n-1}(t)  \right] ,
\label{X_constraint}	
\end{multline}
\begin{multline}
	Y_{n,n} = \frac{S_{n}}{S_{n+1}}\pi_{n}\eta_{n+1} - \frac{S_{n-1}}{S_{n}}\pi_{n-1}\eta_{n}
\\
	- \xi s^{a+1}e^{-s}
	\left[ \frac{S_{n}}{S_{n+1}}P_{n}(s)Q^{(1)}_{n+1}(-s) - \frac{S_{n-1}}{S_{n}}P_{n-1}(s)Q^{(1)}_{n}(-s)  \right]
\\
	- \psi t^{b+1}e^{-t}
	\left[ \frac{S_{n}}{S_{n+1}}P^{(1)}_{n}(-t)Q_{n+1}(t) - \frac{S_{n-1}}{S_{n}}P^{(1)}_{n-1}(-t)Q_{n}(t)  \right] ,
\label{Y_constraint}	
\end{multline}
\begin{multline}
	X_{n,n} -n-a - \frac{S_{n}\eta_{n+1}}{S_{n+1}\eta_{n}} + \frac{S_{n-1}\eta_{n-1}}{S_{n}\eta_{n}}
	\\
	+ \xi s^{a+1}e^{-s}\frac{1}{\pi_{n}\eta_{n}}Q^{(1)}_{n}(-s) 
	\left[ \frac{S_{n}}{S_{n+1}}P_{n+1}(s)-(Y_{n,n}+s)P_{n}(s)-\frac{S_{n-1}}{S_{n}}P_{n-1}(s) \right]
	\\
	+ \psi t^{b+1}e^{-t}\frac{1}{\pi_{n}\eta_{n}}Q_{n}(t) 
	\left[ \frac{S_{n}}{S_{n+1}}P^{(1)}_{n+1}(-t)-(Y_{n,n}-t)P^{(1)}_{n}(-t)-\frac{S_{n-1}}{S_{n}}P^{(1)}_{n-1}(-t) \right] = 0 ,
\label{XY_constraint} 
\end{multline}
and
\begin{multline}
	Y_{n,n} -n-b - \frac{S_{n}\pi_{n+1}}{S_{n+1}\pi_{n}} + \frac{S_{n-1}\pi_{n-1}}{S_{n}\pi_{n}}
	\\
	+ \xi s^{a+1}e^{-s}\frac{1}{\pi_{n}\eta_{n}}P_{n}(s) 
	\left[ \frac{S_{n}}{S_{n+1}}Q^{(1)}_{n+1}(-s)-(X_{n,n}-s)Q^{(1)}_{n}(-s)-\frac{S_{n-1}}{S_{n}}Q^{(1)}_{n-1}(-s) \right]
	\\
	+ \psi t^{b+1}e^{-t}\frac{1}{\pi_{n}\eta_{n}}P^{(1)}_{n}(-t) 
	\left[ \frac{S_{n}}{S_{n+1}}Q_{n+1}(t)-(X_{n,n}+t)Q_{n}(t)-\frac{S_{n-1}}{S_{n}}Q_{n-1}(t) \right] = 0 .
\label{YX_constraint} 
\end{multline}
Note that not all the constraints given above are independent: for example employing \eqref{X_constraint} along with \eqref{Y_constraint} in \eqref{XY_constraint},
and combining the latter result with \eqref{YX_constraint} directly implies \eqref{00_constraint}.
Or alternatively the system of four relations \eqref{X_constraint}--\eqref{YX_constraint}, 
when viewed as a linear inhomogeneous system in $\pi_{n+1},\pi_{n-1},\eta_{n+1},\eta_{n-1}$, has only rank three.
\end{proposition}
\begin{proof}
Setting $ m=n $ in the identity \eqref{sum_BOP-derivatives} and noting $ \langle x\partial_x P_{n}, Q_{n} \rangle = \langle P_{n}, y\partial_y Q_{n} \rangle = n $
we establish \eqref{00_constraint}.
Note that the right-hand sides of \eqref{X_constraint} and \eqref{Y_constraint} are perfect differences in $n$.
One can set $ n\mapsto n+1 $ in \eqref{sum_BOP-derivatives} and $ m=n $, and because $ \langle P_{n+1}, y\partial_y Q_{n} \rangle = 0 $ we have
\begin{equation}
	\langle x\partial_x P_{n+1}, Q_{n} \rangle 
		= \pi_{n+1}\eta_{n} - \xi s^{a+1}e^{-s}P_{n+1}(s)Q^{(1)}_{n}(-s) - \psi t^{b+1}e^{-t}P^{(1)}_{n+1}(-t)Q_{n}(t).
\end{equation}
Now if one expands $ x\partial_x P_{n+1} $ in a $P$-basis thus
\begin{equation}
	x\partial_x P_{n+1}(x) = (n+1)P_{n+1}(x) - \frac{S_{n+1,n}}{S_{n}}P_{n}(x) + \Pi_{n-1}[x],
\end{equation}
and so $ \langle x\partial_x P_{n+1}, Q_{n} \rangle = -S_{n+1,n}/S_{n} $. Consequently
\begin{equation}
	-\frac{S_{n+1,n}}{S_{n}} = \pi_{n+1}\eta_{n} - \xi s^{a+1}e^{-s}P_{n+1}(s)Q^{(1)}_{n}(-s) - \psi t^{b+1}e^{-t}P^{(1)}_{n+1}(-t)Q_{n}(t) .
\end{equation}
Then \eqref{X_constraint} follows by differencing this and using the relation \eqref{P_subleading}.
The constraints \eqref{XY_constraint} and \eqref{YX_constraint} most readily deduced from the traces of $ \mathcal{D}^{(0)}_{n} $, $ \mathcal{A}^{(0)}_{n} $,
after employing identity \eqref{00_constraint} with $ n\mapsto n+1 $.
\end{proof}

As our main results we give a complete list of the fundamental differential equations with respect to $s,t$ satisfied by the 
auxiliary parameters $ \pi_{n}, \eta_{n} $, $ X_{n,n}, Y_{n,n} $,
the specialised polynomials $ P_{n}(s), Q_{n}(t) $ and
the specialised associated functions $ P^{(1)}_{n}(-t), Q^{(1)}_{n}(-s) $.
\begin{proposition}\label{PQ_partial-derivatives}
The partial derivatives of the polynomials and associated functions with respect to the deformations, and then specialised, are given by
\begin{equation}
	\partial_{s} 
		\begin{pmatrix}
			P_{n+1}(s) \\ P_{n}(s) \\ P_{n-1}(s) 
		\end{pmatrix} = 
		\left[ \mathcal{B}^{(\infty,0)}_{n} + \xi s^{a}e^{-s} K_{n}^{(0,1)}(s,-s)\mathbb{1}_3 \right]
		\begin{pmatrix}
			P_{n+1}(s) \\ P_{n}(s) \\ P_{n-1}(s) 
		\end{pmatrix} ,
\label{s-PartialDeriv-P0Ps}
\end{equation}
\begin{equation}
	\partial_{t} 
		\begin{pmatrix}
			P_{n+1}(s) \\ P_{n}(s) \\ P_{n-1}(s) 
		\end{pmatrix} = 
		\mathcal{C}^{(\infty,0)}_{n}
		\begin{pmatrix}
			P_{n+1}(s) \\ P_{n}(s) \\ P_{n-1}(s) 
		\end{pmatrix}
		+ \psi t^{b}e^{-t} K_{n}^{(0,0)}(s,t)
		\begin{pmatrix}
			P^{(1)}_{n+1}(-t) \\ P^{(1)}_{n}(-t) \\ P^{(1)}_{n-1}(-t) 
		\end{pmatrix} ,
\label{t-PartialDeriv-P0Ps}
\end{equation}
\begin{equation}
	\partial_{s} 
		\begin{pmatrix}
			Q_{n+1}(t) \\ Q_{n}(t) \\ Q_{n-1}(t) 
		\end{pmatrix} = 
		\mathcal{B}^{(\infty,\bar{0})}_{n}
		\begin{pmatrix}
			Q_{n+1}(t) \\ Q_{n}(t) \\ Q_{n-1}(t) 
		\end{pmatrix}
		+ \xi s^{a}e^{-s} K_{n}^{(0,0)}(s,t)
		\begin{pmatrix}
			Q^{(1)}_{n+1}(-s) \\ Q^{(1)}_{n}(-s) \\ Q^{(1)}_{n-1}(-s) 
		\end{pmatrix} ,
\label{s-PartialDeriv-Q0Pt}
\end{equation}
\begin{equation}
	\partial_{t} 
		\begin{pmatrix}
			Q_{n+1}(t) \\ Q_{n}(t) \\ Q_{n-1}(t) 
		\end{pmatrix} = 
		\left[ \mathcal{C}^{(\infty,\bar{0})}_{n} + \psi t^{b}e^{-t} K_{n}^{(1,0)}(-t,t)\mathbb{1}_3 \right]
		\begin{pmatrix}
			Q_{n+1}(t) \\ Q_{n}(t) \\ Q_{n-1}(t) 
		\end{pmatrix} ,
\label{t-PartialDeriv-Q0Pt}
\end{equation}
\begin{equation}
	\partial_{s} 
		\begin{pmatrix}
			P^{(1)}_{n+1}(-t) \\ P^{(1)}_{n}(-t) \\ P^{(1)}_{n-1}(-t) 
		\end{pmatrix} = 
		\mathcal{B}^{(\infty,0)}_{n}
		\begin{pmatrix}
			P^{(1)}_{n+1}(-t) \\ P^{(1)}_{n}(-t) \\ P^{(1)}_{n-1}(-t) 
		\end{pmatrix}
		+ \xi s^{a}e^{-s} K_{n}^{(1,1)}(-t,-s)
		\begin{pmatrix}
			P_{n+1}(s) \\ P_{n}(s) \\ P_{n-1}(s) 
		\end{pmatrix} ,
\label{s-PartialDeriv-P1Mt}
\end{equation}
\begin{equation}
	\partial_{t} 
		\begin{pmatrix}
			P^{(1)}_{n+1}(-t) \\ P^{(1)}_{n}(-t) \\ P^{(1)}_{n-1}(-t) 
		\end{pmatrix} = 
		\left[ \mathcal{C}^{(\infty,0)}_{n} + \psi t^{b}e^{-t} K_{n}^{(1,0)}(-t,t)\mathbb{1}_3 \right]
		\begin{pmatrix}
			P^{(1)}_{n+1}(-t) \\ P^{(1)}_{n}(-t) \\ P^{(1)}_{n-1}(-t) 
		\end{pmatrix} ,
\label{t-PartialDeriv-P1Mt}
\end{equation}
\begin{equation}
	\partial_{s} 
		\begin{pmatrix}
			Q^{(1)}_{n+1}(-s) \\ Q^{(1)}_{n}(-s) \\ Q^{(1)}_{n-1}(-s) 
		\end{pmatrix} = 
		\left[ \mathcal{B}^{(\infty,\bar{0})}_{n} + \xi s^{a}e^{-s} K_{n}^{(0,1)}(s,-s)\mathbb{1}_3 \right]
		\begin{pmatrix}
			Q^{(1)}_{n+1}(-s) \\ Q^{(1)}_{n}(-s) \\ Q^{(1)}_{n-1}(-s) 
		\end{pmatrix} ,
\label{s-PartialDeriv-Q1Ms}
\end{equation}
\begin{equation}
	\partial_{t} 
		\begin{pmatrix}
			Q^{(1)}_{n+1}(-s) \\ Q^{(1)}_{n}(-s) \\ Q^{(1)}_{n-1}(-s) 
		\end{pmatrix} = 
		\mathcal{C}^{(\infty,\bar{0})}_{n}
		\begin{pmatrix}
			Q^{(1)}_{n+1}(-s) \\ Q^{(1)}_{n}(-s) \\ Q^{(1)}_{n-1}(-s) 
		\end{pmatrix}
		+ \psi t^{b}e^{-t} K_{n}^{(1,1)}(-t,-s)
		\begin{pmatrix}
			Q_{n+1}(t) \\ Q_{n}(t) \\ Q_{n-1}(t) 
		\end{pmatrix} .
\label{t-PartialDeriv-Q1Ms}
\end{equation}
\end{proposition}
\begin{proof}
All of these follow from the sets of derivatives \eqref{Pdiff-s}-\eqref{Qdiff-t} and \eqref{P1diff-s}-\eqref{Q1diff-t},
making the substitutions $ n\mapsto n\pm 1$ where necessary 
and using the recurrence relations \eqref{P_recur}, \eqref{Q_recur} to shift the indices $n\pm 2$ back into the range $\{n-1,n,n+1\}$.
\end{proof}

Furthermore define the auxiliary matrices
\begin{equation}
	\mathcal{A}^{(0,+)}_{n} :=
		\begin{pmatrix}
			n+1-\frac{S_{n}\pi_{n+1}}{S_{n+1}\pi_{n}}	& \frac{\pi_{n+1}}{\pi_{n}}(Y_{n,n}+s)												& \frac{S_{n-1}\pi_{n+1}}{S_{n}\pi_{n}}	\\
			-\frac{S_{n}}{S_{n+1}}						& Y_{n,n}+s-n-a-b-1																	& \frac{S_{n-1}}{S_{n}}	\\
			-\frac{S_{n}\pi_{n-1}}{S_{n+1}\pi_{n}}		& -\frac{\pi_{n-1}}{\pi_{n}}(X_{n,n}-s)	+ \psi t^{b+1}e^{-t}P^{(1)}_{n-1}(-t)Q_{n}(t)	& \frac{S_{n-1}\pi_{n-1}}{S_{n}\pi_{n}}-n-a-b
		\end{pmatrix} ,
\end{equation}
\begin{equation}
	\mathcal{A}^{(0,-)}_{n} :=
		\begin{pmatrix}
			n+1+a+t-\frac{S_{n}\pi_{n+1}}{S_{n+1}\pi_{n}}	& \frac{\pi_{n+1}}{\pi_{n}}(Y_{n,n}-t)											& \frac{S_{n-1}\pi_{n+1}}{S_{n}\pi_{n}}	\\
			-\frac{S_{n}}{S_{n+1}}						    & Y_{n,n}-n-b-1																	& \frac{S_{n-1}}{S_{n}}	\\
			-\frac{S_{n}\pi_{n-1}}{S_{n+1}\pi_{n}}		& -\frac{\pi_{n-1}}{\pi_{n}}(X_{n,n}+t)	+ \xi s^{a+1}e^{-s}P_{n-1}(s)Q^{(1)}_{n}(-s)	& \frac{S_{n-1}\pi_{n-1}}{S_{n}\pi_{n}}-n-b+t
		\end{pmatrix} ,
\end{equation}
\begin{equation}
	\mathcal{D}^{(0,+)}_{n} :=
		\begin{pmatrix}
			n+1-\frac{S_{n}\eta_{n+1}}{S_{n+1}\eta_{n}}	& \frac{\eta_{n+1}}{\eta_{n}}(X_{n,n}+t)											& \frac{S_{n-1}\eta_{n+1}}{S_{n}\eta_{n}}	\\
			-\frac{S_{n}}{S_{n+1}}						& X_{n,n}+t-n-a-b-1																	& \frac{S_{n-1}}{S_{n}}	\\
			-\frac{S_{n}\eta_{n-1}}{S_{n+1}\eta_{n}}	& -\frac{\eta_{n-1}}{\eta_{n}}(Y_{n,n}-t) + \xi s^{a+1}e^{-s}P_{n}(s)Q^{(1)}_{n-1}(-s)	& \frac{S_{n-1}\eta_{n-1}}{S_{n}\eta_{n}}-n-a-b
		\end{pmatrix} ,
\end{equation}
\begin{equation}
	\mathcal{D}^{(0,-)}_{n} :=
		\begin{pmatrix}
			n+1+b+s-\frac{S_{n}\eta_{n+1}}{S_{n+1}\eta_{n}}	& \frac{\eta_{n+1}}{\eta_{n}}(X_{n,n}-s)										& \frac{S_{n-1}\eta_{n+1}}{S_{n}\eta_{n}}	\\
			-\frac{S_{n}}{S_{n+1}}						& X_{n,n}-n-a-1																		& \frac{S_{n-1}}{S_{n}}	\\
			-\frac{S_{n}\eta_{n-1}}{S_{n+1}\eta_{n}}	& -\frac{\eta_{n-1}}{\eta_{n}}(Y_{n,n}+s) + \psi t^{b+1}e^{-t}P^{(1)}_{n}(-t)Q_{n-1}(t)& \frac{S_{n-1}\eta_{n-1}}{S_{n}\eta_{n}}+s-n-a
		\end{pmatrix} ,
\end{equation}
\begin{equation}
	\mathcal{A}^{(\diagdown)}_{n} := \tfrac{1}{2}\xi s^{a+1}e^{-s}
		\begin{pmatrix}
			P_{n+1}(s)Q^{(1)}_{n+1}(-s)	& 0								& 0								\\
			0							& -P_{n}(s)Q^{(1)}_{n}(-s)		& 0								\\
			0							& 0								& -P_{n-1}(s)Q^{(1)}_{n-1}(-s)	\\
		\end{pmatrix} ,
\end{equation}
\begin{equation}
	\mathcal{D}^{(\diagdown)}_{n} := \tfrac{1}{2}\psi t^{b+1}e^{-t}
		\begin{pmatrix}
			P^{(1)}_{n+1}(-t)Q_{n+1}(t)	& 0								& 0								\\
			0							& -P^{(1)}_{n}(-t)Q_{n}(t)		& 0								\\
			0							& 0								& -P^{(1)}_{n-1}(-t)Q_{n-1}(t)	\\
		\end{pmatrix} .
\end{equation}

\begin{proposition}\label{PQ_total-derivatives}
The total derivatives of the specialised polynomials and associated functions with respect to the deformations are given by
\begin{equation}
	s\frac{d}{ds}
	\begin{pmatrix}
		P_{n+1}(s) \\ P_{n}(s) \\ P_{n-1}(s) 
	\end{pmatrix} = 
	\left[ \mathcal{A}^{(0,+)}_{n} + \mathcal{A}^{(\diagdown)}_{n} \right]
	\begin{pmatrix}
		P_{n+1}(s) \\ P_{n}(s) \\ P_{n-1}(s) 
	\end{pmatrix} 
	- \psi t^{b+1}e^{-t} K_{n}^{(0,0)}(s,t)
	\begin{pmatrix}
		P^{(1)}_{n+1}(-t) \\ P^{(1)}_{n}(-t) \\ P^{(1)}_{n-1}(-t) 
	\end{pmatrix} ,
\label{s-TotalDeriv-P0Ps}
\end{equation}
\begin{equation}
	t\frac{d}{dt} 
	\begin{pmatrix}
		Q_{n+1}(t) \\ Q_{n}(t) \\ Q_{n-1}(t) 
	\end{pmatrix} = 
	\left[ \mathcal{D}^{(0,+)}_{n} + \mathcal{D}^{(\diagdown)}_{n} \right]
	\begin{pmatrix}
		Q_{n+1}(t) \\ Q_{n}(t) \\ Q_{n-1}(t) 
	\end{pmatrix}
	- \xi s^{a+1}e^{-s} K_{n}^{(0,0)}(s,t)
	\begin{pmatrix}
		Q^{(1)}_{n+1}(-s) \\ Q^{(1)}_{n}(-s) \\ Q^{(1)}_{n-1}(-s) 
	\end{pmatrix} ,
\label{t-TotalDeriv-Q0Pt}
\end{equation}
\begin{equation}
	t\frac{d}{dt}
	\begin{pmatrix}
		P^{(1)}_{n+1}(-t) \\ P^{(1)}_{n}(-t) \\ P^{(1)}_{n-1}(-t) 
	\end{pmatrix} = 
	\left[ \mathcal{A}^{(0,-)}_{n} + \mathcal{D}^{(\diagdown)}_{n} \right]
	\begin{pmatrix}
		P^{(1)}_{n+1}(-t) \\ P^{(1)}_{n}(-t) \\ P^{(1)}_{n-1}(-t) 
	\end{pmatrix}
	- \xi s^{a+1}e^{-s} K_{n}^{(1,1)}(-t,-s)
	\begin{pmatrix}
		P_{n+1}(s) \\ P_{n}(s) \\ P_{n-1}(s) 
	\end{pmatrix}  ,
\label{t-TotalDeriv-P1Mt}
\end{equation}
\begin{equation}
	s\frac{d}{ds} 
	\begin{pmatrix}
			Q^{(1)}_{n+1}(-s) \\ Q^{(1)}_{n}(-s) \\ Q^{(1)}_{n-1}(-s) 
	\end{pmatrix} = 
	\left[ \mathcal{D}^{(0,-)}_{n} + \mathcal{A}^{(\diagdown)}_{n} \right]
	\begin{pmatrix}
			Q^{(1)}_{n+1}(-s) \\ Q^{(1)}_{n}(-s) \\ Q^{(1)}_{n-1}(-s) 
	\end{pmatrix}
	- \psi t^{b+1}e^{-t} K_{n}^{(1,1)}(-t,-s)
	\begin{pmatrix}
		Q_{n+1}(t) \\ Q_{n}(t) \\ Q_{n-1}(t) 
	\end{pmatrix} .
\label{s-TotalDeriv-Q1Ms}
\end{equation}
\end{proposition}
\begin{proof}
For the total derivatives reported we require a sum of the partial derivatives \eqref{Pdiff-s}-\eqref{Qdiff-t} and \eqref{P1diff-s}-\eqref{Q1diff-t},
and the spectral derivatives \eqref{P_diff-x},\eqref{Q_diff-y},\eqref{P1_diff-x},\eqref{Q1_diff-y} with $x$ or $y$ locked to $s,-t$ or $-s,t$ respectively.
The relevant derivatives are the following
\begin{align*}
	s\frac{d}{ds}P_{n}(s) = & -(n+a+b+1) P_{n}(s) - \pi_{n}\hat{P}_{n}(s)
\\ 
	& \phantom{=} - \tfrac{1}{2}\xi s^{a+1}e^{-s} \left( P_{n}(s) \right)^2 Q^{(1)}_{n}(-s) - \psi t^{b+1}e^{-t} P^{(1)}_{n}(-t) K^{0,0}_{n}(s,t) ,
\\
	t\frac{d}{dt}Q_{n}(t) = & -(n+a+b+1) Q_{n}(t) - \eta_{n}\check{Q}_{n}(t)
\\ 
	& \phantom{=} - \xi s^{a+1}e^{-s} Q^{(1)}_{n}(-s)K^{0,0}_{n}(s,t)  - \tfrac{1}{2}\psi t^{b+1}e^{-t} \left( Q_{n}(t)\right)^2 P^{(1)}_{n}(-t) ,
\\
	t\frac{d}{dt}P^{(1)}_{n}(-t) = & -(n+b+1-t) P^{(1)}_{n}(-t) + \pi_{n}\left( 1-\hat{P}^{(1)}_{n}(-t) \right)
\\ 
	& \phantom{=} - \xi s^{a+1}e^{-s} P_{n}(s) K^{1,1}_{n}(-t,-s) - \tfrac{1}{2}\psi t^{b+1}e^{-t} \left( P^{(1)}_{n}(-t) \right)^2 Q_{n}(t) ,
\nonumber\\
	s\frac{d}{ds}Q^{(1)}_{n}(-s) = & -(n+a+1-s) Q^{(1)}_{n}(-s) + \eta_{n}\left( 1-\check{Q}^{(1)}_{n}(-s) \right)
\\ 
	& \phantom{=} - \tfrac{1}{2}\xi s^{a+1}e^{-s} \left( Q^{(1)}_{n}(-s) \right)^2 P_{n}(s) - \psi t^{b+1}e^{-t} Q_{n}(t) K^{1,1}_{n}(-t,-s)  .
\end{align*}
In addition we make the appropriate substitutions for the intertwined polynomials using \eqref{hatPsubs},\eqref{checkQsubs},\eqref{hatP1subs} and \eqref{checkQ1subs}.
In constructing the matrix derivatives we make the substitutions $ n\mapsto n\pm 1$ where necessary 
and use the recurrence relations \eqref{P_recur}, \eqref{Q_recur} to restore the indices $n\pm 2$ back into the range $\{n-1,n,n+1\}$.
\end{proof}

\begin{proposition}\label{PiEta_derivatives}
The deformation derivatives of the auxiliary parameters are given in matrix form by
\begin{multline}
	\partial_{s} 
		\begin{pmatrix}
			\pi_{n+1} \\ \pi_{n} \\ \pi_{n-1} 
		\end{pmatrix} = \mathcal{B}^{(\infty,0)}_{n}
		\begin{pmatrix}
			\pi_{n+1} \\ \pi_{n} \\ \pi_{n-1} 
		\end{pmatrix}
\\
		- \xi s^{a}e^{-s} \frac{1}{\eta_{n}}\left[ \frac{S_{n}}{S_{n+1}}Q^{(1)}_{n+1}(-s)-(X_{n,n}-s)Q^{(1)}_{n}(-s)-\frac{S_{n-1}}{S_{n}}Q^{(1)}_{n-1}(-s) \right]
		\begin{pmatrix}
			P_{n+1}(s) \\ P_{n}(s) \\ P_{n-1}(s)
		\end{pmatrix} ,
\label{s-Deriv-pi}
\end{multline}
\begin{multline}
	\partial_{t} 
		\begin{pmatrix}
			\pi_{n+1} \\ \pi_{n} \\ \pi_{n-1} 
		\end{pmatrix} = \mathcal{C}^{(\infty,0)}_{n}
		\begin{pmatrix}
			\pi_{n+1} \\ \pi_{n} \\ \pi_{n-1} 
		\end{pmatrix}
\\
	- \psi t^{b}e^{-t} \frac{1}{\eta_{n}}\left[ \frac{S_{n}}{S_{n+1}}Q_{n+1}(t)-(X_{n,n}+t)Q_{n}(t)-\frac{S_{n-1}}{S_{n}}Q_{n-1}(t) \right]
		\begin{pmatrix}
			P^{(1)}_{n+1}(-t) \\ P^{(1)}_{n}(-t) \\ P^{(1)}_{n-1}(-t)
		\end{pmatrix} ,
\label{t-Deriv-pi}
\end{multline}
\begin{multline}
	\partial_{s} 
		\begin{pmatrix}
			\eta_{n+1} \\ \eta_{n} \\ \eta_{n-1} 
		\end{pmatrix} = \mathcal{B}^{(\infty,\bar{0})}_{n}
		\begin{pmatrix}
			\eta_{n+1} \\ \eta_{n} \\ \eta_{n-1} 
		\end{pmatrix}
\\
		- \xi s^{a}e^{-s} \frac{1}{\pi_{n}}\left[ \frac{S_{n}}{S_{n+1}}P_{n+1}(s)-(Y_{n,n}+s)P_{n}(s)-\frac{S_{n-1}}{S_{n}}P_{n-1}(s) \right]
		\begin{pmatrix}
			Q^{(1)}_{n+1}(-s) \\ Q^{(1)}_{n}(-s) \\ Q^{(1)}_{n-1}(-s)
		\end{pmatrix} ,
\label{s-Deriv-eta}
\end{multline}
\begin{multline}
	\partial_{t} 
		\begin{pmatrix}
			\eta_{n+1} \\ \eta_{n} \\ \eta_{n-1} 
		\end{pmatrix} = \mathcal{C}^{(\infty,\bar{0})}_{n}
		\begin{pmatrix}
			\eta_{n+1} \\ \eta_{n} \\ \eta_{n-1} 
		\end{pmatrix}
\\
	- \psi t^{b}e^{-t} \frac{1}{\pi_{n}}\left[ \frac{S_{n}}{S_{n+1}}P^{(1)}_{n+1}(-t)-(Y_{n,n}-t)P^{(1)}_{n}(-t)-\frac{S_{n-1}}{S_{n}}P^{(1)}_{n-1}(-t) \right]
		\begin{pmatrix}
			Q_{n+1}(t) \\ Q_{n}(t) \\ Q_{n-1}(t)
		\end{pmatrix} .
\label{t-Deriv-eta}
\end{multline}
\end{proposition}
\begin{proof}
We differentiate the definitions \eqref{pi-eta}
and employ derivatives for the weights $w_1$, $w_2$ from \eqref{wCL_deform-diff:a},\eqref{wCL_deform-diff:b} 
and for the polynomials \eqref{Pdiff-s}-\eqref{Qdiff-t}.
Written in the simplest way using the intertwined polynomials we have
\begin{align*}
	\partial_{s}\pi_{n} & = \xi s^{a}e^{-s} P_{n}(s) \left[ 1-\tfrac{1}{2}\pi_{n}Q^{(1)}_{n}(-s)-\check{Q}^{(1)}_{n}(-s) \right] ,
\\
	\partial_{t}\pi_{n} & = \psi t^{b}e^{-t} P^{(1)}_{n}(-t) \left[ -\tfrac{1}{2}\pi_{n}Q_{n}(t)-\check{Q}_{n}(t) \right] ,
\\		
	\partial_{s}\eta_{n} & = \xi s^{a}e^{-s} Q^{(1)}_{n}(-s) \left[ -\tfrac{1}{2}\eta_{n}P_{n}(s)-\hat{P}_{n}(s) \right] ,
\\
	\partial_{t}\eta_{n} & = \psi t^{b}e^{-t} Q_{n}(t) \left[ 1-\tfrac{1}{2}\eta_{n}P^{(1)}_{n}(-t)-\hat{P}^{(1)}_{n}(-t) \right] .
\end{align*}
We next eliminate the intertwined polynomials using \eqref{hatPsubs},\eqref{checkQsubs},\eqref{hatP1subs} and \eqref{checkQ1subs}.
As in the case of the polynomial derivatives we construct the matrix derivatives by making the substitutions $ n\mapsto n\pm 1$ where necessary 
and use the recurrence relations \eqref{P_recur}, \eqref{Q_recur} to restore the indices $n\pm 2$ back into the range $\{n-1,n,n+1\}$.
\end{proof}

\begin{proposition}\label{XY_derivatives}
The deformation derivatives of $ X_{n,n} $ and $ Y_{n,n} $ are given by
\begin{align}
	\partial_{s}X_{n,n} & = \xi s^{a}e^{-s} 
	\left[ -\frac{S_{n}}{S_{n+1}}P_{n+1}(s)Q^{(1)}_{n}(-s) + \frac{S_{n-1}}{S_{n}}P_{n}(s)Q^{(1)}_{n-1}(-s) \right] ,
\\
	\partial_{t}X_{n,n} & = \psi t^{b}e^{-t} 
	\left[ -\frac{S_{n}}{S_{n+1}}P^{(1)}_{n+1}(-t)Q_{n}(t) + \frac{S_{n-1}}{S_{n}}P^{(1)}_{n}(-t)Q_{n-1}(t) \right] ,
\\
	\partial_{s}Y_{n,n} & = \xi s^{a}e^{-s} 
	\left[ -\frac{S_{n}}{S_{n+1}}P_{n}(s)Q^{(1)}_{n+1}(-s) + \frac{S_{n-1}}{S_{n}}P_{n-1}(s)Q^{(1)}_{n}(-s) \right] ,
\\
	\partial_{t}Y_{n,n} & = \psi t^{b}e^{-t} 
	\left[ -\frac{S_{n}}{S_{n+1}}P^{(1)}_{n}(-t)Q_{n+1}(t) + \frac{S_{n-1}}{S_{n}}P^{(1)}_{n-1}(-t)Q_{n}(t) \right] .
\end{align}
\end{proposition}
\begin{proof}
We describe only the derivation of $ X_{n,n} $ with respect to $s$ for brevity.
Differentiating $ X_{n,n} = \langle xP_{n},Q_{n} \rangle $ and substituting for the weight derivatives \eqref{wCL_deform-diff:a} 
and the polynomial derivatives \eqref{Pdiff-s}, \eqref{Qdiff-s} we are faced with integrals of the form $ \langle xP_{n},Q_{m} \rangle $.
However these can be evaluated by the formulae in \eqref{XY-elements} in terms of $ S_{n} $, $ \pi_{j}, \eta_{k} $ or $ X_{n,n} $,
and where sums are necessary these can be performed by \eqref{intertwineDN}, \eqref{intertwineUP}, or their analogues for the associated, intertwined functions.
Finally the intertwined variables can be eliminated using \eqref{hatPsubs},\eqref{checkQsubs},\eqref{hatP1subs} and \eqref{checkQ1subs}.
\end{proof}

\begin{remark}
The right-hand sides of Prop.\ref{XY_derivatives} are observed to be perfect differences in $n$ and therefore
a consequence of this is that the derivatives of the sub-leading coefficients \eqref{P_subleading} of the polynomials have the simple expressions 
\begin{align}
	\partial_{s}\frac{S_{n+1,n}}{S_{n+1}} & = \xi s^{a}e^{-s}\frac{S_{n}}{S_{n+1}}P_{n+1}(s)Q^{(1)}_{n}(-s) ,
\\
	\partial_{t}\frac{S_{n+1,n}}{S_{n+1}} & = \psi t^{b}e^{-t}\frac{S_{n}}{S_{n+1}}P^{(1)}_{n+1}(-t)Q_{n}(t) ,
\\
	\partial_{s}\frac{T_{n+1,n}}{S_{n+1}} & = \xi s^{a}e^{-s}\frac{S_{n}}{S_{n+1}}P_{n}(s)Q^{(1)}_{n+1}(-s) ,
\\
	\partial_{t}\frac{T_{n+1,n}}{S_{n+1}} & = \psi t^{b}e^{-t}\frac{S_{n}}{S_{n+1}}P^{(1)}_{n}(-t)Q_{n+1}(t) .	
\end{align}
\end{remark}

In preparation of further computations we require total derivatives of the anti-incidence $G$-matrices in a simple matrix form,
and this is furnished by the following result.
\begin{lemma}\label{anti-Gn-derivatives}
The total derivatives of the anti-incidence specialised matrices are given by 
\begin{multline}
	s\frac{d}{ds}G_{n}(s,-s) = (s-a)G_{n}(s,-s) + s\frac{d}{ds} \log(\pi_{n}\eta_{n})G_{n}(s,-s)
\\
	- \left[ \mathcal{A}^{(0,+)}_{n}+\mathcal{A}^{(\diagdown)}_{n} \right]^{T} G_{n}(s,-s)
	- G_{n}(s,-s)\left[ \mathcal{D}^{(0,-)}_{n}+\mathcal{A}^{(\diagdown)}_{n} \right]
\\
	- \psi t^{b+1}e^{-t}\frac{1}{\pi_{n}\eta_{n}}\frac{1}{s+t}	
	  \left(G_{n}(s,-s)\mathrm{Q}^{(0)}_{n}(t)\right) \otimes \left(G_{n}(-t,-s)^{T}\mathrm{P}^{(1)}_{n}(-t)\right)^{T}
\\ 
	+ \psi t^{b+1}e^{-t}\frac{1}{\pi_{n}\eta_{n}}\frac{1}{s+t}	
	  \left(G_{n}(s,t)\mathrm{Q}^{(0)}_{n}(t)\right) \otimes \left(G_{n}(s,-s)^{T}\mathrm{P}^{(1)}_{n}(-t)\right)^{T} ,
\label{dsGnPsMs}
\end{multline} 
and
\begin{multline}
	t\frac{d}{dt}G_{n}(-t,t) = (t-b)G_{n}(-t,t) + t\frac{d}{dt}\log(\pi_{n}\eta_{n})G_{n}(-t,t) 
\\
	- \left[ \mathcal{A}^{(0,-)}_{n}+\mathcal{D}^{(\diagdown)}_{n} \right]^{T} G_{n}(-t,t)
	- G_{n}(-t,t)\left[ \mathcal{D}^{(0,+)}_{n}+\mathcal{D}^{(\diagdown)}_{n} \right]
\\
	- \xi s^{a+1}e^{-s}\frac{1}{\pi_{n}\eta_{n}}\frac{1}{s+t}	
	  \left(G_{n}(-t,-s)\mathrm{Q}^{(1)}_{n}(-s)\right) \otimes \left(G_{n}(-t,t)^{T}\mathrm{P}^{(0)}_{n}(s)\right)^{T}
\\ 
	+ \xi s^{a+1}e^{-s}\frac{1}{\pi_{n}\eta_{n}}\frac{1}{s+t}	
	  \left(G_{n}(-t,t)\mathrm{Q}^{(1)}_{n}(-s)\right) \otimes \left(G_{n}(s,t)^{T}\mathrm{P}^{(0)}_{n}(s)\right)^{T} .
\label{dtGnMtPt}
\end{multline}
\end{lemma}

\begin{remark}
The $t$ derivative of \eqref{Bilinear_01} and the $s$ derivative of \eqref{Bilinear_10} both immediately vanish, 
after employing the derivatives from Prop.\ref{PQ_total-derivatives} and \ref{PiEta_derivatives} and using $\pi_{n}\eta_{n} = X_{n,n}+Y_{n,n}$.
However the $s$ derivative of \eqref{Bilinear_01} and the $t$ derivative of \eqref{Bilinear_10} can be shown to vanish but this requires the substitution of
\eqref{dsGnPsMs} into 
\begin{multline*}
	\mathrm{P}_{n}^{(0)}(s)^T
	\left\{ s\frac{d}{ds}G_{n}(s,-s)
	+ \left[ \mathcal{A}^{(0,+)}_{n}+\mathcal{A}^{(\diagdown)}_{n} \right]^{T} G_{n}(s,-s)
	+ G_{n}(s,-s)\left[ \mathcal{D}^{(0,-)}_{n}+\mathcal{A}^{(\diagdown)}_{n} \right] 
	\right.
\\
	+ \psi t^{b+1}e^{-t}\frac{1}{\pi_{n}\eta_{n}}\frac{1}{s+t}	
	  \left(G_{n}(s,-s)\mathrm{Q}^{(0)}_{n}(t)\right) \otimes \left(G_{n}(-t,-s)^{T}\mathrm{P}^{(1)}_{n}(-t)\right)^{T}
\\ 
	\left.
	- \psi t^{b+1}e^{-t}\frac{1}{\pi_{n}\eta_{n}}\frac{1}{s+t}	
	  \left(G_{n}(s,t)\mathrm{Q}^{(0)}_{n}(t)\right) \otimes \left(G_{n}(s,-s)^{T}\mathrm{P}^{(1)}_{n}(-t)\right)^{T}
	\right\} \mathrm{Q}_{n}^{(1)}(-s) ,
\end{multline*}
and \eqref{dtGnMtPt} into 
\begin{multline*}
	\mathrm{P}_{n}^{(1)}(-t)^T
	\left\{ t\frac{d}{dt}G_{n}(-t,t)
	+ \left[ \mathcal{A}^{(0,-)}_{n}+\mathcal{D}^{(\diagdown)}_{n} \right]^{T} G_{n}(-t,t)
	+ G_{n}(-t,t)\left[ \mathcal{D}^{(0,+)}_{n}+\mathcal{D}^{(\diagdown)}_{n} \right]
	\right.
\\
	+ \xi s^{a+1}e^{-s}\frac{1}{\pi_{n}\eta_{n}}\frac{1}{s+t}	
	  \left(G_{n}(-t,-s)\mathrm{Q}^{(1)}_{n}(-s)\right) \otimes \left(G_{n}(-t,t)^{T}\mathrm{P}^{(0)}_{n}(s)\right)^{T}
\\ 
	\left.
	- \xi s^{a+1}e^{-s}\frac{1}{\pi_{n}\eta_{n}}\frac{1}{s+t}	
	  \left(G_{n}(-t,t)\mathrm{Q}^{(1)}_{n}(-s)\right) \otimes \left(G_{n}(s,t)^{T}\mathrm{P}^{(0)}_{n}(s)\right)^{T}
	\right\} \mathrm{Q}_{n}^{(0)}(t) .
\end{multline*}
The $s$ derivative of \eqref{00_constraint} (both the first and second members) is just $ -\xi s^ae^{-s} P_{n}(s)Q^{(1)}_{n}(-s) $ times the expression itself, and similarly
the $t$ derivative of \eqref{00_constraint} is just $ -\psi t^be^{-t} P^{(1)}_{n}(-t)Q_{n}(t) $ times the same expression.
Both of the $s,t$ derivatives of \eqref{X_constraint} and \eqref{Y_constraint} are identically zero 
after employing the derivatives from Prop.\ref{PQ_total-derivatives} and \ref{PiEta_derivatives} and using \eqref{00_constraint}.
\end{remark}

\begin{remark}\label{XY+YXremark}
However the task of establishing the vanishing of the $s,t$ derivatives of \eqref{XY_constraint} and \eqref{YX_constraint} is not automatic and requires more effort.
One example, which is typical of the other cases, is the $s$ derivative of \eqref{YX_constraint}.
After dividing out the overall factor of $\xi s^{a}e^{-s}$ we split this quotient into three terms with individual factors of $\xi, \psi$ and unity.
The last of these three terms is then further separated into groups containing 
\begin{gather*}
	\frac{S_n}{S_{n+1}}\pi_{n+1} + \left(b+n-Y_{n,n}\right)\pi_{n}-\frac{S_{n-1}}{S_n}\pi_{n-1},
	\qquad
	\frac{S_n}{S_{n+1}}\eta_{n+1} + \left(a+n-X_{n,n}\right)\eta_{n}-\frac{S_{n-1} }{S_n}\eta_{n-1},
	\\
	\frac{S_n}{S_{n+1}}\pi_{n+1}\eta_{n}-\frac{S_{n-1}}{S_n}\pi_{n}\eta_{n-1},
	\qquad
	\frac{S_n}{S_{n+1}}\pi_{n}\eta_{n+1}-\frac{S_{n-1}}{S_n}\pi_{n-1}\eta_{n},	
\end{gather*}
and to only these specific groups the identities \eqref{YX_constraint}, \eqref{XY_constraint}, \eqref{X_constraint}, \eqref{Y_constraint} are applied, respectively.
In the final step one has to apply \eqref{00_constraint} and only then do all of the terms cancel.
In conclusion it should be remarked that no further identities are required to verify the consistency of the constraints with the dynamics.
\end{remark}

\subsection{Compatibility of the Spectral and Deformation Structures}
From the three-way compatibility of the first triple appearing in \eqref{xst-LaxPairs} we have three Schlesinger equations
\footnote{We do not examine here the three additional Schlesinger equations arising from the set \eqref{xst-LaxPairs} beyond what is reported in \S\ref{SpectralStructures}.}
\begin{equation}
	\partial_{s}\mathcal{A}_{n} - \partial_{x}\mathcal{B}_{n} = \left[\mathcal{B}_{n},\mathcal{A}_{n}\right] ,
\qquad
	\partial_{t}\mathcal{A}_{n} - \partial_{x}\mathcal{C}_{n} = \left[\mathcal{C}_{n},\mathcal{A}_{n}\right] ,
\qquad
	\partial_{t}\mathcal{B}_{n} - \partial_{s}\mathcal{C}_{n} = \left[\mathcal{C}_{n},\mathcal{B}_{n}\right] .
\end{equation}
Consequently we have the following Schlesinger equations for the $P$-residue matrices contained in the first set
\begin{gather}
	\mathcal{A}^{(s)}_{n}+\mathcal{B}^{(s)}_{n}  = \left[\mathcal{B}^{(s)}_{n},\mathcal{A}^{(s)}_{n}\right] ,
\label{AB_comp}\\
	\partial_{s}\mathcal{A}^{(0)}_{n} = \left[\mathcal{B}^{(\infty)}_{n},\mathcal{A}^{(0)}_{n}\right] - s^{-1}\left[\mathcal{B}^{(s)}_{n},\mathcal{A}^{(0)}_{n}\right] ,
\label{AB_comp_0s}\\
	\partial_{s}\mathcal{A}^{(s)}_{n} = \left[\mathcal{B}^{(\infty)}_{n},\mathcal{A}^{(s)}_{n}\right] + \left[\mathcal{B}^{(s)}_{n},\mathcal{A}^{(\infty)}_{n}\right]
	  + s^{-1}\left[\mathcal{B}^{(s)}_{n},\mathcal{A}^{(0)}_{n}\right] + (s+t)^{-1}\left[\mathcal{B}^{(s)}_{n},\mathcal{A}^{(-t)}_{n}\right] ,
\label{AB_comp_ss}\\
	\partial_{s}\mathcal{A}^{(-t)}_{n} = \left[\mathcal{B}^{(\infty)}_{n},\mathcal{A}^{(-t)}_{n}\right] - (s+t)^{-1}\left[\mathcal{B}^{(s)}_{n},\mathcal{A}^{(-t)}_{n}\right] ,
\label{AB_comp_-ts}\\
	\partial_{s}\mathcal{A}^{(\infty)}_{n} = \left[\mathcal{B}^{(\infty)}_{n},\mathcal{A}^{(\infty)}_{n}\right] ,
\label{AB_comp_inftys}
\end{gather}
the second set
\begin{gather}
	-\mathcal{A}^{(-t)}_{n}+\mathcal{C}^{(-t)}_{n}  = \left[\mathcal{C}^{(-t)}_{n},\mathcal{A}^{(-t)}_{n}\right] ,
\label{AC_comp}\\
	\partial_{t}\mathcal{A}^{(0)}_{n} = \left[\mathcal{C}^{(\infty)}_{n},\mathcal{A}^{(0)}_{n}\right] + t^{-1}\left[\mathcal{C}^{(-t)}_{n},\mathcal{A}^{(0)}_{n}\right] ,
\label{AC_comp_0t}\\
	\partial_{t}\mathcal{A}^{(s)}_{n} = \left[\mathcal{C}^{(\infty)}_{n},\mathcal{A}^{(s)}_{n}\right] + (s+t)^{-1}\left[\mathcal{C}^{(-t)}_{n},\mathcal{A}^{(s)}_{n}\right] ,
\label{AC_comp_tt}\\
	\partial_{t}\mathcal{A}^{(-t)}_{n} = \left[\mathcal{C}^{(\infty)}_{n},\mathcal{A}^{(-t)}_{n}\right] + \left[\mathcal{C}^{(-t)}_{n},\mathcal{A}^{(\infty)}_{n}\right]
	- t^{-1}\left[\mathcal{C}^{(-t)}_{n},\mathcal{A}^{(0)}_{n}\right] - (s+t)^{-1}\left[\mathcal{C}^{(-t)}_{n},\mathcal{A}^{(s)}_{n}\right] ,
\label{AC_comp_-tt}\\
	\partial_{t}\mathcal{A}^{(\infty)}_{n} = \left[\mathcal{C}^{(\infty)}_{n},\mathcal{A}^{(\infty)}_{n}\right] ,
\label{AC_comp_inftyt}
\end{gather}
along with a final set
\begin{gather}
	\partial_{t}\mathcal{B}^{(s)}_{n} = \left[\mathcal{C}^{(\infty)}_{n},\mathcal{B}^{(s)}_{n}\right] + (s+t)^{-1}\left[\mathcal{C}^{(-t)}_{n},\mathcal{B}^{(s)}_{n}\right] ,
\label{BC_comp_st}\\
	\partial_{s}\mathcal{C}^{(-t)}_{n} = \left[\mathcal{B}^{(\infty)}_{n},\mathcal{C}^{(-t)}_{n}\right] + (s+t)^{-1}\left[\mathcal{C}^{(-t)}_{n},\mathcal{B}^{(s)}_{n}\right] ,	
\label{BC_comp_-ts}\\
	\partial_{t}\mathcal{B}^{(\infty)}_{n}-\partial_{s}\mathcal{C}^{(\infty)}_{n} = \left[\mathcal{C}^{(\infty)}_{n},\mathcal{B}^{(\infty)}_{n}\right] .
\label{BC_comp_inftyt-inftys}
\end{gather}

\begin{proposition}
The compatibility conditions \eqref{AB_comp}--\eqref{BC_comp_inftyt-inftys} are identically satisfied given the dynamical equations 
\eqref{s-Deriv-pi}--\eqref{t-Deriv-eta}, \eqref{XY_derivatives}, \eqref{s-PartialDeriv-P0Ps}--\eqref{t-PartialDeriv-Q1Ms}, \eqref{s-TotalDeriv-P0Ps}--\eqref{s-TotalDeriv-Q1Ms}, 
and the constraints \eqref{Bilinear_01}, \eqref{Bilinear_10}, \eqref{00_constraint}, \eqref{X_constraint}--\eqref{YX_constraint}. 
\end{proposition}
\begin{proof}
From the evaluations immediately preceding \eqref{BresidueSing@Infty} and \eqref{CresidueSing@Infty} both \eqref{AB_comp} and \eqref{AC_comp} are trivially satisfied, 
while \eqref{BC_comp_st} is equivalent to \eqref{AC_comp_tt} and \eqref{BC_comp_-ts} is equivalent to \eqref{AB_comp_-ts}.
Instead of examining \eqref{AB_comp_0s} by itself it is simpler to treat $s\partial_{s}\mathcal{A}^{(\Sigma)}_{n}$, 
i.e. the sum of \eqref{AB_comp_0s}, \eqref{AB_comp_ss} and \eqref{AB_comp_-ts}.
We find that all entries of this matrix relation are identically zero using the full set of derivatives except for the $(+1,0)$ and $(-1,0)$ entries.
We detail the steps only for the $(+1,0)$ entry as the case for other entry is handled in a similar manner.
For this entry we have to employ the approach described in Remark \ref{XY+YXremark}, and in particular isolate the groups with following terms
\begin{equation*}
	\frac{S_n}{S_{n+1}}\eta_{n+1} + \left(a+n-X_{n,n}\right)\eta_{n}-\frac{S_{n-1}}{S_n}\eta_{n-1},
\qquad
	\frac{S_n}{S_{n+1}}\left[ \pi_{n+1}\eta_{n}-\pi_{n}\eta_{n+1} \right],
\end{equation*}
and apply \eqref{XY_constraint} to the former and the linear combination $\pi_{n}\eta_{n}\times$\eqref{XY_constraint} plus \eqref{X_constraint} to the latter, respectively.
Then everything mutually cancels.
For the $(-1,0)$ entry one requires the combination $\pi_{n}\eta_{n}\times$\eqref{YX_constraint} plus \eqref{X_constraint} to effect the cancellation.
In order to verify \eqref{AB_comp_ss} by itself we first compute the $s$-derivative of the residue matrix as
\begin{multline*}
	s\frac{d}{ds}\mathcal{A}^{(s)}_{n} = 
	\left[ \mathcal{A}^{(0,+)}_{n}+\mathcal{A}^{(\diagdown)}_{n},\mathcal{A}^{(s)}_{n} \right]
\\
	+ \xi\psi s^{a}e^{-s}t^{b+1}e^{-t}\frac{1}{\pi_{n}\eta_{n}}	K^{0,0}_{n}(s,t)
	  \mathrm{P}^{(1)}_{n}(-t) \otimes \left(G_{n}(s,-s)\mathrm{Q}^{(1)}_{n}(-s)\right)^{T}
\\
	+ \xi\psi s^{a}e^{-s}t^{b+1}e^{-t}\frac{1}{\pi_{n}\eta_{n}}	K^{1,1}_{n}(-t,-s)
	  \mathrm{P}^{(0)}_{n}(s) \otimes \left(G_{n}(s,t)\mathrm{Q}^{(0)}_{n}(t)\right)^{T}
\\ 
	+ \xi\psi s^{a}e^{-s}t^{b+1}e^{-t}\frac{P^{(1)}_{n}(-t)}{\pi_{n}^2\eta_{n}^2}
	  \left[ \frac{S_n}{S_{n+1}}Q^{(1)}_{n+1}(-s) - \left(X_{n,n}-s\right)Q^{(1)}_{n}(-s)-\frac{S_{n-1}}{S_n}Q^{(1)}_{n-1}(-s) \right]	
	  \mathrm{P}^{(0)}_{n}(s) \otimes \left(G_{n}(s,t)\mathrm{Q}^{(0)}_{n}(t)\right)^{T} ,
\end{multline*} 
using the result \eqref{dsGnPsMs} from Lemma \ref{anti-Gn-derivatives} and the derivatives \eqref{s-TotalDeriv-P0Ps}, \eqref{s-TotalDeriv-Q1Ms}.
Substituting this into \eqref{AB_comp_ss} leads to the immediate cancellation of the $(+1,+1),(+1,-1),(0,0),(-1,+1),(-1,-1)$ entries 
while the remainder require further application of \eqref{YX_constraint} to eliminate them.
Lastly \eqref{AB_comp_-ts} and \eqref{AB_comp_inftys} are satisfied identically with the appropriate set of derivatives.
The full set of conditions \eqref{AC_comp}, \eqref{AC_comp_0t}, \eqref{AC_comp_tt}, \eqref{AC_comp_-tt}, \eqref{AC_comp_inftyt} follow in a similar manner,
and the cross conditions \eqref{BC_comp_st}, \eqref{BC_comp_-ts}, and \eqref{BC_comp_inftyt-inftys} are satisfied without any additional use of identities.
\end{proof}

\section{Conclusions and Some Observations}\label{}
\setcounter{equation}{0}

The gap probability of the Bures-Hall ensemble is a refined local statistic of the density matrix spectrum of a truly generic quantum system, especially of the extreme eigenvalues, 
which has hitherto not received sufficient attention.  
Following the advances made here, building upon the initial foundations by Bertola et al, 
it is now accessible with the tools of "integrable probability" and most significantly rigorous asymptotic analysis of the distribution of the lowest eigenvalue via Riemann-Hilbert methods.
In addition the Bures-Hall ensembles, in their fixed trace or unconstrained forms, are rare examples of an integrable formulation of a {\it Pfaffian point processes}, 
whereas virtually all the other known cases are of {\it determinantal point processes}.
From the integrable systems point of view the deformed Cauchy-Laguerre bi-orthogonal polynomial system is an irreducibly rank 3 system and clearly beyond the rank 2 Painlev{\'e} class,
which have arisen in the vast majority of random matrix ensembles with integrable structure.

\bibliographystyle{plain}

\begin{thebibliography}{10}

\bibitem{Ami_1979-80}
S.~A. Amitsur.
\newblock On the characteristic polynomial of a sum of matrices.
\newblock {\em Linear and Multilinear Algebra}, 8(3):177--182, 1979/80.

\bibitem{AZ_2015}
M.~R. Atkin and S.~Zohren.
\newblock Violations of {B}ell inequalities from random pure states.
\newblock {\em Phys. Rev. A (3)}, 92(1):012331, 6, 2015.

\bibitem{Aub_2014}
G.~Aubrun.
\newblock Is a random state entangled?
\newblock In {\em X{VII}th {I}nternational {C}ongress on {M}athematical {P}hysics}, pages 534--541. World Sci. Publ., Hackensack, NJ, 2014.

\bibitem{Aubrun+Szarek_2017}
G.~Aubrun and S.~J. Szarek.
\newblock {\em Alice and {B}ob meet {B}anach}, volume 223 of {\em Mathematical Surveys and Monographs}.
\newblock American Mathematical Society, Providence, RI, 2017.
\newblock The interface of asymptotic geometric analysis and quantum information theory.

\bibitem{ASY_2014}
G.~Aubrun, S.~J. Szarek, and D.~Ye.
\newblock Entanglement thresholds for random induced states.
\newblock {\em Comm. Pure Appl. Math.}, 67(1):129--171, 2014.

\bibitem{Bau_1990}
W.~C. Bauldry.
\newblock Estimates of asymmetric {F}reud polynomials on the real line.
\newblock {\em J. Approx. Theory}, 63(2):225–237, 1990.

\bibitem{Bengtsson+Zyczkowski_2017}
I.~Bengtsson and K.~\v{Z}yczkowski.
\newblock {\em Geometry of {Q}uantum {S}tates}.
\newblock Cambridge University Press, Cambridge, 2017.
\newblock An introduction to quantum entanglement, Second edition.

\bibitem{BGS_2009}
M.~Bertola, M.~Gekhtman, and J.~Szmigielski.
\newblock The {C}auchy two-matrix model.
\newblock {\em Comm. Math. Phys.}, 287(3):983--1014, 2009.

\bibitem{BGS_2010}
M.~Bertola, M.~Gekhtman, and J.~Szmigielski.
\newblock Cauchy biorthogonal polynomials.
\newblock {\em J. Approx. Theory}, 162(4):832--867, 2010.

\bibitem{BGS_2013}
M.~Bertola, M.~Gekhtman, and J.~Szmigielski.
\newblock Strong asymptotics for {C}auchy biorthogonal polynomials with application to the {C}auchy two-matrix model.
\newblock {\em J. Math. Phys.}, 54(4):043517, 25, 2013.

\bibitem{BGS_2014}
M.~Bertola, M.~Gekhtman, and J.~Szmigielski.
\newblock Cauchy-{L}aguerre two-matrix model and the {M}eijer-{G} random point field.
\newblock {\em Comm. Math. Phys.}, 326(1):111--144, 2014.

\bibitem{BC_1990}
S.~S. Bonan and D.~S. Clark.
\newblock Estimates of the {H}ermite and the {F}reud polynomials.
\newblock {\em J. Approx. Theory}, 63(2):210–224, 1990.

\bibitem{Bor_2010b}
F.~Bornemann.
\newblock On the numerical evaluation of {F}redholm determinants.
\newblock {\em Math. Comp.}, 79(270):871--915, 2010.

\bibitem{BN_2011}
G.~Borot and C.~Nadal.
\newblock Purity distribution for generalized random bures mixed states.
\newblock {\em J. Phys. A}, 45:7,(7):075209, 43, 2012.

\bibitem{Bur_1969}
D.~Bures.
\newblock An extension of {K}akutani's theorem on infinite product measures to the tensor product of semifinite {$w^{\ast} $}-algebras.
\newblock {\em Trans. Amer. Math. Soc.}, 135:199--212, 1969.

\bibitem{CLZ_2010}
Y.~Chen, D.-Z. Liu, and D.-S. Zhou.
\newblock Smallest eigenvalue distribution of the fixed-trace {L}aguerre
  beta-ensemble.
\newblock {\em J. Phys. A}, 43(31):315303, 13, 2010.

\bibitem{CN_2016}
B.~Collins and I.~Nechita.
\newblock Random matrix techniques in quantum information theory.
\newblock {\em J. Math. Phys.}, 57(1):015215, 34, 2016.

\bibitem{Dit_1993}
J.~Dittmann.
\newblock On the {R}iemannian geometry of finite-dimensional mixed states.
\newblock {\em Sem. Sophus Lie}, 3(1):73--87, 1993.

\bibitem{FK_2016}
P.~J. Forrester and M.~Kieburg.
\newblock Relating the {B}ures measure to the {C}auchy two-matrix model.
\newblock {\em Comm. Math. Phys.}, 342(1):151--187, October 2016.

\bibitem{FL_2019}
P.~J. Forrester and S.-H. Li.
\newblock {F}ox {H}-kernel and $\theta$-deformation of the {C}auchy {T}wo-{M}atrix {M}odel and {B}ures {E}nsemble.
\newblock {\em Int. Math. Res. Not.}, (8):5791--5824, 02 2019.

\bibitem{Fub_1904}
G.~Fubini.
\newblock Sulle metriche definite da una forma hermitiana. {Nota}.
\newblock {\em Ven. {Ist}. {Atti}}, 63((8) 6):501--513, 1904.

\bibitem{Gir_2007}
O.~Giraud.
\newblock Distribution of bipartite entanglement for random pure states.
\newblock {\em J. Phys. A - Math. Theoret.}, 40(11):2793--2801, 2007.

\bibitem{Hal_1998}
M.~J.~W. Hall.
\newblock Random quantum correlations and density operator distributions.
\newblock {\em Phys. Lett. A}, 242(3):123--129, 1998.

\bibitem{HLW_2006}
P.~Hayden, D.~W. Leung, and A.~Winter.
\newblock Aspects of generic entanglement.
\newblock {\em Comm. Math. Phys.}, 265(1):95--117, 2006.

\bibitem{Hueb_1993}
M.~H\"{u}bner.
\newblock The {B}ures metric and {U}hlmann's transition probability: explicit results.
\newblock In {\em Classical and quantum systems ({G}oslar, 1991)}, pages 406--409. World Sci. Publ., River Edge, NJ, 1993.

\bibitem{Hueb_1993a}
M.~H\"{u}bner.
\newblock Computation of {U}hlmann's parallel transport for density matrices and the {B}ures metric on three-dimensional {H}ilbert space.
\newblock {\em Phys. Lett. A}, 179(4-5):226--230, 1993.

\bibitem{KZM_2002}
V.M. Kendon, K.~Zyczkowski, and W.J. Munro.
\newblock Bounds on entanglement in qudit subsystems.
\newblock {\em Phys. Rev. A}, 66(6), 2002.

\bibitem{LAK_2021}
A.~Laha, A.~Aggarwal, and S.~Kumar.
\newblock Random density matrices: analytical results for mean root fidelity and mean square bures distance.
\newblock arXiv:2105.02743, 2021.

\bibitem{LiLi_2019}
C.~Li and S.-H. Li.
\newblock The {C}auchy two-matrix model, {C}-{T}oda lattice and {CKP} hierarchy.
\newblock {\em J. Nonlinear Sci.}, 29(1):3--27, 2019.

\bibitem{LW_2021a}
S.-H. Li and L.~Wei.
\newblock Moments of quantum purity and biorthogonal polynomial recurrence.
\newblock {\em J. Phys. A. Math. Theor.}, 54(44):Paper No. 445204, 29, 2021.

\bibitem{LZ_2010}
D.-Z. Liu and D.-S. Zhou.
\newblock {Local Statistical Properties of {S}chmidt Eigenvalues of Bipartite Entanglement for a Random Pure State}.
\newblock {\em Int. Math. Res. Not. IMRN}, 2011(4):725--766, 05 2011.

\bibitem{NMV_2011}
C.~Nadal, S.~N. Majumdar, and M.~Vergassola.
\newblock Statistical distribution of quantum entanglement for a random bipartite state.
\newblock {\em J. Stat. Phys.}, 142(2):403--438, 2011.

\bibitem{NDGD_2006}
M.~A. Nielsen, M.~R. Dowling, M.~Gu, and A.~C. Doherty.
\newblock Quantum computation as geometry.
\newblock {\em Science}, 311(5764):1133--1135, 2006.

\bibitem{OSZ_2010}
V.~Al. Osipov, H.-J. Sommers, and K.~\v{Z}yczkowski.
\newblock Random {B}ures mixed states and the distribution of their purity.
\newblock {\em Journal of Physics A: Mathematical and Theoretical}, 43(5):055302, 2010.

\bibitem{Pag_1993}
D.~N. Page.
\newblock Average entropy of a subsystem.
\newblock {\em Phys. Rev. Lett.}, 71(9):1291--1294, 1993.

\bibitem{PS_1996}
D.~Petz and C.~Sud\'{a}r.
\newblock Geometries of quantum states.
\newblock {\em J. Math. Phys.}, 37(6):2662--2673, 1996.

\bibitem{DLMF}
{NIST} {D}igital {L}ibrary of {M}athematical {F}unctions.
\newblock http://dlmf.nist.gov/, Release 1.0.9 of 2014-08-29.

\bibitem{RS_1987}
C.~Reutenauer and M.-P. Sch\"{u}tzenberger.
\newblock A formula for the determinant of a sum of matrices.
\newblock {\em Lett. Math. Phys.}, 13(4):299--302, 1987.

\bibitem{SK_2019}
A.~Sarkar and S.~Kumar.
\newblock Bures-{H}all ensemble spectral densities and average entropies.
\newblock {\em J. Phys. A}, 52(29):295203, 25, 2019.

\bibitem{SZ_2003}
H.-J. Sommers and K.~\v{Z}yczkowski.
\newblock Bures volume of the set of mixed quantum states.
\newblock {\em J. Phys. A}, 36(39):10083--10100, 2003.

\bibitem{SZ_2004}
H.-J. Sommers and K.~\v{Z}yczkowski.
\newblock Statistical properties of random density matrices.
\newblock {\em J. Phys. A}, 37(35):8457--8466, 2004.

\bibitem{Stu_1905}
E.~Study.
\newblock K{\"u}rzeste wege im komplexen gebiet.
\newblock {\em Mathematische Annalen}, 60(3):321--378, 1905.

\bibitem{Sza_2005}
S.~J. Szarek.
\newblock Volume of separable states is super-doubly-exponentially small in the number of qubits.
\newblock {\em Phys. Rev. A}, 72:032304, 2005.

\bibitem{SWZ_2011}
S.~J. Szarek, E.~Werner, and K.~\.{Z}yczkowski.
\newblock How often is a random quantum state {$k$}-entangled?
\newblock {\em J. Phys. A}, 44(4):045303, 15, 2011.

\bibitem{Uhl_1976}
A.~Uhlmann.
\newblock The ``transition probability'' in the state space of a {$\sp*$}-algebra.
\newblock {\em Rep. Mathematical Phys.}, 9(2):273--279, 1976.

\bibitem{Uhl_1986}
A.~Uhlmann.
\newblock Parallel transport and ``quantum holonomy'' along density operators.
\newblock {\em Rep. Math. Phys.}, 24(2):229--240, 1986.

\bibitem{Uhl_1992}
A.~Uhlmann.
\newblock The metric of {B}ures and the geometric phase.
\newblock In {\em Quantum groups and related topics ({W}roc{\l}aw, 1991)}, volume~13 of {\em Math. Phys. Stud.}, pages 267--274. Kluwer Acad. Publ., Dordrecht, 1992.

\bibitem{ZS_2001}
K.~\v{Z}yczkowski and H.-J. Sommers.
\newblock Induced measures in the space of mixed quantum states.
\newblock volume~34, pages 7111--7125. 2001.
\newblock Quantum information and computation.

\bibitem{ZS_2003}
K.~\v{Z}yczkowski and H.-J. Sommers.
\newblock Hilbert-{S}chmidt volume of the set of mixed quantum states.
\newblock {\em J. Phys. A}, 36(39):10115--10130, 2003.

\bibitem{ZS_2005b}
K.~\v{Z}yczkowski and H.-J. Sommers.
\newblock Average fidelity between random quantum states.
\newblock {\em Phys. Rev. A (3)}, 71(3, part A):032313, 11, 2005.

\bibitem{Wei_2020a}
L.~Wei.
\newblock Exact variance of von {N}eumann entanglement entropy over the {B}ures-{H}all measure.
\newblock {\em Phys. Rev. E}, 102(6):062128, 11, 2020.

\bibitem{Wei_2020}
L.~Wei.
\newblock Proof of {S}arkar{\textendash}{K}umar's conjectures on average entanglement entropies over the {B}ures{\textendash}{H}all ensemble.
\newblock {\em J. Phys. A}, 53(23):235203, 2020.

\bibitem{WW_2021}
L.~Wei and N.~S. Witte.
\newblock Quantum interpolating ensemble: Average entropies and orthogonal polynomials.
\newblock arXiv:2103.04231 2021.

\bibitem{Ye_2009}
D.~Ye.
\newblock On the {B}ures volume of separable quantum states.
\newblock {\em J. Math. Phys.}, 50(8):083502, 14, 2009.

\bibitem{Ye_2010}
D.~Ye.
\newblock On the comparison of volumes of quantum states.
\newblock {\em J. Phys. A}, 43(31):315301, 17, 2010.

\bibitem{ZPNC_2011}
K.~\v{Z}yczkowski, K.~A. Penson, I.~Nechita, and B.~Collins.
\newblock Generating random density matrices.
\newblock {\em J. Math. Phys.}, 52(6):062201, 20, 2011.

\end{thebibliography}
\def\cprime{$'$} \def\cprime{$'$} \def\cprime{$'$} \def\cprime{$'$}

\end{document}